\documentclass[10pt,a4paper,reqno]{amsart}
\usepackage[T1]{fontenc}
\usepackage[utf8]{inputenc}
\usepackage[english]{babel}
\usepackage{lmodern,amsmath,amssymb,amsthm,mathtools,mathrsfs}
\usepackage[a4paper,left=28mm,right=28mm,top=30mm,bottom=30mm,headheight=14pt,headsep=7mm,footskip=10mm]{geometry}
\usepackage{microtype,booktabs,array,longtable,enumitem,float}
\usepackage[dvipsnames]{xcolor}
\usepackage{tikz-cd}
\usepackage{hyperref}
\usepackage{bookmark}
\hypersetup{colorlinks=true,linkcolor=blue,citecolor=green,urlcolor=magenta,linktocpage=true,pdftitle={Finite abelian groups acting on rationally connected threefolds III: Classification},pdfauthor={Konstantin Loginov},pdfsubject={Finite abelian actions, terminal Fano threefolds, anticanonical linear systems}}
\allowdisplaybreaks[2]
\numberwithin{equation}{section}
\newtheorem{theorem}{Theorem}[section]
\newtheorem{proposition}[theorem]{Proposition}
\newtheorem{lemma}[theorem]{Lemma}
\newtheorem{corollary}[theorem]{Corollary}
\newtheorem{question}[theorem]{Question}
\theoremstyle{definition}
\newtheorem{definition}[theorem]{Definition}
\newtheorem{example}[theorem]{Example}
\newtheorem{assumption}[theorem]{Assumption}
\newtheorem{remark}[theorem]{Remark}
\newcommand{\CC}{\mathbb C}
\newcommand{\QQ}{\mathbb Q}
\newcommand{\ZZ}{\mathbb Z}

\newcommand{\PP}{\mathbb P}
\renewcommand{\AA}{\mathbb A}
\newcommand{\OO}{\mathscr O}
\newcommand{\MM}{\mathscr M}
\newcommand{\TT}{\mathscr T}
\newcommand{\B}{\mathbf B}
\newcommand{\plur}{\boldsymbol{P}}
\newcommand{\eps}{\varepsilon}
\newcommand{\tors}{\mathrm{tors}}
\newcommand{\reg}{\mathrm{reg}}
\newcommand{\nG}{\mathrm{nG}}
\newcommand{\id}{\mathrm{id}}
\DeclareMathOperator{\Aut}{Aut}
\DeclareMathOperator{\Bir}{Bir}
\DeclareMathOperator{\Cr}{Cr}
\DeclareMathOperator{\Cl}{Cl}
\DeclareMathOperator{\Pic}{Pic}
\DeclareMathOperator{\NS}{NS}
\DeclareMathOperator{\rk}{rk}
\DeclareMathOperator{\ct}{ct}
\DeclareMathOperator{\ord}{ord}
\DeclareMathOperator{\Fix}{Fix}
\DeclareMathOperator{\Sing}{Sing}
\DeclareMathOperator{\Supp}{Supp}
\DeclareMathOperator{\Bs}{Bs}

\DeclareMathOperator{\Alb}{Alb}
\DeclareMathOperator{\MW}{MW}

\DeclareMathOperator{\Spec}{Spec}
\DeclareMathOperator{\Proj}{Proj}
\DeclareMathOperator{\Tr}{Tr}
\DeclareMathOperator{\diag}{diag}

\newcommand{\arxiv}[1]{\href{https://arxiv.org/abs/#1}{\texttt{arXiv:#1}}}
\newlength{\grouplistleftwidth}
\newlength{\grouplistrightwidth}
\newlength{\kthreeleftwidth}
\newlength{\kthreerightwidth}
\AtBeginDocument{%
 \settowidth{\grouplistleftwidth}{$\displaystyle\Gamma_3=(\ZZ/6\ZZ)^2\times(\ZZ/3\ZZ)^2,$}%
 \settowidth{\grouplistrightwidth}{$\displaystyle\Gamma_4=(\ZZ/8\ZZ)^2\times\ZZ/4\ZZ\times\ZZ/2\ZZ.$}%
 \settowidth{\kthreeleftwidth}{$\displaystyle H_3=\ZZ/6\ZZ\times(\ZZ/3\ZZ)^2,$}%
 \settowidth{\kthreerightwidth}{$\displaystyle H_4=\ZZ/8\ZZ\times\ZZ/4\ZZ\times\ZZ/2\ZZ.$}%
}
\setlist[enumerate]{label=\textup{(\roman*)},leftmargin=2em,itemsep=2pt}
\title[Finite abelian groups on threefolds III]{Finite abelian groups acting on rationally connected\\ threefolds~III: Classification}
\author{Konstantin Loginov}
\address{Steklov Mathematical Institute, Moscow, Russia}
\email{loginov@mi-ras.ru}
\date{}
\begin{document}
\begin{abstract}
We complete the classification of finite abelian groups that admit a faithful action on a rationally connected threefold. For a finite abelian group $G$ with invariant factors $1<d_1\mid\cdots\mid d_s$, put $c(G)=\sum_i(1-1/d_i)$. We prove that for $1\le n\le3$ such a group admits a faithful action on a rationally connected variety of dimension $n$ if and only if $c(G)\le n$. We also show that finite abelian subgroups of the Cremona group of rank three are, up to isomorphism, precisely the products of finite abelian subgroups of the Cremona groups of ranks one and two. Four further groups can act faithfully on rationally connected threefolds. For each of these groups, the faithful action is unique up to birational conjugacy and automorphism of the acting group.
\end{abstract}
\maketitle
\tableofcontents

\section{Introduction}\label{sec:intro}
We work over $\CC$. We complete the classification of finite abelian groups admitting faithful action on rationally connected threefolds, building on \cite{Log,LPZ,Lo26}. The existence of these actions has a numerical characterization in terms of the following invariant.

\begin{definition}\label{def:complexity}
For a nontrivial finite abelian group $G$ with invariant-factor decomposition
\begin{equation}\label{eq:groupdecomposition}
 G\simeq\ZZ/d_1\ZZ\times\cdots\times\ZZ/d_s\ZZ,
 \qquad1<d_1\mid\cdots\mid d_s,
\end{equation}
define its \emph{complexity} by
\begin{equation}\label{eq:complexity}
 c(G)=\sum_{i=1}^s\left(1-\frac1{d_i}\right).
\end{equation}
For the trivial group, set $c(G)=0$.
\end{definition}

The sum defining $c(G)$ occurs in minimum-genus formulae for actions of finite abelian groups on curves \cite[Theorem~4]{Maclachlan}. It also appears in \cite[Conjecture~1.2]{Girard}. Our main result is as follows.

\begin{theorem}\label{thm:complexity}
Let $G$ be a finite abelian group and let $n$ be an integer with $1\le n\le3$. Then $G$ admits a faithful action on a rationally connected variety of dimension $n$ if and only if $c(G)\le n$.
\end{theorem}

The cases $n=1,2$ follow from the known classifications recalled below. The new result is the case $n=3$. 
For a variety $X$, denote its group of birational automorphisms by $\Bir(X)$, and write $\Cr_n(\CC)=\Bir(\PP^n)$. A finite group embeds into $\Cr_n(\CC)$ if and only if it admits a faithful action on a rational variety of dimension $n$. Starting in dimension three, the minimal model program for varieties with a finite group action leads naturally to the broader class of rationally connected varieties.

In dimension one, we have $\Cr_1(\CC)=\operatorname{PGL}_2(\CC)$. A finite abelian group embeds into $\Cr_1(\CC)$ if and only if it is isomorphic to $\ZZ/n\ZZ$ for some positive integer $n$ or to $(\ZZ/2\ZZ)^2$. For finite abelian subgroups of $\Cr_1(\CC)$, the maximal complexity is one and is attained only by the Klein four group $(\ZZ/2\ZZ)^2$. 
In dimension two, Blanc determined the isomorphism classes of finite abelian subgroups of $\Cr_2(\CC)$.

\begin{theorem}[{\cite[Theorem~4.6]{Blanc}}]\label{thm:blanc}
A finite abelian group embeds into $\Cr_2(\CC)$ if and only if it is isomorphic to one of the groups
\begin{equation}\label{eq:A12}
 \ZZ/m\ZZ\times\ZZ/n\ZZ,\quad
 \ZZ/2m\ZZ\times(\ZZ/2\ZZ)^2,\quad
 (\ZZ/4\ZZ)^2\times\ZZ/2\ZZ,\quad
 (\ZZ/3\ZZ)^3,\quad (\ZZ/2\ZZ)^4,
\end{equation}
where $m,n$ are positive integers.
\end{theorem}

Another proof of Theorem~\ref{thm:blanc} is given in \cite[\S6]{Log}. 
For $\Cr_2(\CC)$, the maximal complexity of its finite abelian subgroup is two and is attained precisely by the last three groups in \eqref{eq:A12}. 

\begin{theorem}\label{thm:cremona}
A finite abelian group $G$ embeds into $\Cr_3(\CC)$ if and only if it is {of product type}, that is, $G\simeq G_1\times G_2$, where $G_i\subset \Cr_i(\CC)$ for $i=1,2$. Up to isomorphism, the groups of product type are precisely
\begin{equation}\label{eq:productlist}
\begin{alignedat}{2}
 &\mathmakebox[\grouplistleftwidth][l]{\ZZ/k\ZZ\times\ZZ/l\ZZ\times\ZZ/m\ZZ,}
 &\qquad &\mathmakebox[\grouplistrightwidth][l]{\ZZ/2k\ZZ\times(\ZZ/4\ZZ)^2\times\ZZ/2\ZZ,}\\[2pt]
 &\ZZ/3k\ZZ\times(\ZZ/3\ZZ)^3,
 &&\ZZ/2k\ZZ\times\ZZ/2l\ZZ\times(\ZZ/2\ZZ)^2,\\[2pt]
 &\ZZ/2k\ZZ\times(\ZZ/2\ZZ)^4,
 &&(\ZZ/4\ZZ)^2\times(\ZZ/2\ZZ)^3,\\[2pt]
 &(\ZZ/2\ZZ)^6,
 &&
\end{alignedat}
\end{equation}
where $k,l,m$ are positive integers.
\end{theorem}

Every group of product type embeds into $\Cr_3(\CC)$ by acting on the product of a rational curve and a rational surface. Already in dimension four, the group $(\ZZ/4\ZZ)^5$ embeds into $\Cr_4(\CC)$ and admits no decomposition as a product of subgroups of Cremona groups of smaller positive ranks summing to four, see Example~\ref{ex:rank-five}. 
For finite abelian subgroups of $\Cr_3(\CC)$, the maximal complexity is three and is attained precisely by the last two groups in \eqref{eq:productlist}. 
Put
\begin{equation}\label{eq:exceptional-groups}
\begin{alignedat}{2}
 &\mathmakebox[\grouplistleftwidth][l]{\Gamma_1=(\ZZ/4\ZZ)^4,}
 &\qquad &\mathmakebox[\grouplistrightwidth][l]{\Gamma_2=(\ZZ/6\ZZ)^3\times\ZZ/2\ZZ,}\\[2pt]
 &\Gamma_3=(\ZZ/6\ZZ)^2\times(\ZZ/3\ZZ)^2,
 &&\Gamma_4=(\ZZ/8\ZZ)^2\times\ZZ/4\ZZ\times\ZZ/2\ZZ.
\end{alignedat}
\end{equation}

Each group $\Gamma_i$ has complexity three and acts faithfully on a rationally connected Fermat hypersurface $X_i$ by the standard diagonal action described in Example~\ref{ex:exceptional-actions}. 
By Theorem~\ref{thm:complexity} and Lemma~\ref{lem:complexitytypes}, the maximal complexity of a finite abelian group admitting a faithful action on a rationally connected threefold is $3$. Up to isomorphism, the groups attaining the maximum complexity are $\Gamma_1,\ldots,\Gamma_4$ and the last two groups in \eqref{eq:productlist}.
Theorem~\ref{thm:cremona} follows from the next result, which also determines the faithful actions of $\Gamma_1,\ldots,\Gamma_4$ on rationally connected threefolds up to birational conjugacy and automorphisms of the acting group. The case $G\simeq\Gamma_1$ was proved in \cite[Theorem~1.3]{Lo26}.

\begin{theorem}\label{thm:rc-models}
Let a finite abelian group $G$ act faithfully on a rationally connected threefold $X$. Then exactly one of the following holds:
\begin{enumerate}[label=\textup{(\arabic*)}]
\item $G$ is of product type.
\item $G\simeq\Gamma_i$ for some $i\in\{1,2,3,4\}$. Let $X_i$ be the hypersurface in Example~\ref{ex:exceptional-actions} with the diagonal $\Gamma_i$-action given there. There exist an isomorphism $\alpha\colon G\xrightarrow{\sim}\Gamma_i$ and a birational map $\varphi\colon X\dashrightarrow X_i$ satisfying
\[
 \varphi\circ g=\alpha(g)\circ\varphi\qquad\text{for every }g\in G.
\]
The threefold $X$ is nonrational.
\end{enumerate}
\end{theorem}

By \cite[Theorem~1.7]{Log} and \cite[Theorem~1.5]{LPZ}, if a finite abelian group $G$ acts faithfully on a rationally connected threefold $V$, then either $G$ is of product type, or $G\simeq\Gamma_i$ for some $i\in\{1,2,3,4\}$, or $V$ is $G$-birational to a $G\QQ$-Fano threefold with empty anticanonical linear system. For $G$ not of product type, the following theorem produces a $G\QQ$-Fano threefold with nonempty anticanonical linear system in the $G$-birational class of $V$, proving \cite[Conjecture~1.8]{Log}.

\begin{theorem}\label{thm:nonemptymodel}
Let $G$ be a finite abelian group acting faithfully on a rationally connected threefold $V$. If $G$ is not of product type, then $V$ is $G$-birational to a $G\QQ$-Fano threefold $X$ with nonempty anticanonical linear system.
\end{theorem}

The sufficient condition in Theorem~\ref{thm:complexity} holds in every dimension and admits the following characterization by finite covers of projective space from klt Fano varieties.

\begin{proposition}\label{prop:complexityrealization}
Let $G$ be a finite abelian group and let $n\ge1$ be an integer. Then $c(G)\le n$ if and only if there exists a klt Fano variety $X$ of dimension $n$ with a faithful action of $G$ such that
\[
 X/G\simeq\PP^n.
\]
In particular, $c(G)\le n$ is sufficient for the existence of a faithful action of $G$ on a rationally connected variety of dimension~$n$.
\end{proposition}

The complexity $c(G)$ of a finite abelian group $G$ also has an interpretation in the theory of complements \cite{PrShokComplements}. For $n\ge2$ and a nontrivial group $G$ with invariant factors $d_1,\ldots,d_s$ as in \eqref{eq:groupdecomposition}, choose hyperplanes $L_1,\ldots,L_s\subset\PP^{n-1}$ in general position and put $D=\sum_{j=1}^s(1-1/d_j)L_j$. By Proposition~\ref{prop:exceptionalcomplexity}, the pair $(\PP^{n-1},D)$ is an exceptional klt log Calabi--Yau pair if and only if $c(G)=n$. This gives three possible lists of boundary coefficients on $\PP^1$ and six on $\PP^2$. Within this construction, smaller complexity gives non-exceptional log Fano pairs, see Subsection~\ref{subsec:complements}.

Theorem~\ref{thm:complexity} leads to the following question in arbitrary dimension.

\begin{question}\label{quest:complexity}
Let a finite abelian group $G$ act faithfully on a rationally connected variety $X$. Must $c(G)\le\dim X$?
\end{question}

By Proposition~\ref{prop:complexityrealization}, Question~\ref{quest:complexity} is equivalent to asking whether every finite abelian group acting faithfully on a rationally connected variety of dimension $n$ admits a faithful action on a klt Fano variety $X$ of dimension $n$ with $X/G\simeq\PP^n$.
The following proposition answers Question~\ref{quest:complexity} affirmatively for projective linear actions in every dimension.

\begin{proposition}\label{prop:complexityprojective}
Let $n\ge1$ and let $G\subseteq\operatorname{PGL}_{n+1}(\CC)$ be a finite abelian subgroup. Then
\begin{equation}\label{eq:complexityprojective}
 c(G)\le n.
\end{equation}
Equality holds if and only if $n=1$ and $G\simeq(\ZZ/2\ZZ)^2$.
\end{proposition}

Let a finite abelian $p$-group $G$ act faithfully on a rationally connected variety $X$ of dimension $n$. Applied to a smooth projective variety $G$-birational to $X$, the rank bound of Koll\'ar and Zhuang \cite[Corollary~11]{KZ} gives
\begin{equation}\label{eq:KZrank}
 \rk G\le\frac{p}{p-1}n\le2n.
\end{equation}
For $n=3$, this bound was obtained earlier in \cite{PrElementary,PrTwoElementary,PrShrP,Kuznetsova,XuRank,Loginov3}. For $G\simeq(\ZZ/p\ZZ)^s$, we have $c(G)=(1-1/p)s$, so \eqref{eq:KZrank} gives $c(G)\le n$ and answers Question~\ref{quest:complexity} for finite $p$-elementary groups.

For a finite group $G$, let $\exp G$ denote its exponent. In dimension three, the following corollary makes explicit the constants in the cubic bound implied by \cite[Corollary~3]{Golota}.

\begin{corollary}\label{cor:elementorder}
Let a finite group $G$ act faithfully on a rationally connected threefold $X$. Then
\[
 |G|\le2^{15}\,(\exp G)^3.
\]
If $X$ is rational, then the coefficient $2^{15}$ can be replaced by $3^9$.
\end{corollary}

For finite subgroups $G$ of $\Cr_1(\CC)$ and $\Cr_2(\CC)$, the corresponding bounds are $|G|\le2\exp G$ and $|G|\le288\,(\exp G)^2$, respectively, by \cite[\S2.2 and Lemma~3.3.1]{PrShrJordanConstant}.

\subsection*{Idea of the proof}

\emph{The empty anticanonical system.}
To prove Theorem~\ref{thm:nonemptymodel}, let a finite abelian group $G$ not of product type act faithfully on a rationally connected threefold $V$. Suppose that every $G\QQ$-Fano threefold in the $G$-birational class of $V$ has empty anticanonical linear system, and choose such a threefold $X$ by Corollary~\ref{cor:known-reduction}. Write $\tau_3(G)$ for the least index of a subgroup of $G$ generated by at most three elements, see Lemma~\ref{lem:tau}. The estimate of orbit lengths of non-Gorenstein points on $X$ and singular Riemann--Roch formula give $\tau_3(G)\in\{2,3\}$. If $h^0(X,-2K_X)\ne2$, the basket of singularities $\B(X)$ is one of the four baskets in Table~\ref{tab:baskets}, see Proposition~\ref{prop:baskets}. For the basket notation, see Subsection~\ref{subsec:terminal} and \eqref{eq:basketabbreviations}.

For the baskets $2Q_2+4Q_3+2Q_5$, $12Q_2+2Q_3$, and $14Q_2$, we use the equality $\Cl(X)^G=\ZZ[-K_X]$ of Proposition~\ref{prop:classes}. Every $G$-invariant bianticanonical pencil on $X$ has no fixed components and has canonical threshold at least $1/2$, see Subsection~\ref{sec:classes} and Corollaries~\ref{cor:pencils} and~\ref{cor:ct}.

For a $G$-invariant pencil $\Lambda\subseteq|-2K_X|$, Proposition~\ref{prop:nefmodel} gives a terminal threefold $Y$ that is $G$-birational to $X$. Write $\Lambda_Y$ for the transform of $\Lambda$ on $Y$. Then $-K_Y$ is nef and $(Y,\tfrac12\Lambda_Y)$ is terminal. We consider the four baskets separately.
\begin{enumerate}[label=\textup{(\arabic*)}]
\item $2Q_2+4Q_3+2Q_5$. We obtain a contradiction by constructing three $G$-invariant prime divisors $R_1,R_2,R_3\subset X$ and comparing the local contributions to $R_1\cdot R_2\cdot R_3$ with the lower bound obtained from their classes in $\Cl(X)^G$, see Proposition~\ref{prop:B0}.
\item $12Q_2+2Q_3$. We obtain a contradiction by constructing a $G$-invariant bianticanonical pencil with canonical threshold less than $1/2$, see Proposition~\ref{prop:B4}.
\item $14Q_2$. Let $S\to D$ be the minimal resolution of a general element $D\in|-2K_X|$. Then $S$ is a numerical Campedelli surface. The kernel of the action of $G$ on $|-2K_X|$ contains an involution~$\sigma$. We show that the reduced divisorial fixed locus of $\sigma$ on $X$ belongs to $|-K_X|$, giving a contradiction, see Proposition~\ref{prop:14}.
\item $15Q_2$. We use the bianticanonical map of $X$ to construct a $G$-invariant pencil $\Lambda\subseteq|-2K_X|$ with a base point of Gorenstein index two. For this choice of $\Lambda$, the threefold $Y$ satisfies $h^0(Y,-2K_Y)\le4$, see Proposition~\ref{prop:three}. Since $\tau_3(G)=3$, Proposition~\ref{prop:foursections} gives a contradiction to the assumption that $G$ is not of product type.
\end{enumerate}

If $h^0(X,-2K_X)=2$, take $\Lambda=|-2K_X|$. Then $\Lambda_Y=|-2K_Y|$ and $(-K_Y)^3=0$. Proposition~\ref{prop:foursections} again excludes $\tau_3(G)=3$. 
It remains to consider $\tau_3(G)=2$ and $\Lambda_Y=|-2K_Y|$. Put $H=\ker(G\to\Aut(\Lambda_Y))$, and let $T\to D$ be the minimal resolution of a general element $D\in\Lambda_Y$. Lemma~\ref{lem:surfacecoh} gives two cases.
\begin{itemize}
\item $\B(Y)=12Q_2$ and $q(T)=0$. The surface $T$ is Enriques or properly elliptic. The group $H$ embeds into $\Cr_2(\CC)$, so $G$ is of product type, see Corollary~\ref{cor:regularpencil}.
\item $\B(Y)=16Q_2$ and $q(T)=1$. The surface $D$ is smooth, and its Albanese morphism $a_D\colon D\to C_D$ has general fiber of genus one, where $C_D$ is the image curve of $a_D$, see Proposition~\ref{prop:anticanonical-Albanese}. We study the actions of $H$ on $C_D$ and of $\ker(H\to\Aut(C_D))$ on a general fiber of $a_D$. The resulting bounds give a contradiction, see Lemma~\ref{lem:relative-Albanese-surface} and Proposition~\ref{thm:irregular}.
\end{itemize}
This proves Theorem~\ref{thm:nonemptymodel}.

\emph{Reconstruction from an anticanonical K3 surface.}
For $G$ not of product type, take the $G\QQ$-Fano threefold $X$ with nonempty anticanonical linear system given by Theorem~\ref{thm:nonemptymodel}. Theorem~\ref{thm:nonempty} gives $G\simeq\Gamma_i$ for some $i\in\{1,2,3,4\}$. Choose a $G$-invariant anticanonical K3 surface $S=(s=0)$, and let $A$ be the ample $\QQ$-Cartier Weil divisor class on $S$ induced by $-K_X$. We identify $S$ with a hypersurface $F(y_1,\ldots,y_4)=0$ in a weighted projective space, with $A=[\OO_S(1)]$. Lemma~\ref{lem:restrictionring} gives
\[
 R(X,-K_X)/(s)\simeq R(S,A).
\]
The inertia group $I_S=\ker(G\to\Aut(S))$ implies that the defining relation of $R(X,-K_X)$ has the form $F+cs^d$, where $d=|I_S|$ and $c\in\CC$. If $c=0$, the cone vertex is $G$-fixed and $G$ is of product type by Proposition~\ref{prop:invariants}. Rescaling $s$ gives the Fermat hypersurface $X_i$ with its diagonal action, see Sections~\ref{sec:K3geometry} and~\ref{sec:reconstruction}. This proves Theorem~\ref{thm:rc-models}. The nonrationality of the four hypersurfaces and regularization give Theorem~\ref{thm:cremona}, and Lemma~\ref{lem:complexitytypes} gives the criterion in Theorem~\ref{thm:complexity}.

\subsection*{Acknowledgements}
The author is a Young Russian Mathematics award winner and would like to thank its sponsors and jury. OpenAI's ChatGPT was used to explore proof ideas, examine preliminary arguments and computations, search for references, and assist with drafting and editing the manuscript. The author takes full responsibility for its content.

\section{Preliminaries}\label{sec:prelim}
In this section, we recall some standard definitions and results on singularities, linear systems, group actions, and singular Riemann--Roch theorem. We follow \cite{KM} for the terminology of the minimal model program.

\subsection{Contractions}\label{subsec:actions}\label{subsec:contractions}
Let $G$ be a finite group. A $G$-variety is a variety equipped with an action of $G$. For a closed point $P$ of a $G$-variety $X$, we denote the stabilizer of $P$ by $G_P$.

A contraction is a projective morphism $f\colon X\to Z$ between normal varieties such that $f_*\OO_X=\OO_Z$. A fibration is a contraction with $\dim Z<\dim X$. A $G$-contraction is a contraction equivariant with respect to the actions of $G$ on $X$ and $Z$. 
For a contraction $f\colon X\to Z$, let $N^1(X/Z)$ be the real vector space of Cartier divisor classes modulo numerical equivalence on curves contracted by $f$. For a $G$-contraction, put $\rho(X/Z)^G=\dim_{\mathbb R}N^1(X/Z)^G$. 
A normal $G$-variety is $G\QQ$-factorial if every $G$-invariant Weil divisor is $\QQ$-Cartier.

\subsection{Singularities of pairs}\label{subsec:singularities}
A pair $(X,\Delta)$ consists of a normal variety $X$ and an effective Weil $\QQ$-divisor $\Delta$ with coefficients in $[0,1]$ such that $K_X+\Delta$ is $\QQ$-Cartier. Let $f\colon Y\to X$ be a log resolution of $(X,\Delta)$. Choose $K_Y$ so that $f_*K_Y=K_X$, and define $\Delta_Y$ by
\[
 K_Y+\Delta_Y=f^*(K_X+\Delta).
\]
For a prime divisor $E$ on $Y$, its discrepancy is $a(E,X,\Delta)=-\operatorname{coeff}_E\Delta_Y$, and its log discrepancy is $1+a(E,X,\Delta)$. The pair $(X,\Delta)$ is \emph{log canonical} (lc), respectively \emph{Kawamata log terminal} (klt), if $a(E,X,\Delta)\ge-1$, respectively $a(E,X,\Delta)>-1$, for every prime divisor $E$ over $X$.

For a normal variety $X$ such that $K_X$ is $\QQ$-Cartier, write $a(E,X)=a(E,X,0)$. The variety $X$ is \emph{canonical}, respectively \emph{terminal}, if $a(E,X)\ge0$, respectively $a(E,X)>0$, for every exceptional prime divisor $E$ over $X$. The Gorenstein index of a terminal point $P\in X$ is the least positive integer $r$ such that $rK_X$ is Cartier near $P$. For a terminal threefold $X$, we denote by $\Sing_{\nG}X$ the locus where $K_X$ is not Cartier.

\subsection{Linear systems}\label{subsec:linearsystems}
Let $X$ be a normal projective variety. For every $m\in\ZZ$, put $\plur_m(X)=h^0(X,mK_X)$.

For an integral Weil divisor $D$ on $X$, a linear system $\Lambda\subseteq|D|$ is a nonempty projective linear subspace of $|D|$. Write $\Lambda=F+\MM$, where $F$ is the largest effective Weil divisor satisfying $F\le D'$ for every $D'\in\Lambda$, and $\MM=\{D'-F\mid D'\in\Lambda\}$. We call $F$ the \emph{fixed part} and $\MM$ the \emph{moving part} of $\Lambda$. The linear system $\Lambda$ is \emph{mobile} if $F=0$.

Let $E$ be a prime divisor over $X$, and choose a proper birational morphism $f\colon Y\to X$ with $Y$ normal and $E\subset Y$. If the elements of $\Lambda$ are $\QQ$-Cartier, put $\ord_E\Lambda=\operatorname{coeff}_E(f^*D')$, where $D'\in\Lambda$ is a general element. Suppose that $K_X$ is $\QQ$-Cartier and $\MM$ is a mobile linear system of $\QQ$-Cartier divisors. For $c\ge0$, define
\[
 a(E,X,c\MM)=a(E,X)-c\ord_E\MM.
\]
If $X$ is canonical, the \emph{canonical threshold} of $\MM$ is
\[
 \ct(X,\MM)=\sup\{c\ge0\colon (X,c\MM)\text{ is canonical}\}.
\]
Here the mobile pair $(X,c\MM)$ is canonical, respectively terminal, if $a(E,X,c\MM)\ge0$, respectively $a(E,X,c\MM)>0$, for every exceptional prime divisor $E$ over $X$, see \cite[1.9]{Ale}.

\subsection{$G$-Mori fiber spaces}\label{subsec:mori}
For the minimal model program for varieties with finite group actions, see \cite{PrMMP}.

\begin{definition}\label{def:mori}
A $G\QQ$-Mori fiber space is a $G$-equivariant fibration $f\colon X\to Z$ such that the action of $G$ on $X$ is faithful, $X$ is terminal, projective, and $G\QQ$-factorial, the divisor $-K_X$ is $f$-ample, and
$
 \rho(X/Z)^G=1.
$
If $Z$ is a point, we call $X$ a $G\QQ$-Fano variety.
\end{definition}

We use the three-dimensional minimal model program for mobile pairs as in \cite[\S1]{Ale}, see also \cite[\S4.5]{PrMMP}. A canonical mobile pair admits a projective crepant terminalization \cite[Corollary~4.4.4]{PrMMP}. For a terminal or canonical mobile pair $(X,c\MM)$ on a terminal $G\QQ$-factorial threefold, the $(K_X+c\MM)$-negative minimal model program terminates by \cite[Theorem~5.2]{Shok}. The program yields a minimal model if $K_X+c\MM$ is pseudo-effective and a Mori fiber space otherwise. If $(X,c\MM)$ is terminal, every birational transform of $(X,c\MM)$ in the program is terminal.

\subsection{Finite abelian groups}\label{subsec:groups}
For a finite abelian group $G$ and a prime $p$, let $G_p$ denote the $p$-primary subgroup of $G$ and let $\rk G$ be the minimal number of generators of $G$. We call $G$ $p$-elementary if $G\simeq(\ZZ/p\ZZ)^r$ for some integer $r\ge0$. We put $G[n]=\{a\in G\colon na=0\}$ for $n\ge1$. For every integer $m\ge1$, put $\xi_m=\exp(2\pi\mathrm{i}/m)$.

\begin{definition}\label{def:producttype}
A finite abelian group $G$ is \emph{of product type} if $G\simeq G_1\times G_2$, where $G_i$ is a finite abelian subgroup of $\Cr_i(\CC)$ for $i=1,2$.
\end{definition}

\begin{lemma}\label{lem:liftcharacter}
Let
\[
 1\longrightarrow K\longrightarrow\widetilde G
 \longrightarrow G\longrightarrow1
\]
be an extension with $G$ abelian. Suppose there is a group homomorphism $\chi\colon\widetilde G\to\CC^\times$ whose restriction to $K$ is injective. Then $\widetilde G$ is abelian.
\end{lemma}
\begin{proof}
For any $a,b\in\widetilde G$, we have $[a,b]\in K\cap\ker\chi=\{1\}$.
\end{proof}

\begin{proposition}[{\cite[Proposition~3.11]{Log}}]\label{prop:extensions}
Let
\[
 1\longrightarrow N\longrightarrow G\longrightarrow Q\longrightarrow1
\]
be an exact sequence of finite abelian groups. If $N$ embeds into $\Cr_1(\CC)$ and $Q$ embeds into $\Cr_2(\CC)$, or conversely, then $G$ is of product type.
\end{proposition}

\begin{proposition}[{\cite[Corollary~3.14]{Log}}]\label{prop:mori}
Let a finite abelian group $G$ act faithfully on a rationally connected threefold $X$, and let $f\colon X\to Z$ be a $G\QQ$-Mori fiber space with $\dim Z>0$. Then $G$ is of product type.
\end{proposition}

We use the following reduction and classification in the proof of the main theorems.

\begin{corollary}\label{cor:known-reduction}
Let a finite abelian group $G$ not of product type act faithfully on a rationally connected threefold $V$. Then $V$ is $G$-birational to a $G\QQ$-Fano threefold.
\end{corollary}
\begin{proof}
Choose a smooth projective $G$-variety $\widetilde V$ with a $G$-equivariant birational map $\widetilde V\dashrightarrow V$, and run the $G$-equivariant minimal model program on $\widetilde V$. This gives a $G\QQ$-Mori fiber space $f\colon X\to Z$. By Proposition~\ref{prop:mori} we have $\dim Z=0$.
\end{proof}

\begin{theorem}[{\cite[Corollary~1.6]{LPZ}}]\label{thm:nonempty}
Let a finite abelian group $G$ act faithfully on a terminal Fano threefold with nonempty anticanonical linear system. Then $G$ is of product type or $G\simeq\Gamma_i$ for some $i\in\{1,2,3,4\}$.
\end{theorem}

For an integer $n\ge1$ and a finite abelian group $G$ with the decomposition \eqref{eq:groupdecomposition}, define
\begin{equation}\label{eq:taudefinition}
 \tau_n(G)=\begin{cases}d_1\cdots d_{s-n},&s>n,\\1,&s\le n.\end{cases}
\end{equation}

\begin{proposition}[{\cite[Proposition~8]{KZ}}]\label{prop:tauRC}
Let a finite abelian group $G$ act faithfully on a smooth rationally connected projective threefold $X$. Then
\begin{equation}\label{eq:tauRC}
 \tau_3(G)\mid24.
\end{equation}
\end{proposition}
\begin{proof}
Apply \cite[Proposition~8]{KZ} to the $p$-primary subgroups of $G$. See also \cite[Corollary~5.2]{Lo26}.
\end{proof}

\begin{lemma}\label{lem:tau}
Let $G$ be a finite abelian group and let $n\ge1$ be an integer. If $H\subseteq G$ is generated by at most $n$ elements, then $\tau_n(G)\mid[G\colon H]$, with equality for the subgroup corresponding to the product of the last $\min\{n,s\}$ cyclic factors under a chosen isomorphism in \eqref{eq:groupdecomposition}.
\end{lemma}
\begin{proof}
The invariant factors of $H_p$ divide the corresponding largest invariant factors of $G_p$ for every prime $p$; taking orders gives the divisibility.
\end{proof}

\begin{example}[{\cite[Example~1.8]{LPZ}}]\label{ex:exceptional-actions}
For each hypersurface $X_i$ below, let $m_j$ be the exponent of $x_j$ in its equation, for $1\le j\le4$. Identify $\Gamma_i$ with $\prod_{j=1}^4\ZZ/m_j\ZZ$. An element $(a_j)_{j=1}^4$ of $\Gamma_i$ fixes $x_0$ and multiplies each $x_j$ by $\xi_{m_j}^{a_j}$. This gives a faithful action on the rationally connected hypersurface $X_i$.
\begin{equation}\label{eq:exceptional-models}
\begin{aligned}
 X_1&=\{x_0^4+x_1^4+x_2^4+x_3^4+x_4^4=0\}\subset\PP^4,\\
 X_2&=\{x_0^6+x_1^6+x_2^6+x_3^6+x_4^2=0\}\subset\PP(1,1,1,1,3),\\
 X_3&=\{x_0^6+x_1^6+x_2^6+x_3^3+x_4^3=0\}\subset\PP(1,1,1,2,2),\\
 X_4&=\{x_0^8+x_1^8+x_2^8+x_3^4+x_4^2=0\}\subset\PP(1,1,1,2,4).
\end{aligned}
\end{equation}
The smooth quartic $X_1$ is not rational by \cite[Main Theorem]{IskovskikhManin}. The hypersurfaces $X_2,X_3,X_4$ are not rational by \cite[Corollary~1.1.9]{CP}. 
The last two groups in \eqref{eq:productlist} also have complexity three. For actions of these groups on complete intersections in a weighted projective space, see \cite[Example~1.9]{LPZ}.
\end{example}

\subsection{Stabilizers}\label{subsec:stabilizers}
\begin{lemma}[{\cite[Lemma~4]{Pop}}]\label{lem:tangent}
Let a finite group $G$ act faithfully on an irreducible variety $X$ and fix a closed point $P\in X$. Then the induced representation of $G$ on $T_PX$ is faithful.
\end{lemma}

\begin{corollary}\label{cor:smoothstabilizers}
Let a finite abelian group $G$ act faithfully on an irreducible variety $X$. Then for every smooth closed point $P\in X$, the stabilizer $G_P$ is generated by at most $\dim X$ elements.
\end{corollary}
\begin{proof}
By Lemma~\ref{lem:tangent}, the representation of $G_P$ on $T_PX$ is faithful. Since $G_P$ is finite abelian, this representation is diagonalizable, so $G_P$ embeds into $(\CC^\times)^{\dim X}$.
\end{proof}

\begin{proposition}[{\cite[Theorem~7.3]{Log}}]\label{prop:stabilizers}
Let $G$ be a finite abelian group acting faithfully on a terminal three-dimensional germ. Then $\rk G\le3$ unless $G$ has the form
\[
 (\ZZ/2\ZZ)^2\times \ZZ/2a\ZZ\times \ZZ/2b\ZZ
\]
for positive integers $a,b$. The case $\rk G=4$ occurs only at a Gorenstein germ of type $cA$.
\end{proposition}

\begin{proposition}[{\cite[Corollaries~7.4 and~3.12]{Log}}]\label{prop:invariants}
Let a finite abelian group $G$ act faithfully on a terminal threefold. If $G$ fixes a point or preserves a rational curve or a rational surface, then $G$ is of product type.
\end{proposition}

\subsection{Terminal threefold singularities}\label{subsec:terminal}\label{subsec:indexone}
Let $(X,P)$ be a terminal three-dimensional germ of Gorenstein index $r$. Its \emph{index-one cover} $\pi\colon(X^\#,P^\#)\to(X,P)$ is cyclic of degree $r$ with deck group $\ZZ/r\ZZ$. It is étale away from $P$, with $\pi^{-1}(P)=\{P^\#\}$, and the germ $(X^\#,P^\#)$ is Gorenstein and terminal. At a cyclic quotient point $P$, we have $(X^\#,P^\#)\simeq(\CC^3,0)$ in analytic coordinates in which the deck group acts diagonally. If $P^\#$ is singular, then $(X^\#,P^\#)$ embeds $(\ZZ/r\ZZ)$-equivariantly into $(\CC^4,0)$ as a hypersurface $\{\phi=0\}$. The coordinates and $\phi$ can be chosen semi-invariant under the deck group \cite[Theorem~3.2(a) and Proposition~3.6]{Reid87}.

A non-Gorenstein terminal germ $(X,P)$ admits a $\QQ$-Gorenstein deformation whose general fiber has only terminal cyclic quotient singularities. The analytic types of these singularities, counted with multiplicities, form the \emph{local basket} at $P$. We write $[r,b]$ for the cyclic quotient type
\begin{equation}\label{eq:baskettype}
 \tfrac1r(1,-1,b),\qquad r\ge2,\qquad1\le b\le r/2,\qquad\gcd(b,r)=1.
\end{equation}
For a terminal threefold $X$, the multiplicity of each type $[r,b]$ in the \emph{basket} $\B(X)$ is the sum of its multiplicities in the local baskets at all non-Gorenstein points of $X$. By the classification, the Gorenstein index of every element in the local basket at a point of Gorenstein index $r>1$ is $r$, except at a point of type $cAx/4$. In that case exactly one element has Gorenstein index $4$ and all others have Gorenstein index $2$. The local basket at a noncyclic non-Gorenstein point has at least two elements \cite[(6.4)]{Reid87}. In particular, the Gorenstein index of $P$ is the least common multiple of the Gorenstein indices of the elements in its local basket. We use the abbreviations
\begin{equation}\label{eq:basketabbreviations}
 Q_2=[2,1],\qquad Q_3=[3,1],\qquad Q_5=[5,2].
\end{equation}

Let a finite group $G$ act faithfully on a terminal germ $(X,P)$ of Gorenstein index $r$. For each $g\in G$, there is an automorphism $\widetilde g$ of $(X^\#,P^\#)$ such that $\pi\circ\widetilde g=g\circ\pi$. The group $\widetilde G$ of all lifts fits into the exact sequence
\begin{equation}\label{eq:indexcoverlifts}
 1\longrightarrow\ZZ/r\ZZ\longrightarrow\widetilde G
 \longrightarrow G\longrightarrow1.
\end{equation}
Here $\ZZ/r\ZZ$ acts on $(X^\#,P^\#)$ by deck transformations.

\subsection{Singular Riemann--Roch}\label{subsec:baskets}
Let $X$ be a projective terminal threefold and let $D$ be an integral $\QQ$-Cartier Weil divisor on $X$. Let $P\in\Sing_{\nG}X$ have Gorenstein index $r_P$. By \cite[Corollary~5.2]{KawamataCrepant}, there is a unique $i_P\in\{0,\ldots,r_P-1\}$ such that $D\sim i_PK_X$ locally at $P$. Choose a $\QQ$-Gorenstein deformation $\mathcal X\to\Delta$ of $(X,P)$ as in Subsection~\ref{subsec:terminal}. We extend the local class of $D$ by $i_P[K_{\mathcal X/\Delta}]$, whose restriction to a general fiber $\mathcal X_t$ is $i_P[K_{\mathcal X_t}]$. For each element $[r,b]$ in the local basket at $P$, let $P_t\in\mathcal X_t$ be the corresponding cyclic quotient point, where $K_{\mathcal X_t}$ has Cartier index $r$. At $P_t$, let $i_D$ be the least nonnegative residue of $i_P$ modulo $r$. Since $r\mid r_P$, the residue $i_D$ depends only on the class of $i_P$ modulo $r_P$. The class $i_P[K_{\mathcal X_t}]$ has Cartier index $r/\gcd(r,i_D)$ at $P_t$.

Write $\chi(D)=\chi(X,\OO_X(D))$. Singular Riemann--Roch theorem \cite[Theorem~10.2]{Reid87} reads
\begin{equation}\label{eq:RR}
 \chi(D)=\chi(\OO_X)+\frac{D(D-K_X)(2D-K_X)}{12}
 +\frac{Dc_2(X)}{12}+\sum_{[r,b]\in \B(X)}c_{r,b}(i_D),
\end{equation}
where
\begin{equation}\label{eq:correction}
 c_{r,b}(i)=-\frac{i(r^2-1)}{12r}
 +\sum_{j=1}^{i-1}\frac{\overline{bj}(r-\overline{bj})}{2r},
 \qquad0\le i<r,
\end{equation}
and the bar denotes the least nonnegative residue modulo $r$.

Applying \eqref{eq:RR} to $D=K_X$ and using Serre duality gives
\begin{equation}\label{eq:basketc2}
 (-K_X)c_2(X)+\sum_{[r,b]\in \B(X)}\left(r-\frac1r\right)=24\chi(\OO_X).
\end{equation}
If $-K_X$ is nef and big, Kawamata--Viehweg vanishing gives $\chi(\OO_X)=1$ and $\plur_{-m}(X)=\chi(-mK_X)$ for $m\ge0$. Equations~\eqref{eq:RR}--\eqref{eq:correction} then yield
\begin{equation}\label{eq:antipluri}
 \plur_{-m}(X)=\frac{m(m+1)(2m+1)}{12}(-K_X)^3+2m+1
 -\sum_{[r,b]\in \B(X)}\sum_{j=1}^m
 \frac{\overline{jb}(r-\overline{jb})}{2r}.
\end{equation}
In particular, the equality
\begin{equation}\label{eq:VG}
 (-K_X)^3=\sum_{[r,b]\in \B(X)}\frac{b(r-b)}r-6
\end{equation}
holds when $-K_X$ is nef and big and $\plur_{-1}(X)=0$.
If in addition $2K_X$ is Cartier, write $\B(X)=mQ_2$. Equations~\eqref{eq:antipluri} and~\eqref{eq:VG} give
\begin{equation}\label{eq:index-two-antipluri}
 m=12+2(-K_X)^3,\qquad \plur_{-2}(X)=m-10=2+2(-K_X)^3.
\end{equation}
In particular, we have $\plur_{-2}(X)\ge3$.

\section{Complexity of finite abelian groups}\label{sec:complexity}
In this section, we study finite covers $\pi\colon X\to\PP^n$ with abelian Galois group $G$. We express the minimum degree of the branch divisor in terms of $c(G)$ and $\exp G$, where the minimum is taken over all such covers $\pi$ for fixed $G$ and $n$. We show that $\pi$ can be chosen so that $X$ is a klt Fano variety if and only if $c(G)\le n$. We also relate the complexity $c(G)$ to the theory of complements and to exceptional log Calabi--Yau pairs.

\subsection{Group-theoretic properties}\label{subsec:complexityproperties}

\begin{lemma}\label{lem:complexity}
Let $H$ be a subgroup of a finite abelian group $G$. Then we have $c(H)\le c(G)$ and $c(G/H)\le c(G)$. The first inequality is strict if $H\ne G$, and the second is strict if $H\ne\{0\}$. Let
\[
 0\longrightarrow N\longrightarrow G\longrightarrow Q\longrightarrow0
\]
be an exact sequence of finite abelian groups. Then we have
\begin{equation}\label{eq:complexitysubadditivity}
 c(G)\le c(N)+c(Q).
\end{equation}
\end{lemma}
\begin{proof}
The invariant factors of $H$ divide the corresponding largest invariant factors of $G$. If $H\ne G$, then $H$ has fewer invariant factors or at least one of them is strictly smaller than the corresponding factor of $G$. This gives $c(H)\le c(G)$, with strict inequality for a proper subgroup. Every quotient of a finite abelian group is isomorphic to its subgroup, so the same argument applies to quotients.

Suppose first that $G/H$ is cyclic of order $m$. Write the invariant-factor decompositions as $H\simeq\prod_{i=1}^s\ZZ/h_i\ZZ$ and $G\simeq\prod_{i=1}^s\ZZ/g_i\ZZ$, allowing $h_i=1$ so that both decompositions have $s$ factors. Then $h_i\mid g_i\mid h_{i+1}$ for $i<s$ and $h_s\mid g_s$. Indeed, for every prime $p$ and every $a\ge0$, the group $p^aG_p$ is generated by $p^aH_p$ and at most one further element. Put $t_i=g_i/h_i$. Then $h_i\ge\prod_{j<i}t_j$ and $\prod_i t_i=m$, so
\begin{equation}\label{eq:complexitycyclic}
 c(G)-c(H)
 =\sum_{i=1}^s\frac{1-t_i^{-1}}{h_i}
 \le\sum_{i=1}^s\frac{1-t_i^{-1}}{\prod_{j<i}t_j}
 =1-\frac1m.
\end{equation}
Write $Q=Q_1\times\cdots\times Q_t$ with cyclic factors of orders $|Q_1|\mid\cdots\mid |Q_t|$. Let $\varphi\colon G\to Q$ be the quotient homomorphism and set $G_i=\varphi^{-1}(Q_1\times\cdots\times Q_i)$, with $G_0=N$. Since $G_i/G_{i-1}\simeq Q_i$ is cyclic, \eqref{eq:complexitycyclic} gives the bound $c(G_i)-c(G_{i-1})\le1-1/|Q_i|$. Summing over $1\le i\le t$ gives \eqref{eq:complexitysubadditivity}.
\end{proof}

\begin{lemma}\label{lem:complexitytypes}
Let $G$ be a finite abelian group. Then the following hold:
\begin{enumerate}[label=\textup{(\arabic*)}]
\item\label{item:complexity-low-dimension} For $i=1,2$, the group $G$ embeds into $\Cr_i(\CC)$ if and only if $c(G)\le i$.
\item The group $G$ is of product type if and only if $c(G)<3$ or
\[
 G\simeq(\ZZ/4\ZZ)^2\times(\ZZ/2\ZZ)^3
 \quad\text{or}\quad G\simeq(\ZZ/2\ZZ)^6.
\]
\item\label{item:complexity-three} We have $c(G)\le3$ if and only if $G$ is of product type or $G\simeq\Gamma_i$ for some $i\in\{1,2,3,4\}$.
\end{enumerate}
\end{lemma}
\begin{proof}
Let $s$ be the number of invariant factors of $G$ in \eqref{eq:groupdecomposition}. The condition $c(G)\le1$ holds precisely for cyclic groups and $(\ZZ/2\ZZ)^2$. If $c(G)\le2$, then $s\le4$. For $s=3$, the tuple of invariant factors $(d_1,d_2,d_3)$ is $(2,2,d_3)$ with $d_3$ even, $(2,4,4)$, or $(3,3,3)$. For $s=4$, all invariant factors equal two. Together with the groups of rank at most two, these are exactly the groups in \eqref{eq:A12}. This proves~(1).
Taking products of a cyclic group or $(\ZZ/2\ZZ)^2$ with the groups in \eqref{eq:A12} gives precisely the groups in \eqref{eq:productlist}. By Lemma~\ref{lem:complexity}, every group of product type has complexity at most three.

Conversely, suppose that $c(G)\le3$. Then $s\le6$. If $s\le3$, then $G$ is of product type. For $s=4$, the inequality $\sum_i1/d_i\ge1$ and the divisibilities in \eqref{eq:groupdecomposition} give
\[
 (d_1,d_2,d_3,d_4)=(2,2,d_3,d_4),\quad(2,4,4,d_4),\quad(3,3,3,d_4),
\]
or one of the four tuples corresponding to $\Gamma_1,\ldots,\Gamma_4$. For $s=5$, the only possibilities are $(2,2,2,2,d_5)$ and $(2,2,2,4,4)$; for $s=6$, all invariant factors equal two. The four groups $\Gamma_1,\ldots,\Gamma_4$ have complexity three. All the other possibilities are in \eqref{eq:productlist}, where precisely the last two groups have complexity three. This proves~(2) and~(3).
\end{proof}

By Lemmas~\ref{lem:complexity} and~\ref{lem:complexitytypes}, the classes of finite abelian groups admitting an embedding into $\Cr_2(\CC)$ and of groups of product type are closed under taking subgroups and quotients.

\subsection{Abelian covers of projective space}\label{subsec:complexitycovers}
Let $G$ be a finite abelian group and let $n\ge1$. A $G$-cover of $\PP^n$ is a finite quotient morphism $\pi\colon X\to\PP^n$ for a faithful $G$-action on a connected normal variety $X$. Let $B_1,\ldots,B_t\subset\PP^n$ be the prime divisors over which $\pi$ ramifies. For each $j$, choose a prime divisor $R_j\subset X$ with $\pi(R_j)=B_j$. Its \emph{inertia subgroup} $I_j\subseteq G$ consists of the elements fixing $R_j$ pointwise. Since $G$ is abelian, $I_j$ is independent of the choice of $R_j$. The group $I_j$ is cyclic of order $e_j$, the ramification index over $B_j$. The branch $\QQ$-divisor $\Delta_\pi$ on $\PP^n$ and its degree are given by
\begin{equation}\label{eq:branchboundary}
 \Delta_\pi=\sum_{j=1}^t\left(1-\frac1{e_j}\right)B_j,
 \qquad
 \deg\Delta_\pi=\sum_{j=1}^t\deg(B_j)\left(1-\frac1{e_j}\right).
\end{equation}
For the hyperplane class $H$ on $\PP^n$, the Hurwitz formula gives
\begin{equation}\label{eq:coverhurwitz}
 K_X=\pi^*(K_{\PP^n}+\Delta_\pi)\sim_{\QQ}(\deg\Delta_\pi-n-1)\pi^*H.
\end{equation}

We first construct a cover that minimizes $\deg\Delta_\pi$ among all $G$-covers $\pi\colon X\to\PP^n$. Let $G\ne\{0\}$ have invariant factors $d_1,\ldots,d_s$ as in \eqref{eq:groupdecomposition}. Choose hyperplanes $H_i=(\ell_i=0)$ in general position in $\PP^n$ for $0\le i\le s$.

Let $g_i$ generate the $i$th cyclic factor of $G$ for $1\le i\le s$, and put $g_0=-\sum_{i=1}^s g_i$. Over $H_i$, prescribe the inertia subgroup $\langle g_i\rangle$, with $g_i$ acting via multiplication by $\xi_{\ord(g_i)}$ on the normal line at a general point of a prime divisor above $H_i$. For the characters $\chi_i$ defined by $\chi_i(g_j)=\xi_{d_i}^{\delta_{ij}}$ for $1\le i,j\le s$, the construction of Alexeev and Pardini \cite[\S2.1, (2)]{AbelianCovers} requires line bundles $L_i$ on $\PP^n$ with isomorphisms
\[
 L_i^{\otimes d_i}\simeq\OO_{\PP^n}\bigl(H_i+(d_i-1)H_0\bigr),\qquad 1\le i\le s.
\]
Taking $L_i=\OO_{\PP^n}(1)$ gives a $G$-cover $\pi\colon X\to\PP^n$ with $X$ normal. Since $g_1,\ldots,g_s$ generate $G$, the variety $X$ is connected.

Take $\ell_i\ell_0^{d_i-1}$ as the section of $L_i^{\otimes d_i}$ defining the cover and trivialize $L_i$ by $\ell_0$ over $U=\PP^n\setminus H_0$. Then \cite[\S1.2, (4) and (5)]{NonNormalAbelianCovers} give
\begin{equation}\label{eq:hyperplane-cover}
 \CC(X)=\CC(\PP^n)(t_1,\ldots,t_s),
 \qquad t_i^{d_i}=\frac{\ell_i}{\ell_0},\quad 1\le i\le s.
\end{equation}

Since $g_0$ has order $d_s$, the branch divisor of $\pi$ is
\begin{equation}\label{eq:hyperplane-branch}
 \Delta_\pi=\left(1-\frac{1}{d_s}\right)H_0+\sum_{i=1}^s\left(1-\frac1{d_i}\right)H_i.
\end{equation}
The hyperplanes $H_i$ have simple normal crossings, so $(\PP^n,\Delta_\pi)$ is klt. Then $X$ is klt as well \cite[Proposition~5.20(4)]{KM}.

Put $R_i=(\pi^{-1}(H_i))_{\mathrm{red}}$ for $0\le i\le s$. Then $\pi^*H_0=d_sR_0$, and $\pi^*H_i=d_iR_i$ and $\operatorname{div}(t_i)=R_i-(d_s/d_i)R_0$ for $1\le i\le s$. The Hurwitz formula \eqref{eq:coverhurwitz} gives
\begin{equation}\label{eq:hyperplane-canonical}
 K_X\sim-(nd_s+1)R_0+\sum_{i=1}^s(d_i-1)R_i
 \sim\bigl(d_s(c(G)-n)-1\bigr)R_0.
\end{equation}

Let $g_{\min}(G)$ denote the least genus of a smooth connected projective curve admitting a faithful $G$-action where $G$ is a finite abelian group. If $c(G)$ is a positive integer, then
\[
 g_{\min}(G)=1+\frac{|G|}{2}\left(c(G)-1-\frac{1}{d_s}\right).
\]
For $c(G)\ge2$, this is Maclachlan's formula \cite[Theorem~4]{Maclachlan}. For $c(G)=1$, the group $G\simeq(\ZZ/2\ZZ)^2$ acts faithfully on $\PP^1$. 
We compute the minimum degree of the branch divisor among all $G$-covers of $\PP^n$.

\begin{proposition}\label{prop:minbranchdegree}
Let $G$ be a nontrivial finite abelian group, let $d_s$ be its largest invariant factor, and let $n\ge1$ be an integer. Then
\begin{equation}\label{eq:minbranchdegree}
 \min_\pi\deg\Delta_\pi=c(G)+1-\frac{1}{d_s},
\end{equation}
where $\pi$ ranges over all $G$-covers of $\PP^n$. The cover \eqref{eq:hyperplane-cover} attains the minimum.
\end{proposition}
\begin{proof}
Let $\pi\colon X\to\PP^n$ be a $G$-cover, and put $b_j=\deg B_j$. Since $\Cl(\PP^n)=\ZZ$, the inertia subgroups $I_j$ generate $G$ by \cite[\S2.1]{AbelianCovers}. By \cite[Theorem~1(1)]{AbelianCovers}, we can choose generators $g_j\in I_j$ such that $\sum_{j=1}^t b_jg_j=0$. Thus $d_s=\operatorname{lcm}(e_1,\ldots,e_t)$. Put
\[
 A=\prod_{j=1}^t(\ZZ/e_j\ZZ)^{b_j},\qquad
 \delta=(1,\ldots,1)\in A.
\]
Sending each of the $b_j$ standard generators of $(\ZZ/e_j\ZZ)^{b_j}$ to $g_j$ gives a surjection $A/\langle\delta\rangle\twoheadrightarrow G$. For each prime $p\mid d_s$, the $p$-primary component $\delta_p$ projects to a generator of a cyclic factor of $A_p$ of order $p^{v_p(d_s)}$. Denote this projection by $\rho_p$. Its restriction to $\langle\delta_p\rangle$ is an isomorphism, so $A_p=\langle\delta_p\rangle\times\ker\rho_p$. Since $\langle\delta\rangle=\prod_{p\mid d_s}\langle\delta_p\rangle$, we obtain $A=\langle\delta\rangle\times\prod_{p\mid d_s}\ker\rho_p$. The order of $\langle\delta\rangle$ is $d_s=\exp A$, so this factor can be taken as the last cyclic factor in the invariant-factor decomposition of $A$. Lemma~\ref{lem:complexity} yields
\[
 \deg\Delta_\pi
 \ge c(A)
 =c(A/\langle\delta\rangle)+1-\frac{1}{d_s}
 \ge c(G)+1-\frac{1}{d_s}.
\]
The lower bound is attained by \eqref{eq:hyperplane-cover}, since the boundary \eqref{eq:hyperplane-branch} has degree $c(G)+1-1/d_s$.
\end{proof}

\begin{proof}[Proof of Proposition~\ref{prop:complexityrealization}]
For the trivial group, we have $c(G)=0$ and may take $X=\PP^n$. For $G\ne\{0\}$, use the invariant factors in \eqref{eq:groupdecomposition}. 
If $c(G)\le n$, take the cover $\pi\colon X\to\PP^n$ in \eqref{eq:hyperplane-cover}. We have $\deg\Delta_\pi=c(G)+1-1/d_s<n+1$, so $-K_X$ is ample by \eqref{eq:coverhurwitz}. The klt Fano variety $X$ is rationally connected by \cite[Theorem~1]{ZhangRC}. 
Conversely, suppose that a klt Fano variety $X$ admits a faithful $G$-action with quotient morphism $\pi\colon X\to\PP^n$. Since $-K_X$ is ample, \eqref{eq:coverhurwitz} gives $\deg\Delta_\pi<n+1$. Proposition~\ref{prop:minbranchdegree} gives $c(G)+1-1/d_s\le\deg\Delta_\pi$, so $c(G)<n+1/d_s$. Since $d_sc(G)$ is an integer, we obtain $c(G)\le n$.
\end{proof}

\begin{remark}\label{rem:coverfermat}
When $s=n+1$, the covering variety $X$ in \eqref{eq:hyperplane-cover} is a Fermat hypersurface in weighted projective space. The unique linear relation among the $n+2$ forms $\ell_i$ has all coefficients nonzero by general position. Rescale the forms so that $\sum_{i=0}^s\ell_i=0$. Then the map $[\ell_0:\cdots:\ell_s]$ identifies $\PP^n$ with $H=(\sum_{i=0}^s y_i=0)\subset\PP^s$. If $\ell_i$ is replaced by $c_i\ell_i$ with $c_i\in\CC^\times$ for $0\le i\le s$, the roots in \eqref{eq:hyperplane-cover} can be chosen as $t_i'=(c_i/c_0)^{1/d_i}t_i$. Thus $\CC(\PP^n)(t_1',\ldots,t_s')=\CC(\PP^n)(t_1,\ldots,t_s)$. Put $w_i=d_s/d_i$ for $1\le i\le s$, and define the finite morphism
\begin{equation}\label{eq:coverpowermap}
 \Phi\colon \PP(1,w_1,\ldots,w_s)\longrightarrow\PP^s,
 \qquad [x_0:\cdots:x_s]\longmapsto[x_0^{d_s}:x_1^{d_1}:\cdots:x_s^{d_s}].
\end{equation}
The coordinate powers all have degree $d_s$, and each $x_i$ is integral over the subring they generate, so $\Phi$ is finite. Thus the inverse image of $H$ is
\begin{equation}\label{eq:coverfermathypersurface}
 Z=\Phi^{-1}(H)=\Bigl\{x_0^{d_s}+\sum_{i=1}^s x_i^{d_i}=0\Bigr\}
 \subset\PP(1,d_s/d_1,\ldots,d_s/d_{s-1},1).
\end{equation}
The affine cone over $Z$ is a hypersurface smooth away from the origin, hence normal by Serre's criterion. All irreducible components of the cone contain the origin, so the cone is integral and $Z$ is integral and normal.

On the chart $x_0\ne0$, the functions $t_i=x_i/x_0^{w_i}$ generate $\CC(Z)$. By \eqref{eq:coverpowermap}, we have $(\Phi|_Z)^*(y_i/y_0)=t_i^{d_i}$. Thus the finite map $\Phi|_Z\colon Z\to H\simeq\PP^n$ identifies $Z$ with the normalization in the field \eqref{eq:hyperplane-cover}. For $n=3$ and $G\simeq\Gamma_i$, this recovers $X_i$ with the diagonal action in Example~\ref{ex:exceptional-actions}, for $1\le i\le4$.
\end{remark}

\subsection{Exceptional log Calabi--Yau pairs}\label{subsec:complements}
A $\QQ$-complement of a projective pair $(Z,B)$ is a $\QQ$-boundary $B^+\ge B$ such that $(Z,B^+)$ is log canonical and $K_Z+B^+\sim_{\QQ}0$. For an integer $m\ge1$, a $\QQ$-boundary $B^+$ is an $m$-complement of $K_Z+B$ if it satisfies the following conditions \cite[Definition~2.10]{PrShokSecondComplements}.
\begin{enumerate}[label=\textup{(\arabic*)}]
\item The pair $(Z,B^+)$ is log canonical.
\item $m(K_Z+B^+)\sim0$.
\item $mB^+\ge m\lfloor B\rfloor+\lfloor(m+1)\{B\}\rfloor$, where $\{B\}=B-\lfloor B\rfloor$.
\end{enumerate}
The boundary $B$ has \emph{standard coefficients} if each of its coefficients belongs to $\{1-1/d\mid d\in\ZZ_{>0}\}\cup\{1\}$. If $B$ has standard coefficients and $m(K_Z+B^+)\sim0$, then the inequality in~(3) is equivalent to $B^+\ge B$ \cite[Remark~2.11 and Lemma~3.3]{PrShokSecondComplements}.

\begin{proposition}\label{prop:logCYquotient}
Suppose that $c(G)=n$ for an integer $n\ge1$. For the cover $\pi\colon X\to\PP^n$ in \eqref{eq:hyperplane-cover}, put $S=R_0$ and $\Delta=\Delta_\pi+H_0/d_s$.
Then $X$ is terminal Fano, and $(X,S)$ and $(\PP^n,\Delta)$ are log Calabi--Yau pairs. The divisor $S$ gives a $1$-complement of $K_X$, and $\Delta$ gives a $d_s$-complement of $K_{\PP^n}+\Delta_\pi$ of minimal index.
\end{proposition}
\begin{proof}
Since $\pi^*H_0=d_sS$, the Hurwitz formula \eqref{eq:coverhurwitz} gives
$
 K_X+S=\pi^*(K_{\PP^n}+\Delta).
$
The support of $\Delta$ has simple normal crossings and $H_0$ is its only component with coefficient one. Thus $(\PP^n,\Delta)$ is log canonical and $a(F,\PP^n,\Delta)>-1$ for every exceptional prime divisor $F$ over $\PP^n$. By the discrepancy formula for finite covers, the pair $(X,S)$ is log canonical and $a(E,X,S)>-1$ for every exceptional prime divisor $E$ over $X$ \cite[proof of Proposition~5.20]{KM}. Since $c(G)=n$ and $S=R_0$, the linear equivalence \eqref{eq:hyperplane-canonical} gives $K_X+S\sim0$. Thus $S$ gives a $1$-complement of $K_X$, and $-K_X\sim S=\pi^*H_0/d_s$ is ample.

Let $E$ be an exceptional prime divisor over $X$. Since $K_X+S$ is Cartier, the discrepancy $a(E,X,S)$ is an integer. The inequality $a(E,X,S)>-1$ therefore gives $a(E,X,S)\ge0$. If the center of $E$ on $X$ lies in $\Supp(S)$, then $a(E,X)=a(E,X,S)+\ord_E S>0$. Over $U$, the equations for $t_i$ in \eqref{eq:hyperplane-cover} define the cover $\pi^{-1}(U)\to U$. The open subset $\pi^{-1}(U)=X\setminus\Supp(S)$ is smooth because at each $p\in U$ the functions $\ell_i/\ell_0$ with $p\in H_i$ form part of a system of analytic coordinates on $U$ centered at $p$ and each function $\ell_i/\ell_0$ with $p\notin H_i$ admits a holomorphic $d_i$th root near $p$. Thus $X$ is terminal.

We have $\deg\Delta=1+c(G)=n+1$, so $K_{\PP^n}+\Delta\sim_{\QQ}0$. The divisor $d_s\Delta$ is integral, so $\Delta$ gives a $d_s$-complement of $K_{\PP^n}+\Delta_\pi$. If $m<d_s$, an $m$-complement would have coefficient one along both $H_0$ and $H_s$, whose coefficients in $\Delta_\pi$ are $1-1/d_s$. The degree would increase by at least $2/d_s$, whereas $\deg\Delta_\pi=n+1-1/d_s$. Thus the minimal index is $d_s$.
\end{proof}

For $n\ge2$, adjunction to $H_0$ in Proposition~\ref{prop:logCYquotient} gives a boundary on $\PP^{n-1}$ with coefficient $1-1/d_j$ along $H_0\cap H_j$. 
Recall that a projective pair is \emph{exceptional} if it admits a $\QQ$-complement and every $\QQ$-complement is klt \cite[Example~5.8]{PrShokSecondComplements}.

\begin{proposition}\label{prop:exceptionalcomplexity}
Let $d_1,\ldots,d_s$ be a sequence of integers satisfying $1< d_1\mid\cdots\mid d_s$. Let $n\ge2$ and let $L_1,\ldots,L_s\subset\PP^{n-1}$ be hyperplanes in general position. Put $B=\sum_{j=1}^s(1-1/d_j)L_j$.
Then the following conditions are equivalent.
\begin{enumerate}[label=\textup{(\arabic*)}]
\item $\sum_{j=1}^s(1-1/d_j)=n$;
\item $(\PP^{n-1},B)$ is a log Calabi--Yau pair;
\item $(\PP^{n-1},B)$ is exceptional.
\end{enumerate}
\end{proposition}
\begin{proof}
Since the $L_i$ are in general position, $(\PP^{n-1},B)$ is klt. The equality $\deg(K_{\PP^{n-1}}+B)=\deg B-n$ shows that this pair is log Calabi--Yau if and only if $\deg B=n$.

If $\deg B=n$, then $B$ is a $\QQ$-complement. Any $\QQ$-complement $B^+\ge B$ has the same degree as $B$, so $B^+=B$. Since $(\PP^{n-1},B)$ is klt, the pair is exceptional.

If $\deg B>n$, the degree of $K_{\PP^{n-1}}+B$ is positive, so no $\QQ$-complement exists. If $\deg B<n$, then $n-\deg B$ is a positive multiple of $1/d_s$ because every $d_j$ divides $d_s$. Let $B'$ be the boundary obtained from $B$ by raising the coefficient of $L_s$ to one. The remaining degree is $n-\deg B'=n-\deg B-1/d_s\ge0$. We add to $B'$ hyperplanes in general position with the $L_i$, each with coefficient $1/d_s$, to reach degree $n$. The resulting boundary $B^+$ is a $\QQ$-complement since it has simple normal crossings support and degree $n$. The coefficient of $L_s$ is one, so $(\PP^{n-1},B^+)$ is not klt. Thus exceptionality implies $\deg B=n$.
\end{proof}

For a nontrivial finite abelian group $G$, take $d_1,\ldots,d_s$ to be its invariant factors as in \eqref{eq:groupdecomposition}. Proposition~\ref{prop:exceptionalcomplexity} then shows that $c(G)=n$ if and only if the corresponding hyperplane pair $(\PP^{n-1},B)$ is an exceptional log Calabi--Yau pair.

\begin{remark}
For exceptional klt log Calabi--Yau pairs on $\PP^1$ and $\PP^2$ with boundary supported on hyperplanes in general position and coefficients $1-1/d_j$, Lemma~\ref{lem:complexitytypes} gives the following tuples with $1<d_1\mid\cdots\mid d_s$.
\[
\begin{aligned}
 \PP^1\quad & (3,3,3),\qquad (2,4,4),\qquad (2,2,2,2),\\
 \PP^2\quad & (4,4,4,4),\qquad (2,6,6,6),\qquad (3,3,6,6),\\
 & (2,4,8,8),\qquad (2,2,2,2,2,2),\qquad (2,2,2,4,4).
\end{aligned}
\]
Without the divisibility condition, exceptional log Calabi--Yau pairs of this form have one additional tuple of boundary denominators on $\PP^1$, namely $(2,3,6)$, and twelve on $\PP^2$ (see \cite[Example~2.10(i)]{PrShokComplements}). For such a pair, the least positive integer $m$ with $m(K_{\PP^{n-1}}+B)\sim0$ is $\operatorname{lcm}(d_1,\ldots,d_s)$.
\end{remark}

\subsection{Projective linear actions}\label{subsec:complexitylinear}

\begin{proof}[Proof of Proposition~\ref{prop:complexityprojective}]
For every finite abelian group $A$, the inequality $1-1/d\le\frac12\log_2d$ for $d\ge2$ gives
$
 c(A)\le\frac12\log_2|A|,
$
with equality if and only if $A$ is a $2$-elementary group.

Put $r=n+1$. We prove the assertion for finite abelian subgroups of $\operatorname{PGL}_r(\CC)$ by induction on $r$, starting with the trivial group for $r=1$. Let $\widetilde G\subseteq\operatorname{SL}_r(\CC)$ be the preimage of $G$ via the natural map $\operatorname{SL}_r(\CC)\to \operatorname{PGL}_r(\CC)$. The kernel $Z$ of $\widetilde G\to G$ consists of the scalar matrices in $\operatorname{SL}_r(\CC)$ and is isomorphic to $\ZZ/r\ZZ$.

Suppose that the representation of $\widetilde G$ on $\CC^r$ is irreducible, and let $\chi$ be its character. By Schur's lemma, the center of $\widetilde G$ is $Z$. Since $G$ is abelian, for each $\gamma\in\widetilde G\setminus Z$ there is an element $h\in\widetilde G$ with $h\gamma h^{-1}=\lambda\gamma$ for some $\lambda\ne1$. Taking traces gives $\chi(\gamma)=\lambda\chi(\gamma)$, hence $\chi(\gamma)=0$. Character orthogonality gives
\[
 1=\frac1{|\widetilde G|}\sum_{\gamma\in\widetilde G}|\chi(\gamma)|^2
 =\frac{r^2|Z|}{|\widetilde G|}
 =\frac{r^2}{|G|}.
\]
Thus $|G|=r^2$ and $c(G)\le\log_2r\le r-1$. For $r\ge2$, equality requires $r=2$ and $G\simeq(\ZZ/2\ZZ)^2$.

If the representation is reducible, choose a $\widetilde G$-invariant decomposition $\CC^r=V_1\oplus V_2$, where $r_i=\dim V_i>0$. Let $H$ be the image of $G$ in $\operatorname{PGL}(V_1)\times\operatorname{PGL}(V_2)$, let $G_i$ be its image in $\operatorname{PGL}(V_i)$, and let $K$ be the kernel of $G\to H$. Every element $\widetilde g$ of the inverse image of $K$ in $\widetilde G$ acts on each $V_i$ by a scalar $\lambda_i(\widetilde g)$. The homomorphism $g\mapsto\lambda_1(\widetilde g)/\lambda_2(\widetilde g)$, where $\widetilde g$ is any lift of $g\in K$, embeds $K$ into $\CC^\times$. Hence $K$ is cyclic. Lemma~\ref{lem:complexity} and induction give
\[
 c(G)\le c(K)+c(H)
 \le1-\frac1{|K|}+c(G_1)+c(G_2)
 \le r-1-\frac1{|K|}<r-1.
\]
Finally, $(\ZZ/2\ZZ)^2$ embeds into $\operatorname{PGL}_2(\CC)$ and has complexity one, so equality does occur.
\end{proof}

\subsection{Orders of finite groups}\label{subsec:finitegrouporders}

\begin{proof}[Proof of Corollary~\ref{cor:elementorder}]
Put $e=\exp G$. If $e\le2$, then $G$ is abelian and Theorem~\ref{thm:complexity} gives $|G|\le8e^3$. For $e\ge3$, choose an abelian subgroup $A\subseteq G$ with $[G:A]\le10368$ by \cite[Theorem~1.2.4]{PrShrJordanConstant}.

Since $\exp A\mid e$, the lists \eqref{eq:productlist} and \eqref{eq:exceptional-groups} give $|A|\le81$ for $e=3$, $|A|\le256$ for $e=4$, and $|A|\le2e^3$ for $e\ge5$. If $X$ is rational, they give $|A|\le128$ for $e=4$ and $|A|\le e^3$ for $e\ge5$. For $e=3$, the group $G$ is a $3$-group, so $[G:A]\le3^8$. For $e=4$, the group $G$ is a $2$-group, so $[G:A]\le2^{13}$. Consequently,
\[
 \frac{|G|}{e^3}
 \le\max\{8,3^8\cdot3,2^{13}\cdot4,2\cdot10368\}
 =2^{15}.
\]
If $X$ is rational, the same argument gives
\[
 \frac{|G|}{e^3}
 \le\max\{8,3^8\cdot3,2^{13}\cdot2,10368\}
 =3^9.
\]
This completes the proof.
\end{proof}

\subsection{A subgroup of \texorpdfstring{$\Cr_4(\CC)$}{Cr4(C)} not of product type}

\begin{example}\label{ex:rank-five}
We show that the group $(\ZZ/4\ZZ)^5$ embeds into $\Cr_4(\CC)$ but is not isomorphic to a product of two or more finite abelian subgroups of Cremona groups of positive ranks summing to four. Indeed, if $A$ is a direct factor of $(\ZZ/4\ZZ)^5$, then $A\simeq(\ZZ/4\ZZ)^a$ for some $a\ge0$. The classification of finite abelian subgroups of $\operatorname{PGL}_2(\CC)$ and Theorems~\ref{thm:blanc} and~\ref{thm:cremona} give $a\le r$ if $A$ embeds into $\Cr_r(\CC)$ with $r\in\{1,2,3\}$. A decomposition $(\ZZ/4\ZZ)^5\simeq\prod_i A_i$ with $A_i\simeq(\ZZ/4\ZZ)^{a_i}$, $A_i\subseteq\Cr_{r_i}(\CC)$, and $\sum_i r_i=4$ would give $5=\sum_i a_i\le\sum_i r_i=4$.

Let $S=\{w^2=x^4+y^4+z^4\}\subset\PP(1,1,1,2)$ be the degree-two del Pezzo surface, and let $T=\{r^4=x^4+y^4+z^4\}\subset\PP^3$. Let $H\subseteq\Aut(S)$ be the subgroup generated by the automorphisms $x\mapsto\xi_4 x$ and $y\mapsto\xi_4 y$ together with the involution $\sigma([x:y:z:w])=[-x:y:z:-w]$. Thus $H\simeq(\ZZ/4\ZZ)^2\times\langle\sigma\rangle$. Let $\widetilde K\simeq(\ZZ/4\ZZ)^3$ act on $T$ by multiplying $r,x,y$ separately by $\xi_4$, and let $\delta\in\widetilde K$ change the sign of $r$. Then $T/\langle\delta\rangle\simeq S$ via $w=r^2$.

Set $X=(S\times T)/\langle(\sigma,\delta)\rangle$ and $G=(H\times\widetilde K)/\langle(\sigma,\delta)\rangle$. The group $H\times\widetilde K$ acts faithfully on $\CC(S\times T)$, and $\CC(X)=\CC(S\times T)^{\langle(\sigma,\delta)\rangle}$. By the Galois correspondence, an element of $H\times\widetilde K$ fixes $\CC(X)$ pointwise precisely when it belongs to $\langle(\sigma,\delta)\rangle$. Thus $G$ acts faithfully on $X$. In $G$, the relation $(\sigma,\delta)=1$ identifies the direct factor $\langle\sigma\rangle$ of $H$ with $\langle\delta\rangle$ in $\widetilde K$, so $G\simeq(\ZZ/4\ZZ)^5$. The projection $f\colon X\to B=T/\langle\delta\rangle$ is $G$-equivariant. The action on $B$ gives the exact sequence $1\to H\to G\to\widetilde K/\langle\delta\rangle\to1$, where $\widetilde K/\langle\delta\rangle\simeq(\ZZ/4\ZZ)^2\times\ZZ/2\ZZ$. The group $G$ represents the exceptional isomorphism class in \cite[Theorem~1.12]{Filin}.

The morphism $f$ is a $G\QQ$-Mori fiber space with general fiber a del Pezzo surface of degree $2$. The fibers of $f$ are connected because, for every orbit $O$ of $\langle\delta\rangle$ on $T$, the group $\langle(\sigma,\delta)\rangle$ acts transitively on the surfaces $S\times\{t\}$ with $t\in O$. The involution $\sigma$ fixes the four points given by $x=w=0$, and $\delta$ fixes the smooth curve $r=0$ on $T$. Thus $X$ has cyclic quotient singularities of type $\frac12(1,1,1,0)$ and is $\QQ$-factorial. The Reid--Tai criterion gives terminality since $3/2>1$ \cite[Theorem~4.11]{Reid87}. Write $q\colon S\times T\to X$ for the quotient morphism and $\operatorname{pr}_S\colon S\times T\to S$ for the projection. Since $q$ is \'etale in codimension one and $K_T\sim0$, the Hurwitz formula gives $q^*(-K_X)\sim\operatorname{pr}_S^*(-K_S)$. Hence $-K_X$ is $f$-ample.

Pullback along $q$ followed by restriction to $S\times\{t\}$ for general $t\in T$ gives an injection $N^1(X/B)^G\hookrightarrow N^1(S)^H$. The group $H$ contains the Geiser involution $w\mapsto-w$, obtained by composing $\sigma$ with $x\mapsto-x$. The quotient of $S$ by this involution is $\PP^2$, so the subspace of $N^1(S)$ fixed by the involution has dimension $1$. The relative class of $-K_X$ is nonzero and $G$-invariant, hence $\rho(X/B)^G=1$.

By \cite[Theorem~4.4]{Blanc}, we can write $\CC(S)=\CC(u,v)$ with $\sigma(u)=-u$ and $\sigma(v)=v$. Put $a=r/z\in\CC(T)$, so $\delta(a)=-a$. Then
$
 \CC(X)=\CC(T)^{\langle\delta\rangle}(au,v).
$
The functions $au,v$ are invariant. Put $F=\CC(T)^{\langle\delta\rangle}(au,v)$. Since $a^2\in F$, $\CC(T)=\CC(T)^{\langle\delta\rangle}(a)$, and $u=(au)/a$, we have $F(a)=\CC(S\times T)$ and $[\CC(S\times T):F]=2$. Thus $F$ is the fixed field of $(\sigma,\delta)$. Since $T/\langle\delta\rangle$ is rational, this proves that $X$ is rational.
\end{example}

\section{Baskets of singularities}\label{sec:numerical}
In this section, we reduce the empty anticanonical case of Theorem~\ref{thm:nonemptymodel} to the case $\plur_{-2}(X)=2$ and four further cases distinguished by $\B(X)$. 

\begin{lemma}\label{lem:basketorbits}
Let a finite abelian group $G$ act faithfully on a projective terminal threefold $X$. Then $\tau_3(G)$ divides the length of every $G$-orbit in $\Sing_{\nG}X$ and the multiplicity of each singularity type $[r,b]$ in $\B(X)$.
\end{lemma}
\begin{proof}
At a non-Gorenstein point $P\in X$, Proposition~\ref{prop:stabilizers} gives $\rk G_P\le3$. Thus $\tau_3(G)\mid|G\cdot P|$ by Lemma~\ref{lem:tau}. For every $g\in G$, the germs $(X,P)$ and $(X,gP)$ are isomorphic, so their local baskets coincide. If a type $[r,b]$ occurs $e$ times in the local basket at $P$, the orbit $G\cdot P$ contributes $e|G\cdot P|$ to the multiplicity of $[r,b]$ in $\B(X)$. Sum over the $G$-orbits of non-Gorenstein points of~$X$.
\end{proof}

\begin{proposition}\label{prop:numerical}
Let $X$ be a terminal Fano threefold with $\plur_{-1}(X)=0$, and suppose that $e\in\{2,3,6\}$ divides the multiplicity of each cyclic quotient type $[r,b]$ in $\B(X)$. Then $e\ne6$. If $\plur_{-2}(X)\ne2$, the tuple $(e,\B(X),(-K_X)^3,\plur_{-2}(X))$ occurs in the first four columns of Table~\ref{tab:baskets}.

\begin{table}[H]
\normalfont
\centering
\renewcommand{\arraystretch}{1.15}
\begin{tabular}{clrrr}
\toprule
$e$&$\B(X)$&$(-K_X)^3$&$\plur_{-2}(X)$&GRDB\\
\midrule
2&$2Q_2+4Q_3+2Q_5$&$1/15$&0&13\\
2&$14Q_2$&$1$&4&258\\
2&$12Q_2+2Q_3$&$4/3$&4&263\\
3&$15Q_2$&$3/2$&5&264\\
\bottomrule
\end{tabular}
\par\vspace{2pt}
\caption{Baskets and numerical invariants}\label{tab:baskets}
\end{table}
The last column of Table~\ref{tab:baskets} gives the corresponding identifiers in the Graded Ring Database \cite{FanoDatabase,GradedRings}. 
\end{proposition}
\begin{proof}
Put
\[
 B=\sum_{[r,b]\in\B(X)}b,\qquad
 T=\sum_{[r,b]\in\B(X)}\frac{b^2}{r},\qquad
 C=\sum_{[r,b]\in\B(X)}\left(r-\frac1r\right).
\]
Since $2b\le r$, the inner sum in \eqref{eq:antipluri} for $m=2$ equals $3b/2-5b^2/(2r)$. Applying \eqref{eq:antipluri} for $m=1,2$ and using $\plur_{-1}(X)=0$ gives
\begin{equation}\label{num:basic}
 (-K_X)^3=B-T-6>0,\qquad \plur_{-2}(X)=B-10.
\end{equation}
By \eqref{eq:basketc2} and \cite[Corollary~7.3]{IJL}, we have
\begin{equation}\label{num:chern}
 C=24-(-K_X)c_2(X)<24.
\end{equation}
Put $\delta(r,b)=r-1/r-3b/2$. We have $\delta(2,1)=0$, $\delta(3,1)=7/6$, and $\delta(r,b)\ge3/2$ for every other type. Indeed, for $b\ge2$ we have $r\ge2b+1$, so $\delta(r,b)\ge b/2+1/2$, and for $b=1$, $r\ge4$ gives $\delta(r,1)>2$. Thus $3B/2\le C<24$, and \eqref{num:basic} yields
\begin{equation}\label{num:plurigenera}
 0\le \plur_{-2}(X)\le5,\qquad e\mid \plur_{-2}(X)+10.
\end{equation}
Writing $\B(X)=e\B'$, we obtain
\begin{equation}\label{num:delta}
 \sum_{[r,b]\in\B'}b=\frac{\plur_{-2}(X)+10}{e},\qquad
 \sum_{[r,b]\in\B'}\delta(r,b)<\frac{18-3\plur_{-2}(X)}{2e}.
\end{equation}
If $e=6$, then \eqref{num:plurigenera} gives $\plur_{-2}(X)=2$. By \eqref{num:delta}, we have $\sum_{\B'}b=2$ and $\sum_{\B'}\delta(r,b)<1$, so $\B(X)=12Q_2$. This gives $(-K_X)^3=0$ by \eqref{num:basic}, a contradiction.

For $e\in\{2,3\}$, \eqref{num:plurigenera} gives
\begin{equation}\label{num:smallP2}
 e=2\ \Longrightarrow\ \plur_{-2}(X)\in\{0,2,4\},\qquad
 e=3\ \Longrightarrow\ \plur_{-2}(X)\in\{2,5\}.
\end{equation}
If $e=2$ and $\plur_{-2}(X)=4$, then \eqref{num:delta} gives $\sum_{\B'}b=7$ and $\sum_{\B'}\delta(r,b)<3/2$. Thus every element of $\B'$ is $Q_2$ or $Q_3$, and at most one is $Q_3$, giving $\B(X)=14Q_2$ or $12Q_2+2Q_3$. If $e=3$ and $\plur_{-2}(X)=5$, then \eqref{num:delta} gives $\sum_{\B'}b=5$ and $\sum_{\B'}\delta(r,b)<1/2$. Hence every element of $\B'$ is $Q_2$, so $\B(X)=15Q_2$.

If $\plur_{-2}(X)=0$, then $e=2$ by \eqref{num:smallP2}. In \cite[Theorem~3.5, Table~A]{ChenChen}, the only basket in which every type has even multiplicity is the one in row~13, so $\B(X)=2Q_2+4Q_3+2Q_5$.
The volumes in Table~\ref{tab:baskets} follow from \eqref{num:basic}.
\end{proof}

For the proof of Theorem~\ref{thm:nonemptymodel}, we take $e=\tau_3(G)$ and make the following assumption.

\begin{assumption}\label{ass:counterexample}
Fix a faithful action of a finite abelian group $G$ not of product type on a rationally connected threefold $V$. Assume that
\begin{equation}\label{eq:geometriccontradiction}
\plur_{-1}(X)=0
\end{equation}
for every $G\QQ$-Fano threefold $X$ that is $G$-birational to $V$. In statements under this assumption, we let $X$ denote such a threefold.
\end{assumption}

\begin{lemma}\label{lem:parity}
Let $X$ be a $G\QQ$-Fano threefold under Assumption~\ref{ass:counterexample}. Then $4\nmid\tau_3(G)$, so $\tau_3(G)\in\{2,3,6\}$.
\end{lemma}
\begin{proof}
Choose a $G$-equivariant resolution $\pi\colon Z\to X$ that is an isomorphism over $X_{\reg}$. Proposition~\ref{prop:tauRC} applied to $Z$ gives $\tau_3(G)\mid24$. If $\tau_3(G)=1$, then $\rk G\le3$, so $G$ is of product type by \eqref{eq:productlist}, contrary to Assumption~\ref{ass:counterexample}.

Suppose $4\mid\tau_3(G)$. By Proposition~\ref{prop:stabilizers}, we have $2\mid[G\colon G_P]$ for every closed point $P\in X$. For every closed point $P\in Z$, Corollary~\ref{cor:smoothstabilizers} gives $\rk G_P\le3$, hence $4\mid[G\colon G_P]$ by Lemma~\ref{lem:tau}. Since $\OO_Z(-K_Z)$ is naturally $G$-linearized, \cite[Lemma~2.6]{Lo26} gives $4\mid(-K_Z)^3$. The singularities of $X$ are isolated \cite[Corollary~5.18]{KM}, so \cite[Lemma~4.2]{Lo26}, applied with $N=2$, and Kawamata--Viehweg vanishing give
\[
 \chi(Z,\OO_Z(-K_Z))\equiv\chi(X,\OO_X(-K_X))=0\pmod2.
\]
On the other hand, rational connectedness of $Z$ and Riemann--Roch give
\[
 \chi(Z,\OO_Z(-K_Z))=3+\frac{(-K_Z)^3}{2}\equiv1\pmod2,
\]
a contradiction.
\end{proof}

\begin{proposition}\label{prop:baskets}
Let $X$ be a $G\QQ$-Fano threefold under Assumption~\ref{ass:counterexample}. Then we have
\begin{equation}\label{eq:groupreduction}
\begin{gathered}
 G\simeq\ZZ/d_1\ZZ\times\ZZ/d_2\ZZ\times\ZZ/d_3\ZZ\times\ZZ/d_4\ZZ,\\
 d_1\mid d_2\mid d_3\mid d_4,\qquad d_1=\tau_3(G)\in\{2,3\}.
\end{gathered}
\end{equation}
Moreover, if $\plur_{-2}(X)\ne2$, then $(\tau_3(G),\B(X),(-K_X)^3,\plur_{-2}(X))$ occurs in the first four columns of a row of Table~\ref{tab:baskets}.
\end{proposition}
\begin{proof}
Put $e=\tau_3(G)$. Lemma~\ref{lem:parity} gives $e\in\{2,3,6\}$, and Lemma~\ref{lem:basketorbits} shows that $e$ divides the multiplicity of each type $[r,b]$ in $\B(X)$. By Proposition~\ref{prop:numerical}, we have $e\in\{2,3\}$, and if $\plur_{-2}(X)\ne2$, then $(e,\B(X),(-K_X)^3,\plur_{-2}(X))$ occurs in the first four columns of Table~\ref{tab:baskets}.

Write the decomposition of $G$ as in \eqref{eq:groupdecomposition}. Since $e>1$, we have $s>3$. If $s\ge5$, then $e=d_1\cdots d_{s-3}\ge4$, a contradiction. Thus $s=4$ and $d_1=e$.
\end{proof}

\begin{remark}
If $\plur_{-2}(X)=2$, the bounds \eqref{num:basic} and \eqref{num:delta} with $e=\tau_3(G)$ give nine possible baskets $\B(X)$, seven for $\tau_3(G)=2$ and two for $\tau_3(G)=3$.
However, the argument with the bianticanonical pencil $|-2K_X|$ does not require this finer description of $\B(X)$.
\end{remark}

\section{Bianticanonical pencils}\label{sec:divisors-pencils}
In this section, we bound the numerical classes of effective $G$-invariant divisors on $G\QQ$-Fano threefolds $X$ with empty anticanonical linear system and study $G$-invariant pencils $\Lambda\subseteq|-2K_X|$. Under Assumption~\ref{ass:counterexample}, we estimate the canonical thresholds of these pencils. Proposition~\ref{prop:nefmodel} then constructs a terminal rationally connected threefold $Y$ that is $G$-birational to $X$, with $-K_Y$ nef, such that the transform $\Lambda_Y$ of $\Lambda$ satisfies $\Lambda_Y\subseteq|-2K_Y|$ and $(Y,\tfrac12\Lambda_Y)$ is terminal.

\subsection{Invariant divisor classes}\label{sec:classes}

\begin{lemma}\label{lem:effectivedivisors}
Let a finite group $G$ act faithfully on a $G\QQ$-Fano threefold $X$ with $\plur_{-1}(X)=0$. Let $D$ be a nonzero effective Weil divisor with $[D]\in\Cl(X)^G$, and write $D\equiv a(-K_X)$. Then $a\ge1$, and if $a=1$, then $h^0(X,D)\le1$. If $\Cl(X)^G$ is torsion-free, then $a>1$.
\end{lemma}
\begin{proof}
The divisor $D$ is $\QQ$-Cartier, and $a>0$ because $-K_X$ is ample. If $0<a\le1$, applying the local inequality for the Riemann--Roch correction terms in \cite[proof of Lemma~2.12]{Ale} to \eqref{eq:RR} gives
\[
 \chi(-K_X)-\chi(-K_X-D)-\chi(D)+\chi(\OO_X)
 \ge(-K_X)D(-K_X-D).
\]
Kawamata--Viehweg vanishing gives $\chi(B)=h^0(X,B)$ for $B\in\{-K_X,-K_X-D,D,0\}$. Hence
\[
 0\le a(1-a)(-K_X)^3
 \le1-h^0(X,-K_X-D)-h^0(X,D)\le0.
\]
Thus $a=1$ and $h^0(X,D)\le1$. If $\Cl(X)^G$ is torsion-free and $a=1$, then $D+K_X\equiv0$. Since $D+K_X$ is $\QQ$-Cartier and $H^1(X,\OO_X)=0$, the class $[D+K_X]\in\Cl(X)^G$ is torsion. Hence $D\sim-K_X$, contrary to $\plur_{-1}(X)=0$.
\end{proof}

\begin{corollary}\label{cor:pencils}
Let a finite group $G$ act faithfully on a $G\QQ$-Fano threefold $X$ with $\plur_{-1}(X)=0$. Then every $G$-invariant pencil $\Lambda\subseteq|-2K_X|$ is mobile and has an integral general element.
\end{corollary}
\begin{proof}
Write $\Lambda=F+\MM$, where $F\ge0$ is the fixed divisor and $\MM$ is mobile. Suppose that $F\ne0$, and choose $M\in\MM$. The classes $[F]$ and $[M]$ are $G$-invariant. Write
\[
 F\equiv a(-K_X),\qquad M\equiv b(-K_X).
\]
Lemma~\ref{lem:effectivedivisors} gives $a\ge1$ and, since $h^0(X,M)\ge2$, also $b>1$. This contradicts $F+M\sim-2K_X$, which gives $a+b=2$. Therefore $F=0$.

Suppose that a general element of $\Lambda$ is reducible. Take a $G$-equivariant resolution $\pi\colon Z\to X$ of the rational map defined by $\Lambda$, and let $f\colon Z\to\PP^1$ be the induced $G$-equivariant morphism. Let
\[
 Z\xrightarrow{q}C\xrightarrow{\nu}\PP^1
\]
be the Stein factorization of $f$, and put $n=\deg\nu>1$. Since $Z$ is rationally connected, we have $C\simeq\PP^1$. For a general point $t\in C$, put $Z_t=q^{-1}(t)$. The general fiber $Z_t$ of $q$ is smooth and integral. It meets the locus where $\pi$ is an isomorphism, so $S=\pi_*Z_t$ is a prime divisor on $X$. The $G$-action on $C$ induced by the Stein factorization gives $g(S)=\pi_*Z_{g(t)}$. Since $Z_t\sim Z_{g(t)}$, we have $[S]\in\Cl(X)^G$. The divisors $\pi_*Z_t$ vary in $|S|$, so $h^0(X,S)\ge2$, and a general element of $\Lambda$ is linearly equivalent to $nS$. Writing $S\equiv a(-K_X)$, Lemma~\ref{lem:effectivedivisors} gives $a>1$, contrary to $na=2$. General elements of $\Lambda$ are reduced, hence integral.
\end{proof}

\begin{lemma}\label{lem:halftorsion}
Let $X$ be a terminal Fano threefold with $\plur_{-1}(X)=0$, and let a finite abelian group $G$ act faithfully on $X$. Suppose that a class $t\in\Cl(X)^G$ has order two and $i_t=0$ at every element $[r,b]\in \B(X)$ with $r>2$. Then $G$ is of product type.
\end{lemma}
\begin{proof}
The class $t$ is $\QQ$-Cartier and $h^0(X,t)=0$, since a nonzero effective divisor cannot be numerically trivial on a projective variety. Since $t-K_X\equiv-K_X$ is ample, Kawamata--Viehweg vanishing gives $H^i(X,t)=0$ for $i>0$. Thus $\chi(t)=0$. At each element $Q_2$ with $i_t=1$, the residues $i_{-K_X},i_t,i_{-K_X+t}$ are $1,1,0$, respectively. At every other basket element one has $i_t=0$. Hence \eqref{eq:RR} and~\eqref{eq:correction} give
\[
 \chi(-K_X+t)-\chi(-K_X)=\chi(\OO_X)-\chi(t)=1.
\]
Kawamata--Viehweg vanishing and $\plur_{-1}(X)=0$ therefore give $h^0(X,-K_X+t)=1$.

Choose a trivialization $\OO_X(2t)\simeq\OO_X$. Use it to define multiplication on $\OO_X\oplus\OO_X(t)$ by composing $\OO_X(t)\otimes\OO_X(t)\to\OO_X(2t)$ with this trivialization. Let $\pi\colon\widetilde X\to X$ be the relative spectrum of this algebra. Since $H^0(X,t)=0$, the cover is connected. It is normal by Serre's criterion, since the cover is \'etale over $X_{\reg}$ and both summands of its algebra are reflexive. The morphism $\pi$ is quasi-\'etale of degree two, and
\[
 \pi_*\OO_{\widetilde X}=\OO_X\oplus\OO_X(t),\qquad
 K_{\widetilde X}=\pi^*K_X.
\]
Since $\pi$ is finite and quasi-\'etale, the threefold $\widetilde X$ is terminal by \cite[Proposition~5.20]{KM}. The divisor $-K_{\widetilde X}=\pi^*(-K_X)$ is ample, so $\widetilde X$ is a terminal Fano threefold.

Since $t$ is $G$-invariant, each element of $G$ lifts to an automorphism of $\widetilde X$. The group of all lifts fits into a central exact sequence
\[
 1\longrightarrow\langle\sigma\rangle\simeq \ZZ/2\ZZ
 \longrightarrow\widetilde G\longrightarrow G\longrightarrow1.
\]
The decomposition of anticanonical sections induced by $\pi$ is
\[
 H^0(\widetilde X,-K_{\widetilde X})
 =H^0(X,-K_X)\oplus H^0(X,-K_X+t)=0\oplus\CC.
\]
The deck involution $\sigma$ acts on $H^0(\widetilde X,-K_{\widetilde X})$ as $-1$, so Lemma~\ref{lem:liftcharacter} shows that $\widetilde G$ is abelian. 
Since $G$ is a proper quotient of $\widetilde G$, Theorem~\ref{thm:nonempty} and Lemmas~\ref{lem:complexity} and~\ref{lem:complexitytypes} give
$
 c(G)<c(\widetilde G)\le3.
$
Thus $G$ is of product type.
\end{proof}

\begin{proposition}\label{prop:classes}
Let $X$ be a $G\QQ$-Fano threefold under Assumption~\ref{ass:counterexample}. Suppose that
\[
 \B(X)\in\{2Q_2+4Q_3+2Q_5,\ 12Q_2+2Q_3,\ 14Q_2\}.
\]
Then $\Cl(X)^G=\ZZ[-K_X]$.
\end{proposition}
\begin{proof}
Every Weil divisor $D$ on $X$ whose class is $G$-invariant is $\QQ$-Cartier, since $\sum_{g\in G}gD$ is $G$-invariant and linearly equivalent to $|G|D$. For a log terminal Fano variety, the groups $\Pic X$ and $\Cl X$ are finitely generated, the group $\Pic X$ is torsion-free, and numerical equivalence of $\QQ$-Cartier divisors agrees with $\QQ$-linear equivalence, see \cite[\S2.1]{PrDegree}. Thus $\Cl(X)^G$ has rank one and its numerically trivial classes form its torsion subgroup $\TT$.

Suppose that $t\in\TT$ has prime order $p$. For $1\le j<p$, the class $jt$ is nonzero torsion, so $h^0(X,jt)=0$. Since $jt-K_X\equiv-K_X$ is ample, Kawamata--Viehweg vanishing gives $\chi(jt)=h^0(X,jt)=0$. If $i_t\ne0$ at a basket element $[r,b]$, then $p\mid r$. If $p$ is odd, sum \eqref{eq:RR} for the classes $jt$, $0\le j<p$. For a basket element $[r,b]$ with $i_t\ne0$, the residues $i_{jt}$ run through $0,r/p,\ldots,(p-1)r/p$. The classes $jt$ and $(p-j)t$ have residues $i=i_{jt}$ and $r-i$ for $1\le j<p$. By \eqref{eq:correction}, the sum of their corrections is $c_{r,b}(i)+c_{r,b}(r-i)=-\overline{bi}(r-\overline{bi})/(2r)$. Since multiplication by $b$ permutes the nonzero multiples of $r/p$ modulo $r$, the total contribution of $[r,b]$ is $-r(p^2-1)/(24p)$. Hence $\sum_{i_t\ne0}r=24p/(p+1)$. Only elements with $p\mid r$ occur in this sum, and their indices add to at most $12$ in each of the three baskets. This gives
\[
 18\le\frac{24p}{p+1}
 =\sum_{\substack{[r,b]\in \B(X)\\i_t\ne0}}r\le12.
\]
Thus $p=2$, and $i_t=0$ at every basket element of odd Gorenstein index. In each of the three baskets, every element $[r,b]$ with $r$ even is of type $Q_2$. Therefore Lemma~\ref{lem:halftorsion} implies that $G$ is of product type, contrary to Assumption~\ref{ass:counterexample}. Hence $\TT=0$ and $\Cl(X)^G\simeq\ZZ$.

Choose a generator $A$ of $\Cl(X)^G$ and write
$
 -K_X=kA$, $k\in\ZZ_{>0}.
$
Let $r$ be the least common multiple of the Gorenstein indices of the element in $\B(X)$. For every $[q,b]\in\B(X)$, the congruence $ki_A\equiv-1\pmod q$ implies $\gcd(k,q)=1$, so $\gcd(k,r)=1$. The Gorenstein index of each point of $X$ divides $r$, so $rA$ is Cartier by \cite[Corollary~5.2]{KawamataCrepant}. Hence
\[
 (rA)^3=\frac{r^3(-K_X)^3}{k^3}\in\ZZ.
\]
For $\B(X)=2Q_2+4Q_3+2Q_5$, $12Q_2+2Q_3$, and $14Q_2$, the pairs $(r,r^3(-K_X)^3)$ are
$
 (30,1800),$$(6,288),$ and $(2,8),
$
respectively. In all three cases, every prime divisor of $r^3(-K_X)^3$ divides $r$. Since $\gcd(k,r)=1$, we obtain $k=1$.
\end{proof}

\subsection{Canonical thresholds}\label{sec:ct}

\begin{proposition}\label{prop:ct}
Let a finite group $G$ act faithfully on a $G\QQ$-Fano threefold $X$. Suppose that $\plur_{-1}(Y)=0$ for every $G\QQ$-Fano threefold $Y$ that is $G$-birational to $X$. Let $\Lambda\subseteq|-2K_X|$ be a mobile $G$-invariant pencil. Then at least one of the following holds:
\begin{enumerate}
\item $\ct(X,\Lambda)\ge1/2$;
\item there is a $G\QQ$-Mori fiber space $h\colon X'\to Z$ such that $X'$ is $G$-birational to $X$ and $\dim Z>0$.
\end{enumerate}
\end{proposition}
\begin{proof}
Suppose that $(X,\frac12\Lambda)$ is not canonical. Take a $G$-equivariant terminal model $f\colon X_1\to X$ such that $X_1$ is $G\QQ$-factorial, the mobile pair $(X_1,\frac12\Lambda_1)$ is terminal, and $K_{X_1}+\frac12\Lambda_1$ is $f$-nef, see \cite[\S\S1,~3.4]{Ale}. Here $\Lambda_1$ is the strict transform of $\Lambda$. The negativity lemma gives a nonzero effective $G$-invariant $f$-exceptional divisor $E$ such that
\begin{equation}\label{eq:thresholdE}
 K_{X_1}+\tfrac12\Lambda_1
 =f^*(K_X+\tfrac12\Lambda)-E
 \sim_{\QQ}-E.
\end{equation}
Let $D\in\Lambda$ be a general element, and let $D_1$ be its strict transform. Then
\[
 2K_{X_1}+D_1+2E=f^*(2K_X+D).
\]
The divisor $2K_X+D$ is principal, so $2E$ is integral.

The divisor $-E$ is not pseudo-effective. Run a $G$-equivariant $(K_{X_1}+\frac12\Lambda_1)$-MMP. On each threefold $X_j$ in the program, let $\Lambda_j$ and $E_j$ be the transforms of $\Lambda_1$ and $E$. By \eqref{eq:thresholdE}, the contracted curves have negative intersection with $-E_j$. A divisorial contraction $\psi\colon X_j\to X_{j+1}$ with $E_j\ne0$ and $E_{j+1}=0$ would make $E_j$ effective, $\psi$-exceptional, and $\psi$-nef, contradicting the negativity lemma. Hence every $E_j$ is nonzero.

The program ends with a contraction $h\colon X'\to Z$. Let $\Lambda'$ and $E'$ be the transforms of $\Lambda_1$ and $E$ on $X'$. The class $-(K_{X'}+\tfrac12\Lambda')$ is $h$-ample. For a general $D'\in\Lambda'$, choose an integral curve $C$ in a general fiber of $h$ meeting $X'\setminus\Supp D'$. Then $D'\cdot C\ge0$. Since $[D']\in N^1(X'/Z)^G$ and $\rho(X'/Z)^G=1$, the divisor $D'$ is $h$-nef. Adding $\tfrac12 D'$ to the $h$-ample class $-(K_{X'}+\tfrac12 D')$ shows that $-K_{X'}$ is $h$-ample. If $\dim Z>0$, we obtain \textup{(ii)}. Suppose that $Z$ is a point. Then $X'$ is a $G\QQ$-Fano threefold and
\[
 K_{X'}+\tfrac12\Lambda'+E'\sim_{\QQ}0,
 \qquad E'\ne0,
\]
Since $2E'$ is integral and $\plur_{-1}(X')=0$, the $G$-invariant pencil
$
 \Lambda'+2E'\subseteq|-2K_{X'}|
$
has the nonzero fixed divisor $2E'$, contrary to Corollary~\ref{cor:pencils}.
\end{proof}

\begin{corollary}\label{cor:ct}
Let $X$ be a $G\QQ$-Fano threefold under Assumption~\ref{ass:counterexample}. Let $\Lambda\subseteq|-2K_X|$ be a $G$-invariant pencil. Then we have
$
 \ct(X,\Lambda)\ge\frac12.
$
\end{corollary}
\begin{proof}
By Corollary~\ref{cor:pencils}, every $G$-invariant pencil $\Lambda\subseteq|-2K_X|$ is mobile. If Proposition~\ref{prop:ct} gives a $G\QQ$-Mori fiber space $h\colon X'\to Z$ with $\dim Z>0$, then $G$ is of product type by Proposition~\ref{prop:mori}, contrary to Assumption~\ref{ass:counterexample}. Thus $\ct(X,\Lambda)\ge1/2$.
\end{proof}

\begin{lemma}\label{lem:terminalpencil}
Let $Y$ be a projective terminal threefold, and let $\Lambda\subseteq|-2K_Y|$ be a mobile pencil with integral general element. If $(Y,\tfrac12\Lambda)$ is terminal, then $2K_Y$ is Cartier and a general element $D\in\Lambda$ is normal, Gorenstein, and disjoint from $\Sing_{\nG}Y$.
\end{lemma}
\begin{proof}
Arguing as in \cite[proof of Lemma~1.22]{Ale}, let $P\in\Sing_{\nG}Y$ have Gorenstein index $r\ge2$. Choose an exceptional prime divisor $E$ over $Y$ with center $P$ and discrepancy $a(E,Y)=1/r$ \cite[p.~193]{KawamataDiscrepancy}. Let $D\in\Lambda$ be a general element, choose a resolution $\mu\colon W\to Y$ on which $E$ is a prime divisor, and let $D_W$ be the strict transform of $D$. Put $v=\ord_E\Lambda=\operatorname{coeff}_E(\mu^*D)$. Effectivity of $D$ and terminality of $(Y,\tfrac12\Lambda)$ give $0\le v<2/r$. Since $2K_Y+D$ is principal, the coefficient of $E$ in $2K_W+D_W-\mu^*(2K_Y+D)$ is integral. Thus $2/r-v\in\ZZ$, so $r=2$ and $v=0$. Now $D$ is Cartier near $P$, and $v=0$ implies $P\notin D$. Thus $2K_Y$ is Cartier and $D\cap\Sing_{\nG}Y=\varnothing$.

For an irreducible curve $C\subset\Bs\Lambda$, the blowup along $C$ near its generic point has exceptional divisor of discrepancy one over $Y$. Terminality of $(Y,\tfrac12\Lambda)$ gives $\operatorname{mult}_C\Lambda<2$. Hence $\operatorname{mult}_C\Lambda=1$, so $D$ is smooth at the generic point of $C$. Bertini's theorem gives smoothness on $Y_{\reg}\setminus\Bs\Lambda$. As a Cartier divisor in the Gorenstein locus of $Y$, the divisor $D$ is Gorenstein. It is regular in codimension one and Cohen--Macaulay, hence normal.
\end{proof}

\begin{proposition}\label{prop:nefmodel}
Let $X$ be a $G\QQ$-Fano threefold under Assumption~\ref{ass:counterexample}, and let $\Lambda\subseteq|-2K_X|$ be a $G$-invariant pencil. Then there is a projective terminal rationally connected threefold $Y$ with a faithful $G$-action and a $G$-equivariant birational map $\varphi\colon X\dashrightarrow Y$. Let $\Lambda_Y$ be the birational transform of $\Lambda$ under $\varphi$. Then
\begin{enumerate}[label=(\arabic*)]
\item $-K_Y$ is nef,
\item $2K_Y$ is Cartier,
\item $\plur_{-1}(Y)=0$,
\item $2\le\plur_{-2}(Y)\le\plur_{-2}(X)$,
\item $\Lambda_Y\subseteq|-2K_Y|$, and
\item $(Y,\tfrac12\Lambda_Y)$ is terminal.
\end{enumerate}
A general element $D\in\Lambda_Y$ is integral, normal, Gorenstein, and disjoint from $\Sing_{\nG}Y$. If $\plur_{-2}(X)=2$, then
\begin{equation}\label{eq:pencilnef}
 \plur_{-2}(Y)=2,\qquad (-K_Y)^3=0.
\end{equation}
\end{proposition}
\begin{proof}
By Corollary~\ref{cor:ct}, the mobile pair $(X,\frac12\Lambda)$ is canonical. Take a $G$-equivariant crepant terminalization $f\colon Y_0\to X$ of $(X,\frac12\Lambda)$, and let $\Lambda_0$ be the strict transform of $\Lambda$. The relation
\[
 K_{Y_0}+\tfrac12\Lambda_0=f^*(K_X+\tfrac12\Lambda)
\]
implies $\Lambda_0\subseteq|-2K_{Y_0}|$. 

Choose a sufficiently small positive rational $\delta$ such that $(Y_0,(1/2+\delta)\Lambda_0)$ remains canonical, and run a $G$-equivariant $(K_{Y_0}+(1/2+\delta)\Lambda_0)$-negative minimal model program. Since every member of $\Lambda_0$ is effective and $K_{Y_0}+(1/2+\delta)\Lambda_0\sim_{\QQ}\delta\Lambda_0$, the class $K_{Y_0}+(1/2+\delta)\Lambda_0$ is pseudo-effective. The transform of $K_{Y_0}+(1/2+\delta)\Lambda_0$ is nef on the resulting threefold $Y$. Let $\Lambda_j$ be the birational transform of $\Lambda_0$ on a threefold $Y_j$ occurring in the program. A divisorial contraction in the program on $Y_j$ would contract curves covering its exceptional divisor on which $K_{Y_j}+(1/2+\delta)D_j$ has negative degree, where $D_j\in\Lambda_j$ is general. Since $K_{Y_j}+\tfrac12 D_j\sim_{\QQ}0$, these curves would have negative intersection with $D_j$. The exceptional divisor would therefore be a fixed component of $\Lambda_j$. Thus all steps are small. Since $K_{Y_j}+\frac12\Lambda_j\sim_{\QQ}0$, these steps are crepant for the pairs $(Y_j,\frac12\Lambda_j)$, which remain terminal. On $Y$, the divisor $-K_Y$ is nef. We also have $\Lambda_Y\subseteq|-2K_Y|$, since the birational map $Y_0\dashrightarrow Y$ is small.
The construction gives the $G$-equivariant diagram
\[
\begin{tikzcd}[column sep=large]
 (X,\tfrac12\Lambda)
 & (Y_0,\tfrac12\Lambda_0) \arrow[l,"f"'] \arrow[r,dashed,"\mathrm{MMP}"]
 & (Y,\tfrac12\Lambda_Y)
\end{tikzcd}
\]

By Corollary~\ref{cor:pencils}, a general element of $\Lambda$ is integral, so its strict transform on $Y$ is integral. Lemma~\ref{lem:terminalpencil} gives that $2K_Y$ is Cartier and that a general element of $\Lambda_Y$ is normal, Gorenstein, and disjoint from $\Sing_{\nG}Y$.

For each $m\ge0$, the morphism $f$ induces an injection $H^0(Y_0,-mK_{Y_0})\hookrightarrow H^0(X,-mK_X)$. Since $f$ is crepant for the pair $(X,\frac12\Lambda)$, a general member $D_X\in\Lambda$ and its strict transform $D_0$ satisfy $2K_{Y_0}+D_0=f^*(2K_X+D_X)$. Thus the image of the injection for $m=2$ contains the sections defining the general members of $\Lambda$ and hence the two-dimensional subspace defining $\Lambda$. The map $Y_0\dashrightarrow Y$ is an isomorphism in codimension one, so it identifies $H^0(Y_0,-mK_{Y_0})$ with $H^0(Y,-mK_Y)$ for all $m\ge0$. This proves~(3) and~(4).

If $\plur_{-2}(X)=2$, then $\plur_{-2}(Y)=2$ by~(4). If $(-K_Y)^3>0$, then $-K_Y$ is nef and big, and \eqref{eq:index-two-antipluri} gives $\plur_{-2}(Y)>2$. Hence $(-K_Y)^3=0$.
\end{proof}

\section{Excluding the baskets \texorpdfstring{$2Q_2+4Q_3+2Q_5$ and $12Q_2+2Q_3$}{2Q2+4Q3+2Q5 and 12Q2+2Q3}}\label{sec:localexclusions}
In this section, we study the divisorial components of the fixed loci of elements of $G$ on $X$ and the actions of point stabilizers on the corresponding index-one covers. We use these actions to exclude the baskets $2Q_2+4Q_3+2Q_5$ and $12Q_2+2Q_3$.

\subsection{Divisorial fixed loci}\label{sec:fixed}

\begin{lemma}\label{lem:fixeddivisor}
Let a finite abelian group $G$ act faithfully on a $G\QQ$-Fano threefold $X$. For $1\ne g\in G$, let $R(g)$ be the reduced divisorial part of $\Fix(g)$. If $R(g)\ne0$, then $R(g)$ is a single normal $G$-invariant prime divisor.
\end{lemma}
\begin{proof}
Since $G$ is abelian, the divisor $R=R(g)$ is $G$-invariant. Thus $R$ is $\QQ$-Cartier. The numerical class of $R$ is a positive multiple of $-K_X$, since $\rho(X)^G=1$ and $R$ is nonzero and effective. In particular, the divisor $R$ is ample. 
The divisor $R$ is Cohen--Macaulay by \cite[Corollary~5.25]{KM}. The fixed locus of $g$ on $X_{\reg}$ is smooth, and $X$ has isolated singularities. Hence $R$ satisfies $R_1$ and $S_2$, so $R$ is normal. 
Since $R$ is ample and $\QQ$-Cartier, we have $H^1(X,\OO_X(-R))=0$ by \cite[Theorem~2.70]{KM}. The sequence $0\to\OO_X(-R)\to\OO_X\to\OO_R\to0$ gives $H^0(R,\OO_R)=\CC$. Thus $R$ is connected and hence prime.
\end{proof}

Let a finite abelian group $G$ act faithfully on a $G\QQ$-Fano threefold $X$, and let $P\in X$ be a cyclic quotient terminal point of Gorenstein index $r$. Let $\pi\colon(U^\#,P^\#)\to(U,P)$ be the index-one cover of a sufficiently small $G_P$-invariant analytic neighborhood $U$ of $P$, with deck group $\ZZ/r\ZZ$. Let $\widetilde G_P$ be the group of all lifts of $G_P$ to $(U^\#,P^\#)$ in the exact sequence \eqref{eq:indexcoverlifts}.

\begin{lemma}\label{lem:cyclicquotientlifts}
The group $\widetilde G_P$ is abelian, and its action on $(U^\#,P^\#)\simeq(\CC^3,0)$ is linear and diagonal in suitable analytic coordinates. If an integer $q\ge2$ satisfies $|G_P[q]|\ge q^3$, then, in these coordinates, $\widetilde G_P[q]$ is the group of matrices $\operatorname{diag}(\xi_q^{a_1},\xi_q^{a_2},\xi_q^{a_3})$ with $a_i\in\ZZ/q\ZZ$. In particular,
$
 \widetilde G_P[q]\simeq(\ZZ/q\ZZ)^3.
$
\end{lemma}
\begin{proof}
Let $\rho_P\colon\widetilde G_P\to\operatorname{GL}(T_{P^\#}U^\#)$ be the representation induced by the action on the index-one cover (see Subsection~\ref{subsec:indexone}). In cyclic quotient coordinates, choose the generator $\delta$ of the deck group so that the eigenvalues of $\rho_P(\delta)$ are $\xi_r,\xi_r^{-1},\xi_r^b$, where $\gcd(b,r)=1$. Thus $(\det\rho_P)(\delta)=\xi_r^b$, so the character $\det\rho_P$ is injective on the deck group. Lemma~\ref{lem:liftcharacter} implies that $\widetilde G_P$ is abelian. Linearizing the $\widetilde G_P$-action on $(U^\#,P^\#)$ and diagonalizing $\rho_P$ gives the asserted coordinates. 
The surjection $\widetilde G_P\to G_P$ induces a surjection modulo $q$, so
\[
 \begin{aligned}
 q^3&\le|G_P[q]|=|G_P/qG_P|\le|\widetilde G_P/q\widetilde G_P|
     =|\widetilde G_P[q]|\le q^3.
 \end{aligned}
\]
The faithful tangent representation embeds $\widetilde G_P[q]$ into the group of diagonal matrices with entries in $\langle\xi_q\rangle$, which has order $q^3$. Equality of orders proves the assertion.
\end{proof}

Choose $\widetilde G_P$-eigenvector coordinates $x_1,x_2,x_3$ on $(U^\#,P^\#)\simeq(\CC^3,0)$, and put $H=\{x_1=0\}\subset U^\#$.

\begin{lemma}\label{lem:localplane}
Suppose that $1\ne g\in G_P$ has a lift $\widetilde g\in\widetilde G_P$ such that $\widetilde g|_H=\id_H$. Then there is a unique normal $G$-invariant prime divisor $R\subset X$ fixed pointwise by $g$, and
$
 \pi^*(R|_U)=H.
$
\end{lemma}
\begin{proof}
Since $\pi$ is finite and $\widetilde g$ fixes $H$ pointwise, the image $\pi(H)$ is a surface contained in $\Fix(g)\cap U$. The fixed locus $\Fix(g)$ is a proper closed algebraic subset of $X$, so $\pi(H)$ lies in its divisorial part $R(g)$. By Lemma~\ref{lem:fixeddivisor}, the divisor $R=R(g)$ is normal, prime, and $G$-invariant. Normality makes the analytic germ of $R$ at $P$ irreducible, so $R\cap U=\pi(H)$ after shrinking $U$. The deck group preserves $H$, so $\pi^{-1}(\pi(H))=H$. The cover is \'etale at the general point of $H$, giving $\pi^*(R|_U)=H$.
\end{proof}

\begin{lemma}\label{lem:oddtriple}
Let $g_1,g_2,g_3$ be distinct involutions in a finite abelian group $G$ acting faithfully on a terminal threefold $X$. Suppose $g_1g_2g_3=1$ and let $R_i$ be nonzero reduced divisors on $X$ fixed pointwise by $g_i$. Then every point of $R_1\cap R_2\cap R_3$ has even Gorenstein index.
\end{lemma}
\begin{proof}
Suppose that $P\in R_1\cap R_2\cap R_3$ has odd Gorenstein index $r$. Choose a $\langle g_1,g_2\rangle$-invariant analytic neighborhood $U$ of $P$ such that $\mathscr L=\OO_U(rK_U)$ is invertible. Since $\langle g_1,g_2\rangle$ fixes $P$, the canonical $\langle g_1,g_2\rangle$-linearization of $\mathscr L$ induces a character
$
 \chi\colon\langle g_1,g_2\rangle\longrightarrow\CC^\times
$
on the fiber $\mathscr L|_P$. We determine $\chi(g_i)$ by examining how $g_i$ acts on $\mathscr L$ over a divisor it fixes pointwise.

Fix $i\in\{1,2,3\}$, and let $D_i$ be an irreducible component of $R_i\cap U$ through $P$. Since $g_i$ fixes $D_i$ pointwise, it acts on each fiber of $\mathscr L|_{D_i}$ by a scalar in $\{1,-1\}$. Since $D_i$ is irreducible, this scalar has the same value $\eps_i$ on every fiber. At a general point $Q\in D_i$, both $X$ and $D_i$ are smooth. The action of $g_i$ on $T_QX$ has eigenvalues $1,1,-1$ by Lemma~\ref{lem:tangent}. Its action on $\det T_Q^*X$ is therefore $-1$, so its action on $\mathscr L|_Q$ is $(-1)^r=-1$. Thus $\chi(g_i)=\eps_i=-1$ for all three involutions. This contradicts $g_1g_2g_3=1$, since
$
 1=\chi(g_1g_2g_3)=\chi(g_1)\chi(g_2)\chi(g_3)=-1.
$
The proof is completed.
\end{proof}

\subsection{The basket \texorpdfstring{$2Q_2+4Q_3+2Q_5$}{2Q2+4Q3+2Q5}}\label{sec:B0}

\begin{proposition}\label{prop:B0}
Let $X$ be a $G\QQ$-Fano threefold under Assumption~\ref{ass:counterexample}. Then $\B(X)\ne2Q_2+4Q_3+2Q_5$.
\end{proposition}
\begin{proof}
Suppose $\B(X)=2Q_2+4Q_3+2Q_5$. By Assumption~\ref{ass:counterexample} and \eqref{num:basic}, we have $\plur_{-1}(X)=\plur_{-2}(X)=0$. Proposition~\ref{prop:baskets} and Table~\ref{tab:baskets} give $\tau_3(G)=2$ and $(-K_X)^3=1/15$.
Proposition~\ref{prop:classes} gives $\Cl(X)^G=\ZZ[-K_X]$.

The points of $X$ of type $\frac12(1,-1,1)$ form one $G$-orbit $\{P_+,P_-\}$. Indeed, every $G$-orbit in $\Sing_{\nG}X$ has even length by Lemma~\ref{lem:basketorbits} and $\tau_3(G)=2$. Suppose that the local basket at a noncyclic point $P$ contains $Q_2$. Since $\B(X)$ has no element of Gorenstein index $4$, the description in Subsection~\ref{subsec:terminal} gives at least two elements of type $Q_2$ in the local basket at $P$. The orbit $G\cdot P$ would then give at least four elements of type $Q_2$ in $\B(X)$, contrary to the multiplicity of $Q_2$ in $\B(X)$. Thus every point whose local basket contains $Q_2$ is a cyclic quotient point. Every point of even Gorenstein index contributes a basket element of even Gorenstein index, so $P_+,P_-$ are the only points of even Gorenstein index.

The stabilizer $G_{P_+}$ has index two in $G$ and is generated by at most three elements by Proposition~\ref{prop:stabilizers}. Since $\rk G_2=4$, we have $\rk(G_{P_+})_2=3$ and hence $|G_{P_+}[2]|=8$. Let $\pi\colon(U^\#,P^\#)\to(U,P_+)$ be the smooth index-one cover of a $G_{P_+}$-invariant analytic neighborhood $U$ of $P_+$, and let $\widetilde G_{P_+}$ be the group of all lifts of $G_{P_+}$ in \eqref{eq:indexcoverlifts}. By Lemma~\ref{lem:cyclicquotientlifts}, choose $\widetilde G_{P_+}$-eigenvector coordinates $x_1,x_2,x_3$ in which $\widetilde G_{P_+}[2]=\{\operatorname{diag}(\eta_1,\eta_2,\eta_3)\mid\eta_i\in\{1,-1\}\}$. The deck involution is $\delta=-I_3$. In particular, $\widetilde G_{P_+}$ contains the elements
\[
 \eps_1=\diag(-1,1,1),\quad
 \eps_2=\diag(1,-1,1),\quad
 \eps_3=\diag(1,1,-1).
\]
The images $g_i\in G_{P_+}$ of $\eps_i$ are distinct non-trivial involutions with $g_1g_2g_3=1$. Let $R_i=R(g_i)$ be the reduced divisorial part of $\Fix(g_i)$ on $X$. Since $\eps_i$ fixes $\{x_i=0\}$ pointwise, Lemma~\ref{lem:localplane} gives a normal prime divisor $R_i$ with $\pi^*(R_i|_U)=\{x_i=0\}$. The element $g_i$ fixes $R_i$ pointwise. Since $G$ is abelian, it preserves $\Fix(g_i)$ and hence $R_i$. The planes $\pi^*(R_i|_U)=\{x_i=0\}$ are pairwise distinct, so the divisors $R_i$ are pairwise distinct.

Write $R_i\sim-a_iK_X$ with $a_i\in\ZZ_{>0}$. Since $\plur_{-1}(X)=\plur_{-2}(X)=0$, we have $a_i\ge3$ for each $i$. By Lemma~\ref{lem:oddtriple}, every point of $R_1\cap R_2\cap R_3$ has even Gorenstein index, so $R_1\cap R_2\cap R_3\subseteq\{P_+,P_-\}$. The equality $\pi^*(R_i|_U)=\{x_i=0\}$ gives $P_+\in R_i$ for each $i$. Since each $R_i$ is $G$-invariant and $P_+,P_-$ form one $G$-orbit, the triple intersection is zero-dimensional and its support is exactly $\{P_+,P_-\}$. The three planes $\pi^*(R_i|_U)$ have intersection multiplicity one at $P^\#$, and $\pi$ has degree two. Hence
$
 (R_1R_2R_3)_{P_+}=\frac12.
$
By $G$-invariance, the same equality holds at $P_-$. This gives the contradiction
\[
 1=R_1R_2R_3=a_1a_2a_3(-K_X)^3
 =\frac{a_1a_2a_3}{15}\ge\frac{27}{15}>1.
\]
The proof is completed.
\end{proof}

\subsection{The basket \texorpdfstring{$12Q_2+2Q_3$}{12Q2+2Q3}}\label{sec:B4}

\begin{proposition}\label{prop:B4}
Let $X$ be a $G\QQ$-Fano threefold under Assumption~\ref{ass:counterexample}. Then $\B(X)\ne12Q_2+2Q_3$.
\end{proposition}
\begin{proof}
Suppose $\B(X)=12Q_2+2Q_3$. By \eqref{num:basic}, we have $\plur_{-2}(X)=4$, so \eqref{eq:groupreduction} and Table~\ref{tab:baskets} give $d_1=\tau_3(G)=2$.
Assumption~\ref{ass:counterexample} excludes groups of product type, so \eqref{eq:productlist} gives $d_2>2$.

By Lemma~\ref{lem:basketorbits}, every $G$-orbit in $\Sing_{\nG}X$ has even length. The local basket description in Subsection~\ref{subsec:terminal} shows that a noncyclic point whose local basket contains $Q_3$ has Gorenstein index three and at least two copies of $Q_3$ in its local basket. Its $G$-orbit would therefore contribute at least four copies of $Q_3$, contrary to $\B(X)=12Q_2+2Q_3$. Thus the points whose local baskets contain $Q_3$ are two cyclic quotient points of type $\frac13(1,-1,1)$ in a single $G$-orbit.

Fix one of these points $P$ and put $q=d_2$. Since $[G\colon G_P]=2$ and $|G[q]|=2q^3$, we have $|G_P[q]|\ge q^3$. Let $\pi\colon(U^\#,P^\#)\to(U,P)$ be the smooth index-one cover of a $G_P$-invariant analytic neighborhood $U$ of $P$, and let $\widetilde G_P$ be the group of all lifts of $G_P$ to $U^\#$. Identify $(U^\#,P^\#)$ with $(\CC^3,0)$. By Lemma~\ref{lem:cyclicquotientlifts}, choose $\widetilde G_P$-eigenvector coordinates $x_1,x_2,x_3$ on $\CC^3$ in which $\widetilde G_P[q]$ consists of all diagonal matrices with entries in $\langle\xi_q\rangle$. Since $P$ has type $\frac13(1,-1,1)$, we can choose the generator $\delta$ of the deck group so that
\[
 \delta^*x_1=\xi_3 x_1,\qquad
 \delta^*x_2=\xi_3^{-1}x_2,\qquad
 \delta^*x_3=\xi_3 x_3.
\]
For each $i\in\{1,2,3\}$, choose an element $\widetilde g_i\in\widetilde G_P[q]$ acting by $\widetilde g_i^*x_i=\xi_q x_i$ and $\widetilde g_i^*x_j=x_j$ for $j\ne i$. Its image $g_i\in G_P$ is non-trivial, and $\widetilde g_i$ fixes the plane $\{x_i=0\}\subset U^\#\simeq\CC^3$ pointwise. By Lemmas~\ref{lem:fixeddivisor} and~\ref{lem:localplane}, there are distinct normal $G$-invariant prime divisors $R_1,R_2,R_3\subset X$ such that
$ \pi^*(R_i|_U)=\{x_i=0\}$ for $i=1,2,3$. 

Choose a $G$-eigenbasis $s_0,\ldots,s_3$ of $H^0(X,-2K_X)$ and put $D_j=(s_j=0)$. Set $\omega=dx_1\wedge dx_2\wedge dx_3$ and write the reflexive pullback of $s_j$ to $U^\#$ as $\pi^*s_j=f_j\omega^{-2}$, so that $\operatorname{div}(f_j)=\pi^*(D_j|_U)$. Since $\pi\circ\delta=\pi$, the section $\pi^*s_j$ is $\delta$-invariant. The equality $\delta^*\omega=\xi_3\omega$ therefore gives
\begin{equation}\label{eq:cubiccoefficients}
 \delta^*f_j=\xi_3^2f_j.
\end{equation}
Both $\pi^*s_j$ and $\omega$ are eigenvectors for $\widetilde G_P[q]$, so $f_j$ is also an eigenvector. The power series expansion of $f_j$ has the form
\[
 f_j=x_1^{a_{1j}}x_2^{a_{2j}}x_3^{a_{3j}}
 F_j(x_1^q,x_2^q,x_3^q),\qquad 0\le a_{ij}<q.
\]
By Proposition~\ref{prop:classes}, write $R_i\sim-a_iK_X$ with $a_i\in\ZZ_{>0}$. Since $\plur_{-1}(X)=0$, we have $a_i\ge2$.

Assume that $a_{ij}>0$ for some $i,j$. Then $R_i\le D_j$, so $D_j-R_i$ is effective and
\[
 0\le(-K_X)^2(D_j-R_i)=(2-a_i)(-K_X)^3.
\]
Thus $a_i=2$, and the effective linearly trivial divisor $D_j-R_i$ is zero. Hence $f_j=x_i u_j$ with $u_j(0)\ne0$. Comparing linear terms in \eqref{eq:cubiccoefficients} gives $i=2$. If $a_{ij}>0$ for some $i$, then $D_j=R_2$. Since the $s_j$ are linearly independent, this holds for at most one index $j$.

For at least three indices $j$, we have $a_{ij}=0$ for all $i$, so $f_j=F_j(x_1^q,x_2^q,x_3^q)$. Evaluating \eqref{eq:cubiccoefficients} at the origin gives $f_j(0)=0$ for these $j$. Choose two such indices $j_1,j_2$. The pencil $\Lambda$ generated by $s_{j_1},s_{j_2}$ is mobile by Corollary~\ref{cor:pencils}. Let $E$ be the exceptional divisor of the toric blowup at $P$ corresponding to the ray through $\frac13(1,2,1)$. The weights of $x_1,x_2,x_3$ are $1/3,2/3,1/3$, respectively. Each monomial in $f_{j_1}$ or $f_{j_2}$ has weight at least $q/3$, and the discrepancy of $E$ is $(1+2+1)/3-1$. Hence
\[
 a(E,X)=\frac13,\qquad\ord_E\Lambda\ge\frac q3,
 \qquad\ct(X,\Lambda)\le\frac1q<\frac12.
\]
This contradicts Corollary~\ref{cor:ct} and completes the proof.
\end{proof}

\section{Bianticanonical divisors}\label{sec:bianticanonical-surfaces}
In this section, we study general bianticanonical divisors on terminal Fano threefolds and on the terminal rationally connected threefolds from Proposition~\ref{prop:nefmodel}, together with the minimal resolutions of these divisors. We then use the holomorphic Lefschetz formula to restrict finite abelian groups acting on surfaces with $p_g=q=0$.

\subsection{The general bianticanonical divisor}\label{sec:cohom}

\begin{lemma}\label{lem:bigmember}
Let $X$ be a Fano threefold with terminal singularities. Suppose that $2K_X$ is Cartier and that
$ \plur_{-1}(X)=0$, $(-K_X)^3>1/2.$
Then the following hold for a general element $D\in|-2K_X|$.
\begin{enumerate}[label=(\arabic*)]
\item The divisor $D$ is integral and normal with only Du Val singularities and avoids $\Sing_{\nG}X$.
\item The minimal resolution $\nu\colon T\to D$ satisfies
$ p_g(T)=q(T)=0$, $K_T\text{ nef and big}$, and $K_T^2=2(-K_X)^3$.
\item If $s_D\in H^0(X,-2K_X)$ defines $D$, then restriction to $D$, adjunction, and pullback along $\nu$ induce an isomorphism 
\begin{equation}\label{eq:bigrestriction}
 H^0(X,-2K_X)/\langle s_D\rangle\simeq H^0(T,2K_T).
\end{equation}
\end{enumerate}
\end{lemma}
\begin{proof}
The divisor $-2K_X$ is ample Cartier, with $K_X^2(-2K_X)=2(-K_X)^3>1$ and $-K_X=\frac12(-2K_X)$. By \cite[Proposition~6.1]{PrDegree}, a general element $D\in|-2K_X|$ is integral and normal, has only Du Val singularities, and is disjoint from $\Sing_{\nG}X$. Adjunction gives $\omega_D\simeq\OO_D(-K_X)$. The minimal resolution $\nu\colon T\to D$ satisfies $K_T=\nu^*K_D$. Thus $K_T\sim\nu^*((-K_X)|_D)$ is nef and big, and $K_T^2=2(-K_X)^3$.

Kawamata--Viehweg vanishing gives $H^i(X,\OO_X)=0$ for $i>0$ and, together with $\plur_{-1}(X)=0$, gives $H^i(X,-K_X)=0$ for all $i$. Serre duality gives $H^i(X,2K_X)=0$ for all $i$. The sequence
\[
 0\longrightarrow\OO_X(2K_X)\longrightarrow\OO_X\longrightarrow\OO_D\longrightarrow0
\]
then gives $p_g(D)=q(D)=0$. By rationality of Du Val singularities \cite[Theorem~4.20 and Corollary~5.24]{KM}, we have $H^i(T,\OO_T)\simeq H^i(D,\OO_D)$, so $p_g(T)=q(T)=0$. Finally, the exact sequence
\[
 0\longrightarrow\OO_X\xrightarrow{\ s_D\ }\OO_X(-2K_X)
 \longrightarrow\OO_D(2K_D)\longrightarrow0
\]
and $H^1(X,\OO_X)=0$ give the restriction isomorphism.
\end{proof}

\begin{lemma}\label{lem:ruled}
Let a finite abelian group $G$ act faithfully on a rationally connected threefold $X$. Let $\Lambda$ be a $G$-invariant pencil on $X$ with integral general element $D$, and let $\nu\colon T\to D$ be a resolution. If $T$ is birationally ruled, then $G$ is of product type.
\end{lemma}
\begin{proof}
Let $\varphi_\Lambda\colon X\dashrightarrow\PP^1$ be the rational map defined by $\Lambda$. Take a $G$-equivariant resolution $\pi\colon Y\to X$ on which $\varphi_\Lambda$ becomes a $G$-equivariant morphism
\[
 f=\varphi_\Lambda\circ\pi\colon Y\longrightarrow\PP^1.
\]
Since $D$ is integral, we have $f_*\OO_Y=\OO_{\PP^1}$. The strict transform $F\subset Y$ of $D$ is a smooth fiber of $f$ and is birational to $T$. Thus $K_F=K_Y|_F$ is not pseudo-effective, so $K_Y$ is not pseudo-effective over $\PP^1$.

The $G$-equivariant $K_Y$-minimal model program over $\PP^1$ gives
\[
 Y\overset{\psi}{\dashrightarrow}Y'\xrightarrow{g}Z\xrightarrow{h}\PP^1,
 \qquad f=h\circ g\circ\psi,
\]
where $g$ is a $G\QQ$-Mori fiber space and $h$ is surjective. Since $Y'$ is birational to $X$, the threefold $Y'$ is rationally connected. Since $\dim Z\ge1$, Proposition~\ref{prop:mori} implies that $G$ is of product type.
\end{proof}

\begin{lemma}\label{lem:surfacecoh}
Let a finite abelian group $G$ act faithfully on a projective terminal rationally connected threefold $Y$, and let $\Lambda\subseteq|-2K_Y|$ be a $G$-invariant pencil. Assume that
\begin{enumerate}[label=(\arabic*)]
\item $-K_Y$ is nef,
\item $2K_Y$ is Cartier,
\item $(-K_Y)^3=0$,
\item $\plur_{-1}(Y)=0$, and
\item a general element $D\in\Lambda$ is integral, normal, Gorenstein, and disjoint from $\Sing_{\nG}Y$.
\end{enumerate}
Write $\B(Y)=mQ_2$. Then the following hold.
\begin{enumerate}[label=(\arabic*)]
\item We have
\begin{equation}\label{eq:Dcohom}
 \plur_1(D)=0,\qquad
 \chi(\OO_D)=4-\frac m4,\qquad q(D)=\frac m4-3.
\end{equation}
\item If $G$ is not of product type, then $m=12$ or $16$, and $D$ has only Du Val singularities. The minimal resolution $\nu\colon T\to D$ satisfies
\begin{equation}\label{eq:Tcohom}
 p_g(T)=0,\qquad K_T\text{ nef},\qquad K_T^2=0,\qquad
 \plur_2(T)=\plur_{-2}(Y)-1.
\end{equation}
For $m=12$, we have $q(T)=0$, whereas for $m=16$, we have $q(T)=1$.
\end{enumerate}
\end{lemma}
\begin{proof}
Adjunction gives $\omega_D\simeq\OO_D(-K_Y)$ and
\[
 0\longrightarrow\omega_Y\longrightarrow\OO_Y(-K_Y)
 \longrightarrow\omega_D\longrightarrow0.
\]
Since $H^1(Y,\omega_Y)=H^2(Y,\OO_Y)^*=0$ and $\plur_{-1}(Y)=0$, the first equality of \eqref{eq:Dcohom} follows. Equations~\eqref{eq:RR}--\eqref{eq:correction} and~\eqref{eq:basketc2}, with $\chi(\OO_Y)=1$, give $\chi(Y,-K_Y)=3-m/4$. Since $D\sim-2K_Y$, we have an exact sequence
\[
 0\longrightarrow\OO_Y(2K_Y)\longrightarrow\OO_Y
 \longrightarrow\OO_D\longrightarrow0.
\]
Taking Euler characteristics and applying Serre duality gives
\[
 \chi(\OO_D)=\chi(\OO_Y)-\chi(\OO_Y(2K_Y))
 =1+\chi(\OO_Y(-K_Y))=4-m/4.
\]

Let $\nu\colon T\to D$ be the minimal resolution. Since $\nu$ is an isomorphism over $D_{\reg}$ and $\omega_D$ is reflexive, there is a natural inclusion $\nu_*\omega_T\hookrightarrow\omega_D$. Hence
\[
 H^0(T,\omega_T)=H^0(D,\nu_*\omega_T)
 \hookrightarrow H^0(D,\omega_D)=0.
\]
Thus $p_g(T)=0$. If $T$ is birationally ruled, Lemma~\ref{lem:ruled} applies. Otherwise the minimal model $T_{\min}$ of $T$ has nef canonical class, and Noether's formula and Hodge decomposition give
\[
 K_{T_{\min}}^2+b_2(T_{\min})+8q(T)=10.
\]
Since $b_2(T_{\min})\ge1$, we obtain $q(T)\le1$. Together with $0\le q(D)\le q(T)$, equation~\eqref{eq:Dcohom} gives $m\in\{12,16\}$.

Since $h^2(\OO_D)=h^2(\OO_T)=0$, the Leray spectral sequence gives
\begin{equation}\label{eq:nonrationaldefect}
 \sum_{P\in\Sing D}\dim(R^1\nu_*\OO_T)_P=q(T)-q(D).
\end{equation}
For $m=16$, we have $q(D)=1$, so $q(D)\le q(T)\le1$ gives $q(T)=q(D)=1$.

Suppose instead that $m=12$ and $q(T)=1$ for a general element $D\in\Lambda$. By \eqref{eq:Dcohom} and \eqref{eq:nonrationaldefect}, every sufficiently general element $D_t\in\Lambda$ has a unique nonrational singular point $P_t$. Let $U\subseteq\Lambda$ be a nonempty open subset, and let $\mathscr D_U\subset Y\times U$ be the incidence family with fiber $D_t\times\{t\}$ over $t\in U$. Choose a projective resolution $\mu\colon\mathscr T_U\to\mathscr D_U$. After shrinking $U$, the restriction $\mu_t\colon T_t\to D_t$ is a resolution for every $t\in U$, and generic base change gives $(R^1\mu_*\OO_{\mathscr T_U})|_{D_t}\simeq R^1(\mu_t)_*\OO_{T_t}$.

Let $Z\subseteq\mathscr D_U$ be the support of $R^1\mu_*\OO_{\mathscr T_U}$ endowed with the reduced induced scheme structure. The projection $Z\to U$ is proper and has one point $(P_t,t)$ over each $t\in U$, so $Z\to U$ is finite. After shrinking $U$, the projection $Z\to U$ is \'etale of degree one and hence an isomorphism. Thus the projection $Z\to Y$ defines a rational map $\sigma\colon\Lambda\dashrightarrow Y$ with $\sigma(t)=P_t$ for general $t$. For $g\in G$ and general $t\in\Lambda$, the isomorphism $g|_{D_t}\colon D_t\to D_{g(t)}$ sends $P_t$ to $P_{g(t)}$ by uniqueness of the nonrational singular point on $D_{g(t)}$. Thus $\sigma$ is $G$-equivariant. The closure of the image of $\sigma$ is a $G$-fixed point or a $G$-invariant rational curve. Thus $G$ is of product type by Proposition~\ref{prop:invariants}.

Thus, if $G$ is not of product type, either $m=12$ and $q(T)=q(D)=0$, or $m=16$ and $q(T)=q(D)=1$. In both cases, \eqref{eq:nonrationaldefect} gives $R^1\nu_*\OO_T=0$. The singularities of $D$ are rational and Gorenstein, hence Du Val by \cite[Corollary~5.24 and Theorem~4.20]{KM}. Thus $K_T=\nu^*K_D$ is nef and $K_T^2=2(-K_Y)^3=0$. Finally, the exact sequence
\[
 0\longrightarrow\OO_Y\longrightarrow\OO_Y(-2K_Y)
 \longrightarrow\OO_D(2K_D)\longrightarrow0
\]
and $H^1(Y,\OO_Y)=0$ give the formula for $\plur_2(T)$ in \eqref{eq:Tcohom}.
\end{proof}

\subsection{Surfaces with \texorpdfstring{$p_g=q=0$}{pg=q=0}}\label{sec:lefschetz}

Let $u\ne\id$ be a finite-order automorphism of a smooth projective surface $T$. Let $\mathcal P$ be the set of isolated points of $\Fix(u)$ and $\mathcal C$ the set of its curve components. The curves $C\in\mathcal C$ are smooth and pairwise disjoint.

For $P\in\mathcal P$, let $\alpha_P,\beta_P\ne1$ be the eigenvalues of $u^*$ on the cotangent space of $T$ at $P$. For $C\in\mathcal C$, let $\lambda_C\ne1$ be the scalar by which $u^*$ acts on the conormal line bundle $N^*_{C/T}$. Put
\begin{equation}\label{eq:localLefschetz}
 a_P=\frac1{(1-\alpha_P)(1-\beta_P)},\qquad
 b_C=\frac{1-g(C)}{1-\lambda_C}
 -\frac{\lambda_C C^2}{(1-\lambda_C)^2},
\end{equation}
where $g(C)=h^1(C,\OO_C)$ and $C^2=\deg N_{C/T}$. The holomorphic Lefschetz formula \cite[Theorem~4.6]{AS}, \cite[Theorem~2.5]{LL} gives
\begin{equation}\label{eq:holomorphicLefschetz}
 \sum_{i=0}^2(-1)^i\Tr\bigl(u^*\mid H^i(T,\OO_T)\bigr)
 =\sum_{P\in\mathcal P}a_P+\sum_{C\in\mathcal C}b_C.
\end{equation}
In particular, if $p_g(T)=q(T)=0$, the right-hand side of \eqref{eq:holomorphicLefschetz} equals one.

\begin{lemma}\label{lem:curveorbits}
Let a finite abelian group $G$ act on a smooth projective surface $T$. Let $C\subset T$ be a smooth irreducible curve and put $r=|G\cdot C|$. If the cardinality of every $G$-orbit of closed points of $T$ is divisible by a positive integer $d$, then
\[
 d\mid r(2-2g(C)),\qquad d\mid rC^2.
\]
\end{lemma}
\begin{proof}
Put $G_C=\{a\in G\colon a(C)=C\}$. Let $\mathscr L$ be a $G_C$-linearized line bundle on $C$. In the space of rational sections of $\mathscr L$, the $\CC$-linear span of the $G_C$-orbit of any nonzero section contains a nonzero eigensection, since $G_C$ is finite abelian. Its divisor is a $G_C$-invariant zero-cycle of degree $\deg\mathscr L$. Summing its $r$ translates gives a $G$-invariant zero-cycle on $T$ of degree $r\deg\mathscr L$, which is divisible by $d$. Apply this to the naturally $G_C$-linearized bundles $T_C$ and $N_{C/T}$, using
\[
 \deg T_C=2-2g(C),\qquad \deg N_{C/T}=C^2.
\]
The proof is completed.
\end{proof}

Recall that $\tau_2(G)$ is the smallest index of a subgroup of $G$ generated by at most two elements, see \eqref{eq:taudefinition} and Lemma~\ref{lem:tau}.

\begin{lemma}\label{lem:elementarysurfacegroup}
Let a finite abelian group $G$ act faithfully on a smooth projective surface $T$ with $p_g(T)=q(T)=0$. Let $(p,k)=(3,1)$ or $(2,2)$. If $p^k\mid\tau_2(G)$, then $G$ is a $p$-elementary group.
\end{lemma}
\begin{proof}
Put $d=p^k$ and $e=d/\gcd(d,2)$. By Corollary~\ref{cor:smoothstabilizers}, every point stabilizer is generated by at most two elements. Lemma~\ref{lem:tau} therefore gives $d\mid |G\cdot P|$ for every $P\in T$. For a smooth irreducible curve $C\subset T$, put $r_C=|G\cdot C|$. Lemma~\ref{lem:curveorbits} gives
\[
 e\mid r_C(1-g(C)),\qquad d\mid r_C C^2.
\]
For $u\in G\setminus\{\id\}$, the group $G$ permutes the components of $\Fix(u)$, and the eigenvalues of $u^*$ on the cotangent spaces at isolated fixed points and the triples $(g(C),C^2,\lambda_C)$ for fixed curves are preserved under the action of $G$. Thus the contributions in \eqref{eq:localLefschetz} agree within each orbit.

Suppose that $u$ has prime order $\ell\ne p$. Put $R=\ZZ[\xi_\ell,1/\ell]$. The denominators in \eqref{eq:localLefschetz} are units in $R$. For each orbit of isolated fixed points, $p$ divides its size since $p\mid d$. For each orbit of fixed curves, $p$ divides both $r_C(1-g(C))$ and $r_CC^2$ since $p\mid e$ and $p\mid d$. Hence the sum of the Lefschetz terms over any such orbit belongs to $pR$, contrary to \eqref{eq:holomorphicLefschetz}, whose left-hand side is $1\notin pR$. Thus $G$ is a $p$-group.

Suppose that $u\in G$ has order $p^2$. Put $\varpi=1-\xi_{p^2}$ and $R_p=\ZZ[\xi_{p^2}]_{(\varpi)}$. Let $v$ be the discrete valuation on $\QQ(\xi_{p^2})$ with valuation ring $R_p$, normalized by $v(\varpi)=1$ and extended by $v(0)=+\infty$. Then
\[
 v(p)=p(p-1),\qquad
 v(e)=\begin{cases}6,&(p,k)=(3,1),\\2,&(p,k)=(2,2).\end{cases}
\]
At an isolated fixed point $P$, at least one cotangent eigenvalue has order $p^2$ by Lemma~\ref{lem:tangent}. Hence
\[
 v\bigl((1-\alpha_P)(1-\beta_P)\bigr)\le p+1,
 \qquad v\bigl(|G\cdot P|a_P\bigr)\ge kp(p-1)-p-1>0.
\]
For a curve component $C$ of $\Fix(u)$, the conormal eigenvalue $\lambda_C$ has order $p^2$, again by Lemma~\ref{lem:tangent}. Consequently
\[
 v\!\left(\frac{r_C(1-g(C))}{1-\lambda_C}\right)\ge v(e)-1>0,
 \qquad
 v\!\left(\frac{r_C\lambda_C C^2}{(1-\lambda_C)^2}\right)\ge kp(p-1)-2>0.
\]
For each $G$-orbit of isolated points or curve components of $\Fix(u)$, let $s$ be the sum of the terms in \eqref{eq:holomorphicLefschetz} indexed by the orbit. We have $v(s)>0$, so $s$ lies in the maximal ideal $\varpi R_p$. Since the left-hand side of \eqref{eq:holomorphicLefschetz} equals $1$, summing over the orbits gives $1\in\varpi R_p$, a contradiction. Thus $G$ has exponent $p$.
\end{proof}

\begin{proposition}\label{prop:C333}
Let $G$ be a finite abelian group acting faithfully on a smooth projective surface $T$ with $p_g(T)=q(T)=0$. If $G$ contains $E\simeq (\ZZ/3\ZZ)^3$, then $G=E$.
\end{proposition}
\begin{proof}
The inclusion $E\subseteq G$ gives $3\mid\tau_2(G)$. Lemma~\ref{lem:elementarysurfacegroup} implies that $G$ is a $3$-elementary group, while \cite[Corollary~9]{KZ} gives $\rk G\le3$. Hence $G=E$.
\end{proof}

\begin{proposition}\label{prop:surfacegroups}
Let a finite abelian group $G$ act faithfully on a smooth projective surface $T$ with $p_g(T)=q(T)=0$. Then either $G$ embeds into $\Cr_2(\CC)$, or
\begin{equation}\label{eq:badA}
 G=\ZZ/2\ZZ\times \ZZ/e\ZZ\times \ZZ/f\ZZ,\qquad
 2\mid e\mid f,\quad e\ge4,\quad(e,f)\ne(4,4).
\end{equation}
\end{proposition}
\begin{proof}
Since $\chi(\OO_T)=1$, \cite[Proposition~8]{KZ} together with Corollary~\ref{cor:smoothstabilizers} and Lemma~\ref{lem:tau} gives $\tau_2(G)\mid12$. If $3\mid\tau_2(G)$, then $(\ZZ/3\ZZ)^3\subseteq G$ and Proposition~\ref{prop:C333} gives $G=(\ZZ/3\ZZ)^3$. If $4\mid\tau_2(G)$, Lemma~\ref{lem:elementarysurfacegroup} implies that $G$ is a $2$-elementary group. The rank bound in \cite[Corollary~9]{KZ} and $4\mid\tau_2(G)$ give $G=(\ZZ/2\ZZ)^4$. Both groups embed into $\Cr_2(\CC)$ by \eqref{eq:A12}.

The remaining values are $\tau_2(G)=1$ and $2$. If $\tau_2(G)=1$, then $\rk G\le2$. If $\tau_2(G)=2$, the invariant factor decomposition of $G$ is $G\simeq\ZZ/2\ZZ\times\ZZ/e\ZZ\times\ZZ/f\ZZ$, with $2\mid e\mid f$. Comparison with \eqref{eq:A12} shows that among these groups only those in \eqref{eq:badA} do not embed into $\Cr_2(\CC)$.
\end{proof}

\section{Genus-one fibrations}\label{sec:elliptic-pencils}
In this section, we analyze the pencil $\Lambda_Y=|-2K_Y|$ on the threefold $Y$ constructed in Proposition~\ref{prop:nefmodel} when $\plur_{-2}(X)=2$ and $\tau_3(G)=2$. Our aim is to prove that $G$ is of product type. 
Suppose that $G$ is not of product type. For a sufficiently general element $D\in\Lambda_Y$, let $T\to D$ be its minimal resolution. Lemma~\ref{lem:surfacecoh} gives $p_g(T)=0$, $K_T$ is nef, $K_T^2=0$, $\plur_2(T)=1$ and $q(T)\in\{0,1\}$. The proof has two cases.
\begin{enumerate}
\item If $q(T)=0$, then $T$ is Enriques or properly elliptic with $\chi(\OO_T)=1$. Let $H$ be the kernel of the action of $G$ on $\Lambda_Y$. Its action on $T$ is faithful, and Corollary~\ref{cor:regular} gives an embedding of $H$ into $\Cr_2(\CC)$. Since $G/H$ acts faithfully on the base $\PP^1$ of the pencil, Proposition~\ref{prop:extensions} implies that $G$ is of product type.
\item If $q(T)=1$, then $D$ is smooth, so $T=D$, and its Albanese morphism has a general fiber of genus one by Proposition~\ref{prop:anticanonical-Albanese}. If $D$ is properly elliptic, we use the action on the point of the Albanese curve corresponding to the unique multiple fiber. If $D$ is bielliptic, the pencil defines a morphism $f\colon Y\to\PP^1$, and we use torsion sections of the Jacobian of the genus-one fibration $b\colon B\to\PP^1$ obtained from the Albanese curves of the general fibers of $f$. Proposition~\ref{thm:irregular} treats both cases.
\end{enumerate}

\subsection{Surfaces with \texorpdfstring{$q=0$}{q=0}}\label{sec:regular}

\begin{lemma}\label{lem:Enriquesgroups}
Let $T$ be an Enriques surface. Then every finite abelian subgroup of $\Aut(T)$ admits an embedding into $\Cr_2(\CC)$.
\end{lemma}
\begin{proof}
Let $G\subseteq\Aut(T)$ be finite abelian, and let $\pi_T\colon\widetilde T\to T$ be the canonical double cover, with deck involution $\iota$, where $\widetilde T$ is a smooth K3 surface. The group $\widetilde G$ of all lifts of elements of $G$ fits into a central exact sequence
\[
 1\longrightarrow\langle\iota\rangle\simeq \ZZ/2\ZZ
 \longrightarrow\widetilde G\longrightarrow G\longrightarrow1.
\]
The deck involution acts as $-1$ on $H^0(\widetilde T,K_{\widetilde T})$, so $\widetilde G$ is abelian by Lemma~\ref{lem:liftcharacter}.

Suppose that $G$ does not embed into $\Cr_2(\CC)$. Proposition~\ref{prop:surfacegroups} gives $|G|\ge64$. Since $G$ is noncyclic, so is $\widetilde G$. Its action on the K3 surface $\widetilde T$ is faithful. The list of maximal abelian groups in \cite[Theorem~3.7]{LPZ}, obtained from the classification in \cite{BrandhorstHofmann}, gives 
$
 128\le2|G|=|\widetilde G|\le72.
$
This contradiction proves the lemma.
\end{proof}

\smallskip\noindent\emph{Properly elliptic surfaces.}
A smooth minimal projective surface $T$ is called \emph{properly elliptic} if $\kappa(T)=1$. By the classification of surfaces, for a sufficiently large divisible integer $m$, the linear system $|mK_T|$ defines the canonical genus-one fibration $h\colon T\to B$ over a smooth projective curve $B$ \cite[\S E.7.3]{Reid}.

Let $E$ be a genus-one curve over an algebraically closed field of characteristic zero with a chosen origin $o\in E$. Write $t_a$ for translation by $a\in E$, and let $\Aut(E,o)$ denote the subgroup of automorphisms fixing $o$. For $u\in\Aut(E)$, put $\ell(u)=t_{-u(o)}\circ u$. The exact sequence
\[
 1\longrightarrow E\xrightarrow{\,a\mapsto t_a\,}\Aut(E)
 \xrightarrow{\ell}\Aut(E,o)\longrightarrow1
\]
is split by the inclusion of $\Aut(E,o)$. Thus $\Aut(E)=E\rtimes\Aut(E,o)$, with $\ell$ the projection onto the second factor.

\begin{lemma}[{cf.~\cite[Proposition~4.2]{Log}}]\label{lem:ellipticgroup}
Let $G\subseteq\Aut(E)$ be a finite abelian subgroup. If $\ell(G)\ne\{1\}$, then
$
 |G|\le9$ and $|G_2|\le8.
$
If $|\ell(G)|$ is $2,3,4,6$, then $|G|\le8,9,8,6$, respectively. If $|\ell(G)|=2$, then $G[2]=G$.
\end{lemma}

Let $f\colon J\to C$ be a relatively minimal elliptic surface with zero section $O\colon C\to J$. The sections of $f$ form the Mordell--Weil group $\MW(J)$; write $\MW(J)_{\tors}$ for its torsion subgroup. For a complex projective variety $Z$, let $e(Z)$ be its topological Euler characteristic. Then $e(J)=\sum_{c\in C}e(f^{-1}(c))$, and a fiber of Kodaira type $I_n$, $n\ge1$, has Euler characteristic $n$ (see \cite[Theorem~6.10]{SS}). At an additive fiber, that is, a fiber of Kodaira type $\mathrm{II}$, $\mathrm{III}$, $\mathrm{IV}$, $\mathrm{I}_n^*$ with $n\ge0$, $\mathrm{IV}^*$, $\mathrm{III}^*$, or $\mathrm{II}^*$, specialization injects $\MW(J)_{\tors}$ into the component group of that fiber, which has order at most four by \cite[Lemmas~7.3 and~7.8]{SS}. For a rational elliptic surface $J\to\PP^1$, one has $e(J)=12$ by \cite[Corollary~6.11]{SS}, and \cite[Corollary~11.8]{SS} gives
\begin{equation}\label{eq:MWbound}
 |\MW(J)_{\tors}|\le9.
\end{equation}
Let $B\to\PP^1$ be a genus-one fibration on a rational surface, and let $J_B\to\PP^1$ be the Jacobian of $B\to\PP^1$. Then $\chi(\OO_{J_B})=\chi(\OO_B)=1$, and $J_B$ has no multiple fibers. Hence $K_{J_B}\equiv-F$, where $F$ is a fiber of $J_B\to\PP^1$, so $J_B$ is rational (see \cite[\S2.4]{CLS} and \cite[Theorem~6.8]{SS}).

\begin{proposition}\label{prop:translationbase}
Let a finite abelian group $G$ act faithfully on a smooth projective surface $B$ and preserve a relatively minimal genus-one fibration $f\colon B\to\PP^1$ whose Jacobian is a rational elliptic surface. If a subgroup $T\subseteq G$ with $|T|>4$ acts on $B$ fiberwise by translations, then $G$ acts trivially on~$\PP^1$.
\end{proposition}
\begin{proof}
Let $j\colon J\to\PP^1$ be the Jacobian of $f$, with zero section $O\colon\PP^1\to J$. Put $K=\CC(\PP^1)$ and let $J_\eta$ be the generic fiber of $j$. Each $t\in T$ acts on the generic fiber of $f$ as translation by $P_t\in J_\eta(K)_{\tors}$. Denote the corresponding section by $s_t\colon\PP^1\to J$. Identify $T$ with the subgroup $\{s_t\mid t\in T\}$ of $\MW(J)_{\tors}$. Fix $g\in G$. Let $u\in\Aut(J)$ and $\bar u\in\Aut(\PP^1)$ be the automorphisms induced by $g$ on $J$ and $\PP^1$, respectively. By \cite[Theorem~3.3(1) and~(3)]{DolgachevMartin}, we have $j\circ u=\bar u\circ j$ and $u\circ O=O\circ\bar u$. Since $gt=tg$ for each $t\in T$, we have $u\circ s_t=s_t\circ\bar u$.

In short Weierstrass coordinates, an elliptic curve automorphism fixing the origin has the form $(x,y)\mapsto(\zeta^2x,\zeta^3y)$ for a root of unity $\zeta$, see \cite[proof of Theorem~III.10.1]{Silverman}. If $\zeta\ne1$, either $\zeta^2\ne1$ or $\zeta^3\ne1$, so all affine fixed points lie on one of the axes $x=0$ and $y=0$. Each axis meets the Weierstrass cubic in at most three points. Thus a non-trivial origin-preserving automorphism fixes at most four points, including the origin. Put $d=\ord(\bar u)$. Let $J_{\bar\eta}=J_\eta\times_K\overline K$, where $\overline K$ is an algebraic closure of $K$. The automorphism $u^d$ fixes the $|T|>4$ distinct points $P_t$ of $J_{\bar\eta}$, including the origin. Hence $u^d|_{J_{\bar\eta}}=\id$, so $u^d=\id_J$. Thus $\ord(u)=d$.

Let $\rho\colon\PP^1\to C_0=\PP^1/\langle\bar u\rangle$ be the quotient morphism of degree $d$, and let $j_0\colon J_0\to C_0$ be the relatively minimal model of the quotient fibration $J/\langle u\rangle\to C_0$. The induced rational map $q_J\colon J\dashrightarrow J_0$ fits into the commutative diagram
\[
\begin{tikzcd}[column sep=large,row sep=large]
 J \arrow[r,dashed,"q_J"] \arrow[d,"j"']
 & J_0 \arrow[d,"j_0"] \\
 \PP^1 \arrow[r,"\rho","d\,{:}\,1"']
 & C_0
\end{tikzcd}
\]
Since $\ord(u)=\ord(\bar u)=d$, the $\langle u\rangle$-orbit of a general fiber $j^{-1}(p)$ consists of $d$ distinct fibers. Thus $q_J$ restricts to an isomorphism $j^{-1}(p)\to j_0^{-1}(\rho(p))$. Hence $J$ is the relatively minimal model of $J_0\times_{C_0}\PP^1$.

For $t\in T$, the identity $u\circ s_t=s_t\circ\bar u$ makes $q_J\circ s_t$ invariant under $\langle\bar u\rangle$. Thus $q_J\circ s_t=s_{t,0}\circ\rho$ for a rational section $s_{t,0}$ of $j_0$, which extends over $C_0$ by properness. Applying the same construction to $O$ gives the zero section $O_0$ of $j_0$. On general fibers $q_J$ induces isomorphisms sending $O(p)$ to $O_0(\rho(p))$, so $t\mapsto s_{t,0}$ identifies $T$ with $T_0=\{s_{t,0}\mid t\in T\}\subseteq\MW(J_0)_{\tors}$.

At an additive fiber of $j_0$, specialization would embed $T_0$ into the component group of that fiber, which has order at most four \cite[Lemmas~7.3 and~7.8]{SS}. Since $|T_0|=|T|>4$, every singular fiber of $j_0$ has Kodaira type $I_n$ for some $n\ge1$. For $p\in\PP^1$, let $a_p$ be the ramification index of $\rho$ at $p$. If $j_0^{-1}(\rho(p))$ has type $I_n$, then $j^{-1}(p)$ has type $I_{a_pn}$, and if $j_0^{-1}(\rho(p))$ is smooth, then $j^{-1}(p)$ is smooth \cite[\S5.2]{SS}. Since $\sum_{\rho(p)=q}a_p=d$ for every $q\in C_0$, summing the Euler characteristics of the fibers gives $e(J)=d\,e(J_0)$ \cite[Theorem~6.10]{SS}. As $J$ is rational and $e=12\chi(\OO)$ for elliptic surfaces with section \cite[Corollary~6.11]{SS}, we obtain
\[
 1=\chi(\OO_J)=d\,\chi(\OO_{J_0}).
\]
Since $\chi(\OO_{J_0})$ is an integer, we have $d=1$. Thus each $g\in G$ acts trivially on $\PP^1$.
\end{proof}

\begin{lemma}\label{lem:properlyellipticgroups}
Let $T$ be a properly elliptic surface with $p_g(T)=q(T)=0$ and $\plur_2(T)\le6$. Then every finite abelian subgroup of $\Aut(T)$ admits an embedding into $\Cr_2(\CC)$.
\end{lemma}
\begin{proof}
Let $G\subseteq\Aut(T)$ be finite abelian, and suppose that $G$ does not embed into $\Cr_2(\CC)$. Proposition~\ref{prop:surfacegroups} gives
\begin{equation}\label{eq:properA}
 G=\ZZ/2\ZZ\times \ZZ/e\ZZ\times \ZZ/f\ZZ,\qquad
 2\mid e\mid f,\quad e\ge4,\quad(e,f)\ne(4,4).
\end{equation}
The canonical genus-one fibration $h\colon T\to B$ is preserved by every automorphism. Since $q(T)=0$, we have $B=\PP^1$. Let $p_1,\ldots,p_s\in B$ be the points over which $h$ has multiple fibers. Write $h^*p_i=m_iF_i'$, where $m_i$ is the greatest common divisor of the component multiplicities of $h^*p_i$. Let $F$ be a general fiber of $h$. Since $\chi(\OO_T)=1$, the canonical bundle formula \cite[\S E.7.3]{Reid} gives
\[
 K_T\sim-F+\sum_i(m_i-1)F_i',\qquad
 2K_T\sim(s-2)F+\sum_i(m_i-2)F_i'.
\]
Pushing forward $\OO_T(2K_T)$ along $h$ gives $h_*\OO_T(2K_T)\simeq\OO_{\PP^1}(s-2)$. Since $\kappa(T)=1$, we have $s\ge2$. Thus $s=\plur_2(T)+1\le7$. This computation is also made in \cite[proof of Proposition~5.1(3)(a)]{CLS}.

Let $j\colon J\to\PP^1$ be the relatively minimal Jacobian of $h$. We have $\chi(\OO_J)=\chi(\OO_T)=1$ by \cite[\S2.4]{CLS}, and the zero section of $j$ excludes multiple fibers. The canonical bundle formula \cite[Theorem~6.8]{SS} gives $K_J\sim-F_J$, where $F_J$ is a general fiber of $j$. Hence $\kappa(J)=-\infty$, and $\chi(\OO_J)=1$ implies that $J$ is rational.

Set $N=\ker(G\to\Aut(B))$ and $Q=G/N$. Put $K=\CC(B)$ and fix an algebraic closure $\overline K$ of $K$. Let $\eta=\Spec K$ be the generic point of $B$, and write $F_\eta=T\times_B\eta$. Choose an origin on the geometric generic fiber $\overline F_\eta=F_\eta\times_K\overline K$. The action of $N$ on $\overline F_\eta$ is faithful. If this action consists of translations, then $N$ embeds into $\MW(J)_{\tors}$, so \eqref{eq:MWbound} applies. Otherwise apply Lemma~\ref{lem:ellipticgroup}. In either case, we have
\begin{equation}\label{eq:properN}
 |N|\le9.
\end{equation}

By \eqref{eq:properA} and \eqref{eq:properN}, we have $|Q|\ge64/9>7$, so the abelian subgroup $Q\subseteq\Aut(\PP^1)$ is cyclic of order at least eight. It has two fixed points, and every other orbit has length $|Q|$. Since $Q$ preserves the $s\le7$ points of $B$ corresponding to the multiple fibers of $h$, we obtain $s=2$. The order of $Q$ divides $f$, hence $2e\mid|N|$. Since $e\ge4$ and $|N|\le9$, we obtain $e=4$ and $|N|=8$. If $N$ acted by translations, Proposition~\ref{prop:translationbase} would give $Q=1$ and $G=N$, contrary to \eqref{eq:properA}.

Moreover, $2N\ne0$: otherwise $2G$ would be a quotient of the cyclic group $G/N$, whereas \eqref{eq:properA} shows that $2G$ is noncyclic. For the action of $N$ on $\overline F_\eta$, Lemma~\ref{lem:ellipticgroup} therefore gives $|\ell(N)|=4$. Choose $u\in N$ such that $\ell(u)$ has order four. With respect to the chosen origin on $\overline F_\eta$, write $u=t_b\circ\ell(u)$. The fixed scheme $Z_\eta=\Fix_{F_\eta}(u)$ becomes the fiber $(1-\ell(u))^{-1}(b)$ over $\overline K$. Hence $\deg Z_\eta=\deg(1-\ell(u))=2$. The closure of $Z_\eta$ on $T$ is an effective horizontal divisor $Z$ with $Z\cdot F=2$. Thus
$
 2=Z\cdot(m_iF_i')=m_i(Z\cdot F_i'),
$
so $m_1=m_2=2$. The canonical bundle formula gives $2K_T\sim0$, contrary to $\kappa(T)=1$. This completes the proof.
\end{proof}

\begin{corollary}\label{cor:regular}
Let $T$ be a smooth projective surface with
\begin{equation}\label{eq:regularsurface}
 p_g(T)=q(T)=0,\qquad\quad K_T\text{ nef},\qquad\quad K_T^2=0,\qquad\quad \plur_2(T)\le6.
\end{equation}
Then every finite abelian subgroup of $\Aut(T)$ admits an embedding into $\Cr_2(\CC)$.
\end{corollary}
\begin{proof}
Since $K_T$ is nef, the surface $T$ is minimal. By the classification of surfaces \cite[\S\S E.1, E.4, and E.7.3]{Reid}, the surface $T$ is Enriques or properly elliptic. The assertion follows from Lemmas~\ref{lem:Enriquesgroups} and~\ref{lem:properlyellipticgroups}.
\end{proof}

\begin{corollary}\label{cor:regularpencil}
Let a finite abelian group $G$ act faithfully on a normal projective threefold $Y$ and preserve a pencil $\Lambda$ with integral normal general element $D$. Suppose that the minimal resolution $T\to D$ satisfies \eqref{eq:regularsurface}.
Then $G$ is of product type.
\end{corollary}
\begin{proof}
Let $H$ be the kernel of the action of $G$ on $\Lambda$. Since $G$ is finite and acts faithfully on $Y$, the induced action of $H$ on a sufficiently general element $D\in\Lambda$ is faithful. It lifts to the minimal resolution $T\to D$, so Corollary~\ref{cor:regular} gives an embedding of $H$ into $\Cr_2(\CC)$. The quotient $G/H$ acts faithfully on the base $\PP^1$ of the pencil. Proposition~\ref{prop:extensions} therefore implies that $G$ is of product type.
\end{proof}

\subsection{Surfaces with \texorpdfstring{$q=1$}{q=1}}\label{sec:anticanonical-Albanese}

A \emph{bielliptic surface} is a smooth minimal projective surface $T$ with $\kappa(T)=0$, $p_g(T)=0$, and $q(T)=1$. By \cite[Definition~VI.19 and List~VI.20]{BeauvilleSurfaces}, there is an isomorphism
\begin{equation}\label{eq:bielliptic-quotient}
 T\simeq(E\times F)/\Gamma,
\end{equation}
where $E,F$ are elliptic curves and $\Gamma$ is a finite group acting diagonally, faithfully by translations on $E$, and by automorphisms on $F$, with $F/\Gamma\simeq\PP^1$. The projection $\operatorname{pr}_E\colon E\times F\to E$ induces an Albanese morphism $a_T\colon T\to E/\Gamma$, whose fibers are isomorphic to $F$.

We use the construction and universal property of the Albanese variety from \cite[Theorem~V.13]{BeauvilleSurfaces}. For a smooth projective surface with one-dimensional Albanese image, the Albanese morphism has connected fibers \cite[Proposition~V.15]{BeauvilleSurfaces}.

\begin{lemma}[{cf.~\cite[Corollary~4.1.2 and Theorem~4.2.6]{DeligneHodge}}]\label{lem:basecurve-Albanese}
Let $W$ be a normal projective variety with a resolution $\widetilde W\to W$ satisfying $q(\widetilde W)=0$, and let $\Lambda$ be a pencil of Cartier divisors with smooth connected general element $D$. Let $a_D\colon D\to\Alb(D)$ be an Albanese morphism. Let $C\subseteq\Bs\Lambda$ be an integral curve with normalization $\nu\colon\widetilde C\to C$. Then the restriction
\[
 (a_D|_C)\circ\nu\colon\widetilde C\longrightarrow\Alb(D)
\]
is constant.
\end{lemma}
\begin{proof}
Put $\PP^1=\Lambda$, and choose a nonempty open subset $U\subseteq\PP^1$ over which the incidence family $\mathscr D_U=\{(x,t)\in W\times U\mid x\in D_t\}$ is smooth over $U$. Let $p\colon\mathscr D_U\to U$ be the projection $p(x,t)=t$, and let $t_0\in U$ be the parameter with $D=D_{t_0}$. Pullback along $\gamma\colon\widetilde C\times U\to\mathscr D_U$, $\gamma(c,t)=(\nu(c),t)$, induces the following morphism of local systems on $U$
\[
 \varphi\colon R^1p_*\QQ\longrightarrow H^1(\widetilde C,\QQ)_U,
\]
where the target is the constant local system with fiber $H^1(\widetilde C,\QQ)$. Thus $\mathscr I=\operatorname{im}\varphi$ is a constant local system. By semisimplicity \cite[Theorem~4.2.6]{DeligneHodge}, the exact sequence
\[
 0\longrightarrow\ker\varphi\longrightarrow R^1p_*\QQ
 \xrightarrow{\varphi}\mathscr I\longrightarrow0
\]
splits in the category of local systems. Choose a splitting $\sigma\colon\mathscr I\to R^1p_*\QQ$, so that $\varphi\circ\sigma=\id_{\mathscr I}$. Its image
\[
 \mathscr L=\sigma(\mathscr I)\subseteq R^1p_*\QQ
\]
is a constant sub-local system isomorphic to $\mathscr I$. Choose a smooth projective compactification $\overline{\mathscr D}$ of $\mathscr D_U$ on which $p$ extends to a morphism $\bar p\colon\overline{\mathscr D}\to\PP^1$. The global invariant cycle theorem applied to $\bar p$ identifies $H^1(D,\QQ)^{\pi_1(U,t_0)}$ with the image of $H^1(\overline{\mathscr D},\QQ)\to H^1(D,\QQ)$ \cite[Corollary~4.1.2]{DeligneHodge}. The projection $\mathscr D_U\to W$ is birational because a general point of $W$ belongs to a unique element of $\Lambda$. Thus $\overline{\mathscr D}$ is birational to $\widetilde W$. By birational invariance of irregularity and $q(\widetilde W)=0$, we have $H^1(\overline{\mathscr D},\QQ)=0$. Hence $\mathscr L=0$ and $\mathscr I=0$. Thus the fiber map $\varphi_{t_0}\colon H^1(D,\QQ)\to H^1(\widetilde C,\QQ)$ is zero, so $(a_D|_C)\circ\nu$ is constant.
\end{proof}

\begin{proposition}\label{prop:anticanonical-Albanese}
Let $Y$ be a projective terminal threefold, let $m\ge2$ be an integer, and let $D$ be a general element of a pencil $\Lambda\subseteq|-mK_Y|$. Assume that
\begin{enumerate}[label=(\arabic*)]
\item $H^1(Y,\OO_Y)=0$, $-K_Y$ is nef, $(-K_Y)^3=0$, and $mK_Y$ is Cartier,
\item $D$ is integral and normal with only Du Val singularities, and
\item the minimal resolution $T\to D$ satisfies $p_g(T)=0$ and $q(T)=1$.
\end{enumerate}
Then $D$ is smooth and minimal, and a general fiber of its Albanese morphism is an elliptic curve.
\end{proposition}
\begin{proof}
Let $\nu\colon T\to D$ be the minimal resolution. Since $D$ has only Du Val singularities, we have $K_T=\nu^*K_D$. By adjunction,
\[
 K_D\sim_{\QQ}-(m-1)K_Y|_D,\qquad
 K_T^2=m(m-1)^2(-K_Y)^3=0.
\]
Thus $K_T$ is nef and $\chi(\OO_T)=0$, so $T$ is bielliptic or properly elliptic. In either case, there is a finite \'etale cover
\[
 \rho\colon C_1\times C_2\longrightarrow T,
 \qquad g(C_j)\ge1,\quad j=1,2,
\]
where $C_1,C_2$ are smooth projective curves, see \cite[Definition~VI.19 and List~VI.20]{BeauvilleSurfaces} and \cite[Lemma~1.1]{CFGLS}. Let $\phi\colon\PP^1\to T$ be a morphism. Since $\PP^1$ is simply connected, the finite \'etale projection $\PP^1\times_T(C_1\times C_2)\to\PP^1$ admits a section. Thus there is a morphism $\widetilde\phi\colon\PP^1\to C_1\times C_2$ with $\rho\circ\widetilde\phi=\phi$. For $j=1,2$, the composition $\operatorname{pr}_{C_j}\circ\widetilde\phi\colon\PP^1\to C_j$ is constant because $g(C_j)\ge1$. Hence $\phi$ is constant, and $T$ contains no rational curves. Every curve contracted by $\nu$ would be rational, so $\nu$ is an isomorphism. Since the Cartier divisor $D$ is smooth, we have $D\subseteq Y_{\reg}$.

Write $a_D\colon D\to C_D=\Alb(D)$ for the Albanese morphism. Since $q(D)=1$, the Albanese variety $C_D$ is an elliptic curve. If $D$ is bielliptic, the fibers of $a_D$ are smooth curves of genus one by \eqref{eq:bielliptic-quotient}. Suppose that $D$ is properly elliptic, and let $h_D\colon D\to B_D$ be its canonical genus-one fibration. The reduced fibers of $h_D$ are smooth curves of genus one by \cite[Lemma~1.1]{CFGLS}. For a general element $D'\in\Lambda$ distinct from $D$, adjunction gives
$
 (m-1)(D'|_D)\sim mK_D.
$
Since $\kappa(D)=1$, the effective divisor $D'|_D$ is nonzero. Let $F$ be a general fiber of $h_D$. Adjunction on the smooth elliptic fiber $F$ gives $K_D\cdot F=0$. Intersecting $(m-1)(D'|_D)\sim mK_D$ with $F$ gives $(D'|_D)\cdot F=0$. Thus $D'|_D$ is supported on fibers of $h_D$. Choose an irreducible component $C_0$ of $D'|_D$. Then $C_0$ is a smooth elliptic curve contained in $\Bs\Lambda$.

Terminal singularities are rational, so a resolution $\widetilde Y\to Y$ has $q(\widetilde Y)=h^1(Y,\OO_Y)=0$. Lemma~\ref{lem:basecurve-Albanese} shows that $a_D|_{C_0}$ is constant. Let $r$ be the multiplicity of the fiber of $h_D$ supported on $C_0$. Then
$
 rC_0\equiv F$ and $ r\in\ZZ_{>0}.
$
For a point $p\in C_D$, we have $a_D^*p\cdot F=r(a_D^*p\cdot C_0)=0$. Thus $a_D|_F$ is constant. Since $h_D$ has connected fibers, $a_D$ factors through a rational map $B_D\dashrightarrow C_D$. This map extends to a morphism $\gamma_D\colon B_D\to C_D$ with $a_D=\gamma_D\circ h_D$. Connectedness of the fibers of $a_D$ implies that $\gamma_D$ is an isomorphism. We obtain the commutative diagram
\[
\begin{tikzcd}[column sep=large,row sep=large]
 D \arrow[r,"h_D"] \arrow[dr,"a_D"']
 & B_D \arrow[d,"\gamma_D"',"{\rotatebox[origin=c]{90}{$\scriptstyle\sim$}}"{inner sep=2pt}] \\
 & C_D
\end{tikzcd}
\]
Thus a general fiber of $a_D$ is elliptic, and the proof is completed.
\end{proof}

\begin{lemma}\label{lem:Albanese-fiber-group}
Let $D$ be a smooth minimal projective surface with $p_g(D)=0$ and $q(D)=1$, and suppose that its Albanese morphism $a_D\colon D\to C_D=\Alb(D)$ has general fiber of genus one. Let $N\subseteq\Aut(D)$ be a finite abelian subgroup acting trivially on $C_D$. Then $|N|\le9$.
\end{lemma}
\begin{proof}
Put $K=\CC(C_D)$, fix an algebraic closure $\overline K$ of $K$, and let $F_D=D\times_{C_D}\Spec\overline K$ be the geometric generic fiber of $a_D$. The group $N$ acts faithfully on $F_D$. Since $\chi(\OO_D)=0$, the fibration $a_D$ is isotrivial, see \cite[Lemma~1.1]{CFGLS} and \eqref{eq:bielliptic-quotient}. If some element of $N$ acts on $F_D$ by an automorphism that is not a translation, then $|N|\le9$ by Lemma~\ref{lem:ellipticgroup}. We may therefore assume that every element of $N$ acts on $F_D$ by translation.

Let $j_D\colon J_D\to C_D$ be the relative Jacobian of $a_D$. Write $J_{D,\eta}$ for its generic fiber. Each $n\in N$ acts on $F_D$ as translation by $P_n\in J_{D,\eta}(K)_{\tors}$. Denote the corresponding section by $\sigma_n\colon C_D\to J_D$. Let $V\subseteq C_D$ be the open subset over which $a_D$ is smooth, choose $v_0\in V$, and put $E_0=j_D^{-1}(v_0)$. Write
\[
 \rho\colon\pi_1(V,v_0)\longrightarrow\Aut(E_0,0)
\]
for the monodromy representation of $j_D$ over $V$. Its image is nontrivial. Indeed, if $\rho$ were trivial, the product-quotient descriptions in \cite[Lemma~1.1]{CFGLS} and \eqref{eq:bielliptic-quotient} would give a form $\omega_D\in H^0(D,\Omega_D^1)$ with $\omega_D|_{F_D}\ne0$. Choose $0\ne\omega_C\in H^0(C_D,\Omega_{C_D}^1)$. Since $(a_D^*\omega_C)|_{F_D}=0$, the forms $\omega_D$ and $a_D^*\omega_C$ are linearly independent in $H^0(D,\Omega_D^1)$, contrary to $h^0(D,\Omega_D^1)=q(D)=1$.

Let $p_V\colon V'\to V$ be the connected finite \'etale cover corresponding to $\ker\rho$. Choose an isomorphism preserving the zero sections
\[
 \theta\colon J_D\times_{C_D}V'\xrightarrow{\sim}E_0\times V'
\]
over $V'$. Under $\theta$, the pullbacks of the sections $\sigma_n$ define an injective homomorphism $\iota_N\colon N\hookrightarrow E_0(\CC)_{\tors}$. Choose $\alpha\in\operatorname{im}\rho\setminus\{1\}$. Since each $\sigma_n$ is defined over $V$, the automorphism $\alpha$ fixes every point of $\iota_N(N)$. Hence $\iota_N(N)\subseteq\ker(1-\alpha)$ and $|N|\le\deg(1-\alpha)\le4$, which proves the assertion.
\end{proof}

\begin{lemma}\label{lem:bielliptic32}
Let $D$ be a bielliptic surface with $2K_D\sim0$ and $K_D\not\sim0$, and let $a_D\colon D\to C_D=\Alb(D)$ be its Albanese morphism. Let $G\subset\Aut(D)$ be finite abelian. If some $g\in G$ induces an automorphism of $C_D$ that is not a translation, then $|G|\le32$.
\end{lemma}
\begin{proof}
By \cite[Definition~VI.19 and List~VI.20]{BeauvilleSurfaces}, write
\[
 D\simeq(E\times F)/\Gamma,\qquad
 \Gamma\simeq\ZZ/2\ZZ\text{ or }(\ZZ/2\ZZ)^2,
\]
where $\Gamma$ acts faithfully on both elliptic curves and by translations on $E$. Let $\Gamma_F\subseteq\Aut(F)$ be its image on $F$, and write $q\colon E\times F\to D$ and $q_F\colon F\to F/\Gamma_F$ for the quotient morphisms. By \cite[Corollaries~1.2 and~4.7]{Martin}, each $u\in\Aut(D)$ has a lift $\widetilde u=(u_E,u_F)$ with $q\circ\widetilde u=u\circ q$, and translations on $E$ identify $E$ with $\Aut^0(D)$.

Let $\rho_D\colon\Aut(D)\to\Aut(C_D)$ be the homomorphism determined by $\rho_D(u)\circ a_D=a_D\circ u$ for each $u\in\Aut(D)$. For $e\in E$, write $t_e$ for translation by $e$ on $E$ and for the induced automorphism of $D$. Let $p_E\colon E\to C_D$ be the isogeny determined by $\rho_D(t_e)=t_{p_E(e)}$ for $e\in E$. Put $\lambda=\ell\circ\rho_D$ and $r=|\lambda(G)|$. If $u\in\ker\lambda$, then $u_E=t_{e_u}$ for some $e_u\in E$, since $u_E$ induces a translation on $C_D$ via $p_E$. Conjugation by $\widetilde u$ preserves $\Gamma$ and fixes the action of $\Gamma$ on $E$. Since $\Gamma$ acts faithfully on $E$, the automorphism $u_F$ commutes with $\Gamma_F$.

Define $\psi\colon\ker\lambda\to\Aut(F/\Gamma_F)$ by $\psi(u)\circ q_F=q_F\circ u_F$. Any two lifts of $u$ differ by an element of $\Gamma$, whose component on $F$ belongs to $\Gamma_F$. Thus $\psi$ is a well-defined homomorphism. If $\psi(u)=\id$, then $u_F\in\Gamma_F$. Choose $\gamma=(t_{e_\gamma},u_F)\in\Gamma$. The lift $\widetilde u\gamma^{-1}=(t_{e_u-e_\gamma},\id_F)$ shows that $u\in E$. Conversely, every translation of $E$ alone induces the identity on $F/\Gamma_F$, so $\ker\psi=E$. The finite abelian group $\psi(G\cap\ker\lambda)$ preserves the four branch points of $q_F$. Since $F/\Gamma_F\simeq\PP^1$, this group has order at most four. Thus
\[
 [G\colon G\cap E]=r\,|\psi(G\cap\ker\lambda)|\le4r.
\]

Let $\gamma_0\in\Gamma$ be the unique non-trivial element if $|\Gamma|=2$, and the element whose image in $\Gamma_F$ is the unique non-trivial translation if $|\Gamma|=4$. On $E$, the element $\gamma_0$ acts by $t_a$ for some $a\in E[2]\setminus\{0\}$. For $u\in\Aut(D)$, conjugation by $\widetilde u=(u_E,u_F)$ fixes $\gamma_0$, since $u_F$ normalizes $\Gamma_F$ and preserves its subgroup of translations under conjugation. Thus $u_Et_a u_E^{-1}=t_a$, so $\ell(u_E)(a)=a$. The origin-preserving automorphisms of $E$ fixing $a$ form a group of order dividing four, and $p_E\circ\ell(u_E)=\lambda(u)\circ p_E$. Since $\lambda(G)\ne\{1\}$ by hypothesis, we obtain $r\in\{2,4\}$.

Choose $g\in G$ such that $\lambda(g)$ has order $r$, and let $\alpha\colon E\to E$ be conjugation by $g$ on $E=\Aut^0(D)$. Then $p_E\circ\alpha=\lambda(g)\circ p_E$. Since $p_E$ is surjective with finite kernel and $E$ is connected, $\alpha$ has order $r$. Since $G$ is abelian, $G\cap E\subseteq\ker(1-\alpha)$. Moreover, $\deg(1-\alpha)$ is four for $r=2$ and two for $r=4$. Therefore
\[
 |G|=[G\colon G\cap E]\,|G\cap E|
 \le4r\,\deg(1-\alpha)=32.
\]
This completes the proof.
\end{proof}

\subsection{The relative Albanese fibration}\label{sec:irregular}

Let $V$ be a smooth projective geometrically integral variety over a field $K$. Its Albanese torsor $\Alb^1(V)$ is a torsor under $\Alb(V)$, with the action defined over $K$ and free transitive over an algebraic closure $\overline K$. Write $a_V\colon V\to\Alb^1(V)$ for the Albanese morphism. For every $K$-morphism $f\colon V\to P$ to a torsor under an abelian variety over $K$, there is a unique $K$-morphism $g\colon\Alb^1(V)\to P$ such that $f=g\circ a_V$ \cite[\S2]{Wittenberg}. If $v\in V(K)$, taking $a_V(v)$ as the origin identifies $\Alb^1(V)$ with $\Alb(V)$.

\begin{lemma}\label{lem:relative-Albanese-surface}
Let $G$ be a finite group and let $f\colon Y\to\PP^1$ be a surjective $G$-equivariant morphism from a normal projective rationally connected threefold. Suppose that a general fiber of $f$ is smooth and connected with irregularity one. Then, for sufficiently general $t\in\PP^1$, there is a commutative diagram
\[
\begin{tikzcd}[column sep=large,row sep=large]
 D_t \arrow[r,"a_t"] \arrow[d,hook]
 & B_t \arrow[d,hook] \\
 Y \arrow[r,dashed,"a"] \arrow[dr,"f"']
 & B \arrow[d,"b"] \\
 & \PP^1
\end{tikzcd}
\]
where
\begin{itemize}
\item $D_t=f^{-1}(t)$ and $B_t=b^{-1}(t)$;
\item $B$ is a smooth rational $G$-surface, and $b$ is a relatively minimal genus-one fibration with Jacobian fibration $j\colon J\to\PP^1$, where $J$ is a rational elliptic surface;
\item $a$ is a dominant $G$-equivariant rational map, $b$ is $G$-equivariant, and $a_t$ is an Albanese morphism extending $a|_{D_t}$.
\end{itemize}
\end{lemma}
\begin{proof}
Choose a nonempty $G$-invariant open subset $U\subseteq\PP^1$ over which $f$ is smooth with geometrically connected fibers of irregularity one. Put $Y_U=f^{-1}(U)$. By \cite[Theorem~3.3(iii), \S3.5, and Corollary~3.9]{GroPicard}, the relative Albanese torsor $b_U\colon B_U=\Alb^1(Y_U/U)\to U$ is a smooth proper family of genus-one curves. Formation of the relative Albanese morphism $a_U\colon Y_U\to B_U$ commutes with base change, so $a_t=a_U|_{D_t}$ is the Albanese morphism of $D_t$ for every $t\in U$. Functoriality gives a $G$-action on $B_U$ making $a_U$ and $b_U$ $G$-equivariant.

Choose a smooth projective surface $B$ and a relatively minimal genus-one fibration $b\colon B\to\PP^1$ such that the generic fibers of $b$ and $b_U$ are isomorphic over $\CC(\PP^1)$. After shrinking $U$ to a $G$-invariant open subset, identify $b^{-1}(U)$ with $B_U$ over $U$. For $g\in G$, let $g_B\colon B\dashrightarrow B$ be the birational map extending the action of $g$ on $B_U$, and let $\bar g$ be the automorphism of $\PP^1$ given by $g$. Then $b\circ g_B=\bar g\circ b$. By uniqueness of relatively minimal smooth models of genus-one fibrations, $g_B$ is an automorphism of $B$. The automorphisms $g_B$ define the $G$-action on $B$. The morphism $a_U$ defines a dominant $G$-equivariant rational map $a\colon Y\dashrightarrow B$. Thus $B$ is rationally connected and hence rational. Write
\begin{equation}\label{eq:JacB}
 j\colon J\longrightarrow\PP^1
\end{equation}
for the Jacobian fibration of $b$. Since $B$ is rational, $J$ is a rational elliptic surface, and $\MW(J)_{\tors}$ satisfies \eqref{eq:MWbound}.
This completes the proof.
\end{proof}

\begin{proposition}\label{thm:irregular}
Let a finite abelian group $G$ with $\tau_3(G)=2$ act faithfully on a projective terminal rationally connected threefold $Y$. Assume that
\begin{enumerate}[label=(\arabic*)]
\item $-K_Y$ is nef,
\item $(-K_Y)^3=0$,
\item $2K_Y$ is Cartier,
\item $\plur_{-1}(Y)=0$,
\item $\plur_{-2}(Y)=2$, and
\item $\B(Y)=mQ_2$ for $m\in\{12,16\}$.
\end{enumerate}
Assume that a general $D\in|-2K_Y|$ is integral and normal with only Du Val singularities and avoids $\Sing_{\nG}Y$. Then $G$ is of product type.
\end{proposition}
\begin{proof}
Suppose that $G$ is not of product type.
If $m=12$, Lemma~\ref{lem:surfacecoh} and Corollary~\ref{cor:regularpencil} imply that $G$ is of product type. Hence $m=16$. 
Put $\Lambda_Y=|-2K_Y|$ and identify $\PP^1$ with $\Lambda_Y$. Let $H$ be the kernel of the action of $G$ on $\Lambda_Y$.

\smallskip\noindent\textbf{Step 1. Lower bound for $|H|$.}
The representation of $G$ on $H^0(Y,-2K_Y)$ is diagonalizable, so the image $G/H$ of $G$ on $\Lambda_Y$ lies in a one-dimensional diagonal torus and is cyclic. Since $\tau_3(G)=2$, we have
\begin{equation}\label{eq:irregularG}
 G=\ZZ/2\ZZ\times\ZZ/d_2\ZZ\times\ZZ/d_3\ZZ\times\ZZ/d_4\ZZ,\qquad
 2\mid d_2\mid d_3\mid d_4.
\end{equation}
The list \eqref{eq:productlist} gives $d_2>2$ and $(d_2,d_3)\ne(4,4)$. Since $|G/H|\mid d_4$, we obtain
\begin{equation}\label{eq:Alower}
 |H|\ge2d_2d_3\ge64.
\end{equation}

Choose a sufficiently general $D\in|-2K_Y|$. Lemma~\ref{lem:surfacecoh} and Proposition~\ref{prop:anticanonical-Albanese} show that $D$ is smooth and minimal and that its Albanese morphism $a_D\colon D\to C_D=\Alb(D)$ has a general fiber of genus one. Moreover, we have
\begin{equation}\label{eq:irregularsurface}
 p_g(D)=0,\qquad q(D)=1,\qquad K_D\text{ nef},\qquad
 K_D^2=0,\qquad\plur_2(D)=1.
\end{equation}
The surface $D$ is therefore bielliptic or properly elliptic.

The group $H$ acts faithfully on a sufficiently general $D$, since each $h\in H\setminus\{1\}$ fixes pointwise only finitely many divisors in $\Lambda_Y$. Put $N=\ker(H\to\Aut(C_D))$ and $Q=H/N$. Lemma~\ref{lem:Albanese-fiber-group} gives
\begin{equation}\label{eq:Nbound}
 |N|\le9.
\end{equation}

\smallskip\noindent\textbf{Step 2. The unique multiple fiber.}
If $\kappa(D)=1$, the Albanese morphism $a_D$ is the canonical genus-one fibration of $D$. Write the multiple fibers of $a_D$ as $m_iF_i'$, $i=1,\ldots,s$, where $F_i'$ is reduced, $m_i\ge2$, and $p_i=a_D(F_i')\in C_D$. By the canonical bundle formula with $\chi(\OO_D)=0$ of \cite[\S E.7.3]{Reid}, there is a line bundle $M$ on $C_D$ such that
\[
 K_D\sim a_D^*M+\sum_i(m_i-1)F_i',\qquad\deg M=0.
\]
If $s=0$, then $h^0(D,mK_D)=h^0(C_D,M^{\otimes m})\le1$ for every $m\ge1$, since $a_D$ has connected fibers and $\deg M=0$. This contradicts $\kappa(D)=1$. Thus $s>0$, and pushforward of $\OO_D(2K_D)$ gives
\[
 (a_D)_*\OO_D(2K_D)
 \simeq M^{\otimes2}\otimes\OO_{C_D}(p_1+\cdots+p_s).
\]
Since $\deg\bigl(M^{\otimes2}\otimes\OO_{C_D}(p_1+\cdots+p_s)\bigr)=s>0$, we have
\[
 \plur_2(D)=h^0\bigl(C_D,M^{\otimes2}\otimes\OO_{C_D}(p_1+\cdots+p_s)\bigr)=s.
\]
Hence $s=\plur_2(D)=1$.

The group $Q$ fixes $p_1\in C_D$, the image of the unique multiple fiber of $a_D$. Taking $p_1$ as the origin of $C_D$, we have $Q\subseteq\Aut(C_D,p_1)$, which is cyclic of order two, four, or six \cite[Theorem~III.10.1]{Silverman}. Thus $Q$ is cyclic and $|Q|\le6$. By \eqref{eq:Nbound}, we have $|H|=|N||Q|\le54$, contradicting \eqref{eq:Alower}.

\smallskip\noindent\textbf{Step 3. Torsion sections of the relative Jacobian.}
Thus a sufficiently general $D$ is bielliptic. Since $K_D$ is torsion, \eqref{eq:irregularsurface} gives $2K_D\sim0$ and $K_D\not\sim0$. By adjunction, $\OO_D(-2K_Y)\simeq\OO_D(2K_D)\simeq\OO_D$. Let $s_D$ define $D$. The kernel of restriction from $H^0(Y,-2K_Y)$ to $H^0(D,\OO_D(-2K_Y))$ is $\CC s_D$. Thus the section defining any $D'\in\Lambda_Y\setminus\{D\}$ restricts to a nonzero constant under the isomorphism $\OO_D(-2K_Y)\simeq\OO_D$. Hence $D'\cap D=\varnothing$, so $\Lambda_Y$ is base point free and defines a $G$-equivariant morphism $f\colon Y\to\PP^1$.

Lemma~\ref{lem:relative-Albanese-surface} gives a $G$-equivariant genus-one fibration $b\colon B\to\PP^1$ with Jacobian fibration $j\colon J\to\PP^1$, where $J$ is a rational elliptic surface. For a sufficiently general $D=f^{-1}(t)$, we identify $B_t=b^{-1}(t)$ with $C_D$. The kernel of $G\to\Aut(B)$ is contained in $H$ because $b$ is surjective and $G$-equivariant. Each element of $H$ acting nontrivially on $B$ fixes pointwise only finitely many fibers of $b$. Thus for general $t$, the kernel $N$ of the action of $H$ on $C_D\simeq B_t$ equals $\ker(G\to\Aut(B))$. Put $\overline G=G/N\subseteq\Aut(B)$.

If some element of $Q$ is not a translation of $C_D$, applying Lemma~\ref{lem:bielliptic32} to $H\subseteq\Aut(D)$ gives $|H|\le32$, contrary to \eqref{eq:Alower}. Thus $Q$ acts by translations on the general fibers of $b$.

By \cite[Lemma~3.1]{CLS}, we have $Q\hookrightarrow\MW(J)_{\tors}$. In particular,
$
 |Q|\le9$ and $|H|\le81.
$
On the other hand, equations~\eqref{eq:Alower} and \eqref{eq:Nbound} give
$
 |Q|=\frac{|H|}{|N|}\ge\frac{64}{9}>4.
$
Apply Proposition~\ref{prop:translationbase} to $\overline G\subseteq\Aut(B)$ and its translation subgroup $Q$. The action of $G$ on $\PP^1$ is trivial, so $H=G$. But \eqref{eq:irregularG} and $d_2>2$ give
\[
 |G|=2d_2d_3d_4\ge2\cdot4^3=128>81\ge|H|,
\]
a contradiction. This completes the proof.
\end{proof}

\section{Excluding the baskets \texorpdfstring{$15Q_2$ and $14Q_2$}{15Q2 and 14Q2}}\label{sec:emptycase}
In this section, we exclude the case $\tau_3(G)=3$, including the basket $15Q_2$, and the basket $14Q_2$. We then prove Theorem~\ref{thm:nonemptymodel}.

\subsection{The case \texorpdfstring{$\tau_3(G)=3$}{tau3(G)=3}}\label{sec:three}

\begin{proposition}\label{prop:foursections}
Let a finite abelian group $G$ containing $(\ZZ/3\ZZ)^4$ act faithfully on a projective terminal rationally connected threefold $Y$. Let $\Lambda$ be a $G$-invariant pencil on $Y$ with integral general element. Assume that
\begin{enumerate}[label=(\arabic*)]
\item $-K_Y$ is nef,
\item $\plur_{-1}(Y)=0$,
\item $\plur_{-2}(Y)\le4$,
\item $\Lambda\subseteq|-2K_Y|$, and
\item $(Y,\tfrac12\Lambda)$ is terminal.
\end{enumerate}
Then $G$ is of product type.
\end{proposition}
\begin{proof}
Two distinct general elements of $\Lambda$ are integral and share no component, so $\Lambda$ is mobile. Lemma~\ref{lem:terminalpencil} gives that $2K_Y$ is Cartier and that a general element $D\in\Lambda$ is normal, Gorenstein, and disjoint from $\Sing_{\nG}Y$. Write $\B(Y)=mQ_2$ since $2K_Y$ is Cartier. Choose a subgroup $H\subseteq G$ isomorphic to $(\ZZ/3\ZZ)^4$. Since $\tau_3(H)=3$, Lemma~\ref{lem:basketorbits} gives $3\mid m$.

If $(-K_Y)^3>0$, then $-K_Y$ is nef and big. Equation~\eqref{eq:index-two-antipluri} and $3\mid m$ give $m\ge15$ and $\plur_{-2}(Y)=m-10\ge5$, contrary to the hypothesis. Hence $(-K_Y)^3=0$. Suppose that $G$ is not of product type. Lemma~\ref{lem:surfacecoh} and $3\mid m$ give $m=12$. The minimal resolution $T\to D$ then satisfies $p_g(T)=q(T)=0,$ $K_T$ is nef, $K_T^2=0$ and
$
 \plur_2(T)=\plur_{-2}(Y)-1\le3.
$
Corollary~\ref{cor:regularpencil} implies that $G$ is of product type, a contradiction.
\end{proof}

\begin{proposition}\label{prop:three}
Suppose that Assumption~\ref{ass:counterexample} holds. Then the equality $\tau_3(G)=3$ is impossible.
\end{proposition}
\begin{proof}
Suppose that $\tau_3(G)=3$. By Proposition~\ref{prop:baskets}, the group $G$ contains $(\ZZ/3\ZZ)^4$, and either $\plur_{-2}(X)=2$ or $\B(X)=15Q_2$.

\smallskip\noindent\textbf{Case 1. $\plur_{-2}(X)=2$.}
Apply Proposition~\ref{prop:nefmodel} to the pencil $\Lambda=|-2K_X|$. The resulting threefold $Y$ and pencil $\Lambda_Y$ satisfy the hypotheses of Proposition~\ref{prop:foursections}, so $G$ is of product type, a contradiction.

\smallskip\noindent\textbf{Case 2. $\B(X)=15Q_2$.}
Table~\ref{tab:baskets} gives $(-K_X)^3=3/2$ and $\plur_{-2}(X)=5$. Note that $2K_X$ is Cartier. By Lemma~\ref{lem:bigmember}, a general $D\in|-2K_X|$ is integral and normal with only Du Val singularities and avoids $\Sing_{\nG}X$. Its minimal resolution $S\to D$ is a minimal surface of general type with $p_g(S)=q(S)=0$ and $K_S^2=3$. Choose a section $s_D\in H^0(X,-2K_X)$ defining $D$.
Since $K_S^2=3$, Xiao's bicanonical-map theorem \cite[Theorem~1]{Xiao} implies that the bicanonical map of $S$ is generically finite.

Let $\phi_{-2K_X}\colon X\dashrightarrow Z\subset\PP^4$ be the rational map defined by $|-2K_X|$ to its image closure $Z$. By \eqref{eq:bigrestriction}, the image of $D$ under $\phi_{-2K_X}$ has dimension two, since the bicanonical map of $S$ is generically finite. This image lies in the hyperplane section of $Z$ defined by $s_D$. Since $Z$ is nondegenerate, we have $\dim Z=3$, so $\phi_{-2K_X}$ is generically finite. Write $e=\deg\phi_{-2K_X}$ and $d=\deg Z$. Choose a resolution $f\colon\widetilde X\to X$ of the base ideal of $|-2K_X|$, and write $f^*|-2K_X|=F_{\mathrm{fix}}+\MM$, where $F_{\mathrm{fix}}$ is the fixed divisor and $\MM$ is base point free. Write $f^*D=F_{\mathrm{fix}}+M$ with $M\in\MM$. The divisors $M$ and $f^*D$ are nef, and $f^*D-M=F_{\mathrm{fix}}$ is effective. Hence $ed=M^3\le(-2K_X)^3=12$. Since $Z$ is nondegenerate and has dimension three in $\PP^4$, we have $d\ge2$ and $e\le6$.

Choose a subgroup $H\subseteq G$ isomorphic to $(\ZZ/3\ZZ)^4$, and put $W=H^0(X,-2K_X)$ and $H_0=\ker(H\to\operatorname{PGL}(W))$. Since $H_0$ acts trivially on $Z$, its faithful action on $X$ induces an embedding $H_0\hookrightarrow\Aut_{\CC(Z)}(\CC(X))$. Thus $|H_0|=[\CC(X):\CC(X)^{H_0}]$ divides $[\CC(X):\CC(Z)]=e$. Therefore
\begin{equation}\label{eq:projectivekernelbound}
 |H_0|\le3.
\end{equation}
There is a point $P\in\Sing_{\nG}X$ with $|H\cdot P|=3$. Indeed, we have $\Fix_X(H)\cap\Sing_{\nG}X=\varnothing$, while $9\mid|H\cdot P'|$ for every $P'\in\Sing_{\nG}X$ would force the multiplicity $15$ of $Q_2$ in $\B(X)$ to satisfy $9\mid15$. The stabilizer $H_P=\{h\in H\mid h(P)=P\}$ has index three in $H$, so $H_P\simeq(\ZZ/3\ZZ)^3$.

The evaluation map $\operatorname{ev}_P\colon W\to\OO_X(-2K_X)|_P$ is $H_P$-equivariant and is nonzero since a general element of $|-2K_X|$ avoids $P$. Let $\chi_P$ be the character of $H_P$ on $\OO_X(-2K_X)|_P$. Choose a $G$-eigenbasis $s_1,\ldots,s_5$ of $W$ and let $\chi_i$ be the character of $H_P$ on $\CC s_i$. If $s_i(P)\ne0$, then $\operatorname{ev}_P$ identifies $\CC s_i$ with $\OO_X(-2K_X)|_P$, so $\chi_i=\chi_P$. The characters $\chi_i\chi_P^{-1}$ therefore have common kernel $H_P\cap H_0$, and only those with $s_i(P)=0$ can be nontrivial. Consequently \eqref{eq:projectivekernelbound} gives
\[
 \#\{i\mid s_i(P)=0\}\ge\rk\bigl(H_P/(H_P\cap H_0)\bigr)\ge2.
\]

After reordering the basis, let $\Lambda$ be the $G$-invariant pencil generated by $s_1,s_2$, both vanishing at $P$. By Corollary~\ref{cor:pencils}, the pencil $\Lambda$ has no fixed components and its general element is integral. By Corollary~\ref{cor:ct}, the pair $(X,\frac12\Lambda)$ is canonical. Choose an exceptional prime divisor $F$ over $X$ with center $P$ and discrepancy $a(F,X)=1/2$ \cite[p.~193]{KawamataDiscrepancy}. Since $-2K_X$ is Cartier and $P$ is a base point, we have $\ord_F\Lambda\in\ZZ_{\ge1}$. Canonicity gives $\ord_F\Lambda=1$, so $a(F,X,\frac12\Lambda)=\frac12-\frac12\ord_F\Lambda=0$. Thus every terminalization of $(X,\frac12\Lambda)$ extracts $F$.

Let $f_0\colon Y_0\to X$ be a $G$-equivariant projective crepant terminalization of $(X,\frac12\Lambda)$, and let $\Lambda_0$ be the strict transform of $\Lambda$. Since $f_0$ is crepant for $(X,\frac12\Lambda)$, we have $\Lambda_0\subseteq|-2K_{Y_0}|$, so $\plur_{-2}(Y_0)\ge2$. Let $\jmath\colon H^0(Y_0,-2K_{Y_0})\hookrightarrow W$ be the inclusion induced by $f_0$. For $s=\jmath(s_0)\ne0$, let $D_s$ and $D_{s_0}$ be the zero divisors of $s$ and $s_0$ on $X$ and $Y_0$, respectively. Then $\operatorname{coeff}_F(f_0^*D_s)=\operatorname{coeff}_F(D_{s_0})+2a(F,X)\ge1$. Since $F$ is centered at $P$, we have $s(P)=0$. Thus $\operatorname{im}\jmath\subseteq\ker\operatorname{ev}_P$, which has dimension $4$, and $\plur_{-2}(Y_0)\le4$.
Apply the minimal model program in the proof of Proposition~\ref{prop:nefmodel} to $(Y_0,\frac12\Lambda_0)$. The resulting map $Y_0\dashrightarrow Y$ is an isomorphism in codimension one, so it identifies $H^0(Y_0,-2K_{Y_0})$ with $H^0(Y,-2K_Y)$. Thus $2\le \plur_{-2}(Y)\le4$. Proposition~\ref{prop:foursections} now gives a contradiction.
\end{proof}

\subsection{The basket \texorpdfstring{$14Q_2$}{14Q2}}\label{sec:14}

\begin{proposition}\label{prop:14}
Let $X$ be a $G\QQ$-Fano threefold under Assumption~\ref{ass:counterexample}. Then $\B(X)\ne14Q_2$.
\end{proposition}
\begin{proof}
Assume that $\B(X)=14Q_2$. Then by Table~\ref{tab:baskets}, we have $(-K_X)^3=1$ and $\plur_{-2}(X)=4$. The Gorenstein index of every point of $X$ divides two, so $2K_X$ is Cartier. Lemma~\ref{lem:bigmember} gives a general element $D\in|-2K_X|$ that is integral and normal, has only Du Val singularities, and is disjoint from $\Sing_{\nG}X$. Let $\nu\colon S\to D$ be the minimal resolution of $D$. Then
$p_g(S)=q(S)=0$, $K_S^2=2$, and $K_S$ is nef and big. 
Thus $S$ is a numerical Campedelli surface \cite[\S1]{CMLP}. Choose a section $s_D\in H^0(X,-2K_X)$ defining $D$.

Put $W=H^0(X,-2K_X)$ and $H=\ker(G_2\to\operatorname{PGL}(W))$. The group $G_2$ has rank four and its image lies in a diagonal torus of dimension $\dim W-1=3$. Thus $H$ contains an involution $\sigma$. It preserves every element of $|-2K_X|$ and acts nontrivially on a sufficiently general $D$, since otherwise $\Fix_X(\sigma)$ would contain infinitely many distinct divisors.

Since $\sigma\in H$ and $\sigma^2=\id$, we have $\sigma^*|_W=\lambda\,\id_W$ for some $\lambda\in\{1,-1\}$. Put $L=\CC s_D$. Multiplication by sections in $L$ identifies $\OO_X(2K_X)\otimes L$ with the ideal sheaf $\mathscr I_D$ as $\langle\sigma\rangle$-linearized sheaves. Since $\nu$ is $\langle\sigma\rangle$-equivariant and crepant, adjunction and restriction in \eqref{eq:bigrestriction} give an isomorphism of $\langle\sigma\rangle$-representations
\[
 H^0(S,2K_S)\simeq(W/L)\otimes(L^\vee)^{\otimes2}.
\]
The involution $\sigma$ acts on $(L^\vee)^{\otimes2}$ by $\lambda^{-2}=1$, so $\sigma^*$ acts on $H^0(S,2K_S)$ by multiplication by $\lambda$. Thus the bicanonical map of $S$ factors through $S/\langle\sigma\rangle$. By \cite[(5) and Corollary~3.6]{CCM}, the $(-1)$-eigenspace of $\sigma^*$ on $H^0(S,2K_S)$ vanishes, so $\lambda=1$. Let $R_S$ be the reduced divisorial part of the fixed locus of $\sigma$ on $S$. By \cite[Proposition~3.1(i)]{CMLP}, we have $K_SR_S=2$.

At any fixed point $Q\in X_{\reg}$ of $\sigma$, choose analytic coordinates $x_1,x_2,x_3$ centered at $Q$ in which $\sigma$ acts by sign changes. The local generator $(dx_1\wedge dx_2\wedge dx_3)^{-2}$ of $\OO_X(-2K_X)$ is fixed by $\sigma^*$. Write
$
 s_D=f_D(dx_1\wedge dx_2\wedge dx_3)^{-2}.
$
Since $\lambda=1$, we have $f_D\circ\sigma=f_D$, so $d(f_D)_Q\circ d\sigma_Q=d(f_D)_Q$. Hence $d(f_D)_Q$ vanishes on the $-1$-eigenspace of $d\sigma_Q\colon T_QX\to T_QX$.

Suppose that $D$ contains a one-dimensional component $C$ of $\Fix_X(\sigma)$. At a general point $Q\in C$, the $+1$-eigenspace of $d\sigma_Q$ is $T_QC$. Since $f_D$ vanishes along $C$, the differential $d(f_D)_Q$ also vanishes on $T_QC$. Hence $D$ is singular along $C$, contrary to normality.

Let $R$ be the reduced divisorial fixed locus of $\sigma$ on $X$. A general $D$ contains no component of $R$. At a general point $Q$ of an irreducible component of $D\cap R$, the $+1$-eigenspace of $d\sigma_Q$ is $T_QR$. Since $D$ is smooth at $Q$ and $d(f_D)_Q$ vanishes on the $-1$-eigenspace of $d\sigma_Q$, the restriction of $d(f_D)_Q$ to $T_QR$ is nonzero. Thus $D$ and $R$ meet transversely at $Q$, and $D\cap R$ is generically reduced. The components of $R_S$ that map onto curves under $\nu$ are exactly the strict transforms on $S$ of the irreducible components of $D\cap R$. Every curve contracted by $\nu$ has intersection zero with $K_S$. Hence
\[
 2=K_SR_S=(-K_X)\cdot R\cdot(-2K_X)=2K_X^2R.
\]
Thus $K_X^2R=1=(-K_X)^3$. Since $R$ is $G$-invariant, Proposition~\ref{prop:classes} gives $R\sim-aK_X$ for an integer $a>0$. The equality $K_X^2R=(-K_X)^3$ implies $a=1$, contrary to $\plur_{-1}(X)=0$.
\end{proof}

\subsection{A nonempty anticanonical linear system in the given birational class}\label{sec:nonemptymodel}\label{sec:assembly}

\begin{proof}[Proof of Theorem~\ref{thm:nonemptymodel}]
By Corollary~\ref{cor:known-reduction}, the threefold $V$ is $G$-birational to a $G\QQ$-Fano threefold $X$. Suppose that \eqref{eq:geometriccontradiction} holds for every such threefold. Assumption~\ref{ass:counterexample} applies. Proposition~\ref{prop:baskets} gives $\tau_3(G)\in\{2,3\}$ and either $\plur_{-2}(X)=2$ or one of the four baskets in Table~\ref{tab:baskets}. 
Propositions~\ref{prop:B0} and~\ref{prop:B4} exclude $2Q_2+4Q_3+2Q_5$ and $12Q_2+2Q_3$. Proposition~\ref{prop:three} excludes $\tau_3(G)=3$, including the basket $15Q_2$. Proposition~\ref{prop:14} excludes $14Q_2$. Thus $\plur_{-2}(X)=2$ and $\tau_3(G)=2$. 
Apply Proposition~\ref{prop:nefmodel} to the complete pencil $|-2K_X|$. Together with Lemma~\ref{lem:surfacecoh}, it gives a threefold $Y$ satisfying the hypotheses of Proposition~\ref{thm:irregular}. Hence $G$ is of product type, contrary to Assumption~\ref{ass:counterexample}.
\end{proof}

\section{Anticanonical K3 surfaces}\label{sec:K3geometry}
In this section, we determine the inertia groups of $G$-invariant anticanonical K3 surfaces for $G\simeq\Gamma_i$, $1\le i\le4$. When $G\simeq\Gamma_2$, we describe the polarized K3 surface $(S,h)$, where $S$ is a $G$-invariant anticanonical divisor and $h$ is the ample generator of $\Pic(S)^G$ with $h^2=2$. When $G$ is isomorphic to $\Gamma_3$ or $\Gamma_4$, we describe the $G$-equivariant contractions of the minimal resolutions of the $G$-invariant anticanonical surfaces.

\subsection{Inertia groups}\label{sec:K3inertia}

Let $X$ be a $G\QQ$-Fano threefold with $G\simeq\Gamma_i$ for some $i\in\{1,2,3,4\}$ and with nonempty anticanonical linear system.
For an integral ample $\QQ$-Cartier divisor $D$ on a normal projective variety $Z$, write
\[
 R(Z,D)=\bigoplus_{m\ge0}H^0(Z,\OO_Z(mD)),
\]
where $\OO_Z(mD)$ denotes the associated reflexive divisorial sheaf.

\begin{proposition}\label{prop:K3inputs}
\leavevmode
\begin{enumerate}[label=\textup{(\arabic*)}]
\item\label{item:K3-anticanonical} If a finite abelian group $G$ which is not of product type acts faithfully on a $G\QQ$-Fano threefold $X$ with $|{-}K_X|\ne\varnothing$, then every $G$-invariant element of $|{-}K_X|$ is an integral normal K3 surface with at worst Du Val singularities.
\item\label{item:K3-groups} A finite abelian group acting faithfully on a smooth K3 surface and admitting no embedding into $\Cr_2(\CC)$ is isomorphic to one of
\[
\begin{gathered}
 (\ZZ/4\ZZ)^3,\qquad (\ZZ/6\ZZ)^2\times\ZZ/2\ZZ,\qquad
 \ZZ/6\ZZ\times(\ZZ/3\ZZ)^2,\\
 \ZZ/8\ZZ\times\ZZ/4\ZZ\times\ZZ/2\ZZ,\qquad
 (\ZZ/2\ZZ)^5,\qquad \ZZ/4\ZZ\times(\ZZ/2\ZZ)^3.
\end{gathered}
\]
\item\label{item:K3-lattices} Write
\[
\begin{alignedat}{2}
 &\mathmakebox[\kthreeleftwidth][l]{H_1=(\ZZ/4\ZZ)^3,}
 &\qquad &\mathmakebox[\kthreerightwidth][l]{H_2=(\ZZ/6\ZZ)^2\times\ZZ/2\ZZ,}\\
 &H_3=\ZZ/6\ZZ\times(\ZZ/3\ZZ)^2,
 &&H_4=\ZZ/8\ZZ\times\ZZ/4\ZZ\times\ZZ/2\ZZ.
\end{alignedat}
\]
Let one of the groups $H_1,\ldots,H_4$ act faithfully on a smooth K3 surface $T$. Then the invariant integral cohomology lattices are
\begin{equation}\label{eq:K3lattices}
\begin{alignedat}{2}
 &\mathmakebox[\kthreeleftwidth][l]{H^2(T,\ZZ)^{H_1}\simeq\langle4\rangle,}
 &\qquad &\mathmakebox[\kthreerightwidth][l]{H^2(T,\ZZ)^{H_2}\simeq\langle2\rangle,}\\
 &H^2(T,\ZZ)^{H_3}\simeq U(3),
 &&H^2(T,\ZZ)^{H_4}\simeq U(2),
\end{alignedat}
\end{equation}
where $U(e)$ has Gram matrix $\left(\begin{smallmatrix}0&e\\e&0\end{smallmatrix}\right)$. In the last two cases the group contains a purely nonsymplectic element of order six and eight, respectively.
\end{enumerate}
\end{proposition}
\begin{proof}
Part~(1) follows from \cite[\S12, proof of Theorem~1.7, and Proposition~8.2]{Log}. Part~(2) is \cite[Proposition~1.7]{LPZ}. Every faithful action of a group $H$ on a Du Val K3 surface induces a faithful $H$-action on its minimal resolution. Part~(3) follows from the computations with the database of \cite{BrandhorstHofmann} in \cite[Corollary~3.9(1)--(4) and Proposition~4.4(1)--(4)]{LPZ}. Each group $H_i$ contains an element acting nonsymplectically on $T$, so the invariant lattices in \eqref{eq:K3lattices} are contained in $\NS(T)$.
\end{proof}

Let $S\in|{-}K_X|$ be $G$-invariant. By Proposition~\ref{prop:K3inputs}\ref{item:K3-anticanonical}, the surface $S$ is a normal K3 surface with at worst Du Val singularities. Write
\[
 I_S=\ker(G\to\Aut(S)),\qquad H_S=G/I_S
\]
for the inertia subgroup of $G$ and the quotient group which faithfully acts on $S$, respectively.

\begin{lemma}\label{lem:inertiatypes}
Let $X$ be a $G\QQ$-Fano threefold with $G\simeq\Gamma_i$, where $1\le i\le4$, and let $S\in|-K_X|$ be $G$-invariant. Then $H_S\simeq H_i$, and $I_S$ is cyclic of order $4$ for $i=1$, of order $6$ for $i=2,3$, and of order $8$ for $i=4$.
\end{lemma}
\begin{proof}
For $i=1$, the assertion is \cite[Proposition~2.11]{Lo26}. Suppose that $i\in\{2,3,4\}$.
At a general point $Q\in S\cap X_{\reg}$ where $S$ is smooth, Lemma~\ref{lem:tangent} gives a faithful representation of $I_S$ on $T_QX$. Since $I_S$ fixes $S$ pointwise, its representation on the normal line $T_QX/T_QS$ is faithful. Thus $I_S$ is cyclic. If $H_S$ admitted an embedding into $\Cr_2(\CC)$, Proposition~\ref{prop:extensions} would imply that $G$ is of product type. Hence Proposition~\ref{prop:K3inputs}\ref{item:K3-groups} applies.

The orders of $\Gamma_2,\Gamma_3,\Gamma_4$ are $432,324,512$, and the exponents of $\Gamma_2,\Gamma_3,\Gamma_4$ are $6,6,8$, respectively. The six groups in \cite[Proposition~1.7]{LPZ} have orders $64,72,54,64,32,32$. The order of $I_S$ divides the exponent of $G$. Comparing orders gives $|I_S|=6$ and $H_S\simeq H_i$ for $i=2,3$. For $i=4$, the kernel $I_S$ has order eight and $H_S$ is isomorphic to either $H_4$ or $(\ZZ/4\ZZ)^3$. If $H_S\simeq(\ZZ/4\ZZ)^3$, then
$
 4\Gamma_4\simeq(\ZZ/2\ZZ)^2\subset I_S,
$
which is impossible for a cyclic kernel.
\end{proof}

\begin{lemma}\label{lem:distinctinertia}
We have
\begin{equation}\label{eq:nonemptysectionbounds}
 \plur_{-1}(X)\le
 \begin{cases}
 5,&i=1,\\
 4,&i=2,\\
 3,&i=3,4.
 \end{cases}
\end{equation}
If $i=1$, then in any $G$-eigenbasis of $H^0(X,-K_X)$ the image of $G$ in $\operatorname{PGL}(H^0(X,-K_X))$ consists of all diagonal elements of order dividing $4$.
\end{lemma}
\begin{proof}
Put $W=H^0(X,-K_X)$ and $n=\dim W$. Both assertions are immediate for $n=1$, so assume $n\ge2$. Choose a $G$-eigenbasis $s_1,\ldots,s_n$ of $W$ for the action induced by the canonical $G$-linearization of $\OO_X(-K_X)$. By Proposition~\ref{prop:K3inputs}\ref{item:K3-anticanonical}, the divisors $S_j=(s_j=0)$ are distinct integral normal K3 surfaces. Put $I_j=\ker(G\to\Aut(S_j))$. By Lemma~\ref{lem:inertiatypes}, each $I_j$ is cyclic of order $d=|I_S|$.

For each $j$, choose a general point $Q\in S_j\cap X_{\reg}$ where $S_j$ is smooth, and define the faithful character $\nu_j\colon I_j\to\CC^\times$ by $g^*|_{N^*_{S_j/X,Q}}=\nu_j(g)\id$ for $g\in I_j$. Since $I_j$ fixes $S_j$ pointwise, it acts on $\OO_X(-K_X)|_Q$ by $\nu_j^{-1}$. The linear term of $s_j$ at $Q$ spans $N^*_{S_j/X,Q}\otimes\OO_X(-K_X)|_Q$, on which $I_j$ acts by $\nu_j\nu_j^{-1}=1$. For $\ell\ne j$, the value $s_\ell(Q)$ is nonzero in $\OO_X(-K_X)|_Q$. Since the $s_\ell$ are $G$-eigenvectors, for every $g\in I_j$ we have
\[
 g^*s_j=s_j,\qquad g^*s_\ell=\nu_j(g)^{-1}s_\ell,\quad \ell\ne j.
\]

Let $\overline G$ be the image of $G$ in $\operatorname{PGL}(W)$. In the basis $s_1,\ldots,s_n$, the image of $I_j$ consists of the classes of diagonal matrices that multiply $s_j$ by $\xi_d^a$ for $a\in\ZZ/d\ZZ$ and fix $s_\ell$ for $\ell\ne j$. The images of $I_1,\ldots,I_n$ in $\overline G$ generate a subgroup isomorphic to
\[
 (\ZZ/d\ZZ)^n/\langle(1,\ldots,1)\rangle\simeq(\ZZ/d\ZZ)^{n-1}.
\]
For $i=1$, this subgroup consists of all diagonal elements of order dividing $4$. Since $\overline G$ is diagonal and its exponent divides $4$, it equals this subgroup. Hence $4^{n-1}=|\overline G|\le|G|=4^4$, and $n\le5$. For $i=2$, the $3$-primary subgroup of $\overline G$ is a quotient of $(\ZZ/3\ZZ)^3$, so $n-1\le3$. For $i=3$, the $2$-primary subgroup of $\overline G$ is a quotient of $(\ZZ/2\ZZ)^2$, so $n-1\le2$. For $i=4$, the group $4\overline G$ is a quotient of $4G\simeq(\ZZ/2\ZZ)^2$ and contains $(\ZZ/2\ZZ)^{n-1}$. Hence $n-1\le2$ in this case as well.
\end{proof}

\subsection{Polarized K3 surfaces}\label{sec:markedK3}

\subsubsection{The degree-two polarization}
\begin{lemma}\label{lem:sexticK3}
Let $S$ be a Du Val K3 surface with a faithful $H_2$-action. Then $S$ is smooth, and $\Pic(S)^{H_2}$ is generated by an ample class $h$ with $h^2=2$. The section ring of $(S,h)$ is
\begin{equation}\label{eq:sexticK3ring}
 R(S,h)\simeq
 \CC[x_1,x_2,x_3,w]/(w^2+x_1^6+x_2^6+x_3^6),
 \qquad\deg x_i=1,\quad\deg w=3.
\end{equation}
In the coordinates of \eqref{eq:sexticK3ring}, let $D_2\subseteq\Aut(S)$ be the subgroup consisting of the maps
\[
 [x_1:x_2:x_3:w]\longmapsto[\xi_6^a x_1:\xi_6^b x_2:x_3:\xi_2^c w],
 \qquad a,b\in\ZZ/6\ZZ,\quad c\in\ZZ/2\ZZ.
\]
The image of $H_2$ in $\Aut(S)$ is $D_2$.
\end{lemma}
\begin{proof}
The surface $S$ is smooth by \cite[Corollary~4.6]{LPZ}. By \eqref{eq:K3lattices}, the group $\Pic(S)^{H_2}$ is generated by a class $h$ with $h^2=2$. Choose the sign of $h$ so that an $H_2$-invariant ample class is a positive multiple of $h$. Then $h$ is ample.

By \cite[Proposition~8.1 and \S5.1]{SaintDonat}, the complete linear system $|h|$ defines a double cover $\pi_h\colon S\to\PP^2$ branched along a sextic $B_6=(F_6=0)$. The curve $B_6$ is smooth because $S$ is smooth.

Since $h$ is $H_2$-invariant, the action of $H_2$ descends along $\pi_h$ and gives a homomorphism $\rho_h\colon H_2\to\operatorname{PGL}_3(\CC)$. Put $\overline H_2=\operatorname{im}\rho_h$. Since $c(H_2)=13/6>2$, Proposition~\ref{prop:complexityprojective} gives $\ker\rho_h\ne\{1\}$. The elements of $\ker\rho_h$ act over $\PP^2$, so $\ker\rho_h=\langle\iota_h\rangle$, where $\iota_h$ is the deck involution of $\pi_h$. Thus $\overline H_2$ has order $36$ and exponent dividing $6$, hence $\overline H_2\simeq(\ZZ/6\ZZ)^2$.
By \cite[Proposition~6.1.2]{BlancThesis}, the subgroup $\overline H_2\simeq(\ZZ/6\ZZ)^2$ is diagonalizable. After a change of coordinates, it is therefore the full $6$-torsion of the diagonal torus.

Since $\overline H_2$ preserves $B_6$, $F_6$ is semi-invariant under the diagonal action of $\overline H_2$ on $\PP^2$. The group $\overline H_2$ is generated by $[x_1:x_2:x_3]\mapsto[\xi_6x_1:x_2:x_3]$ and $[x_1:x_2:x_3]\mapsto[x_1:\xi_6x_2:x_3]$. Thus for $i=1,2$, the exponent of $x_i$ in each monomial of $F_6$ has the same residue $r_i\in\{0,\ldots,5\}$ modulo $6$. Since $x_1^{r_1}x_2^{r_2}$ divides $F_6$ and $B_6$ is irreducible, we have $r_1=r_2=0$. Hence $F_6=a_1x_1^6+a_2x_2^6+a_3x_3^6$, with $a_1a_2a_3\ne0$ because $a_i=0$ would make the $i$th coordinate point singular on $B_6$. Writing the double cover as $w^2+F_6=0$ and rescaling $x_1,x_2,x_3$ and $w$ gives \eqref{eq:sexticK3ring}.

For a quadratic extension $L/K$, each automorphism of $K$ has at most two lifts to $L$, since two such lifts differ by an element of $\operatorname{Gal}(L/K)$. For the extension $\CC(S)/\CC(\PP^2)$ defined by $\pi_h$, both $H_2$ and $D_2$ induce $\overline H_2$ and have order $2|\overline H_2|$. Hence the image of $H_2$ in $\Aut(S)$ is $D_2$.
\end{proof}

\subsubsection{The actions of $H_3$ and $H_4$}
Consider the weighted projective hypersurfaces in \cite[Example~3.11, Table~4, rows~(3) and~(2)]{LPZ}.
\begin{equation}\label{eq:ranktwoK3models}
\begin{aligned}
 S_3&=\{x^6+y^6+z^3+w^3=0\}\subset\PP(1,1,2,2),\\
 S_4&=\{x^8+y^8+z^4+w^2=0\}\subset\PP(1,1,2,4).
\end{aligned}
\end{equation}
The surfaces $S_3,S_4$ are K3 surfaces with Du Val singularities. For $i=3,4$, let $\mu_i\colon T_i\to S_i$ be the minimal resolution. The group $H_i$ acts on $S_i$ by fixing $x$ and multiplying each of $y,z,w$ by $\xi_m^a$, where $m$ is the exponent of that coordinate in the equation of $S_i$ and $a\in\ZZ/m\ZZ$ is chosen independently for each coordinate. The $H_i$-action on $S_i$ lifts to $T_i$.

\begin{proposition}\label{prop:markedK3}
If a smooth projective K3 surface $T$ admits a faithful action of $H_3$ or $H_4$, then the given action is isomorphic to the standard action of $H_3$ on $T_3$ or of $H_4$ on $T_4$, after precomposing by an automorphism of the acting group.
\end{proposition}
\begin{proof}
Put $H=H_i$ for $i\in\{3,4\}$ and let $d=6$ for $i=3$ and $d=8$ for $i=4$. Let $\chi\colon H\to\CC^\times$ be the character on $H^0(T,\Omega_T^2)$. By Proposition~\ref{prop:K3inputs}\ref{item:K3-lattices}, we can choose an element $g\in H$ of order $d$ with $\chi(g)=\xi_d$. Put $\iota=g^{d/2}$. By formula~\eqref{eq:holomorphicLefschetz}, the holomorphic Lefschetz number of $g$ on $\OO_T$ is $1+\xi_d^{-1}\ne0$. Hence $g$, and thus $\iota$, has a fixed point. Since $\chi(\iota)=-1$, the quotient $Q=T/\langle\iota\rangle$ is a smooth rational surface and the quotient morphism $\pi\colon T\to Q$ is a double cover with nonempty reduced branch divisor $B\sim-2K_Q$ \cite[Lemma~1.2(1), (2), and (4)]{ZhangK3}.

The quotient group $\overline H=H/\langle\iota\rangle$ is $(\ZZ/3\ZZ)^3$ for $H_3$ and $(\ZZ/4\ZZ)^2\times\ZZ/2\ZZ$ for $H_4$. Run a $\overline H$-equivariant minimal model program on $Q$. If the program ends with a conic bundle $f\colon Y'\to C$, then $C\simeq\PP^1$ because $Y'$ is rational. Write $\theta\colon\overline H\to\Aut(C)$ for the action on $C$, and put $N=\ker\theta$ and $J=\operatorname{im}\theta$. Since $N$ acts faithfully on the geometric generic fiber of $f$, which is a projective line, each of $N$ and $J$ is cyclic or isomorphic to $(\ZZ/2\ZZ)^2$. The exponents of $N$ and $J$ divide $3$ when $H=H_3$ and $4$ when $H=H_4$, so we obtain
$
 |\overline H|=|N|\,|J|\le9$ or $16,
$
contrary to $|\overline H|=27$ or $32$. The program therefore produces a smooth del Pezzo surface $Y$ with $\rk\Pic(Y)^{\overline H}=1$.

If $K_Y^2=1$, the unique base point $P$ of $|-K_Y|$ is fixed by $\overline H$ \cite[\S8.4]{BlancThesis}. Lemma~\ref{lem:tangent} gives a faithful representation of $\overline H$ on $T_PY$. Since $\overline H$ is abelian, this implies $\rk\overline H\le2$, whereas $\rk\overline H=3$ in both cases. We have $K_Y^2\le3$ by \cite[Propositions~6.2.4, 6.4.1, and~8.1.14]{BlancThesis}. The lists in \cite[Propositions~8.2.1 and~8.3.3]{BlancThesis} identify the remaining pairs, up to isomorphism, with
\[
\begin{aligned}
 Y_3&=\{x^3+y^3+z^3+w^3=0\}\subset\PP^3,\\
 Y_2&=\{w^2+x^4+y^4+z^4=0\}\subset\PP(1,1,1,2),
\end{aligned}
\]
with the standard diagonal actions on $Y_3$ and $Y_2$. The two actions are cases $3.333$ and $2.G44$ in \cite[Propositions~8.2.1 and~8.3.3]{BlancThesis}, respectively.

Put $Y=Y_3$ for $i=3$ and $Y=Y_2$ for $i=4$. Let $\beta\colon Q\to Y$ be the composition of the $\overline H$-equivariant contractions. Pushing forward $B\sim-2K_Q$ along $\beta$ gives $B_Y=\beta_*B\sim-2K_Y$. The divisor $B_Y$ is reduced and $\overline H$-invariant, and ampleness of $-K_Y$ gives $B_Y\ne0$. The morphism $\beta$ is an isomorphism over the generic point of every prime divisor on $Y$, so the normalization of $Y$ in $\CC(T)$ has branch divisor $B_Y$. Since $Y$ is a smooth rational surface, we have $\Pic(Y)[2]=0$. Thus $B_Y$ determines the normal double cover of $Y$ up to isomorphism over $Y$ \cite[\S2.1]{AbelianCovers}.

Let $P_Y$ be the ambient weighted projective space of $Y$, and let $F_Y$ be the defining polynomial of $Y$ of degree $d_Y\in\{3,4\}$. Since $Y$ is contained in the smooth locus of $P_Y$ and the weights of $P_Y$ sum to $d_Y+1$, adjunction gives $\OO_Y(K_Y)\simeq\OO_Y(-1)$. Multiplication by $F_Y$ gives the exact sequence
\[
 0\longrightarrow\OO_{P_Y}(2-d_Y)
 \xrightarrow{\,F_Y\,}\OO_{P_Y}(2)
 \longrightarrow\OO_Y(2)\longrightarrow0.
\]
Since $d_Y>2$ and $H^1(P_Y,\OO_{P_Y}(2-d_Y))=0$, restriction gives an isomorphism $H^0(P_Y,\OO_{P_Y}(2))\simeq H^0(Y,\OO_Y(2))$. Thus $B_Y\in|-2K_Y|$ is cut out on $Y$ by a polynomial $Q_Y$ of degree $2$, unique up to scalar.

For $Y=Y_k$ with $k\in\{2,3\}$, the image of $\overline H$ in $\Aut(Y_k)$ is generated by the automorphisms of $Y_k$ induced by $D_{k,y},D_{k,z},D_{k,w}$. In the coordinates $(x,y,z,w)$, these matrices are
\[
\begin{aligned}
 D_{3,y}&=\diag(1,\xi_3,1,1),&
 D_{3,z}&=\diag(1,1,\xi_3,1),&
 D_{3,w}&=\diag(1,1,1,\xi_3),\\
 D_{2,y}&=\diag(1,\xi_4,1,1),&
 D_{2,z}&=\diag(1,1,\xi_4,1),&
 D_{2,w}&=\diag(1,1,1,-1).
\end{aligned}
\]
Since $B_Y$ is $\overline H$-invariant, we have $D_{k,u}^*Q_{Y_k}\in\CC^\times Q_{Y_k}$ for $u\in\{y,z,w\}$.
Thus the exponents of $y,z,w$ in the monomials of $Q_Y$ with nonzero coefficients are constant modulo $3,3,3$ for $Y=Y_3$ and modulo $4,4,2$ for $Y=Y_2$. The exponents of $y,z,w$ are smaller than the corresponding moduli because $Q_Y$ has degree $2$ and the coordinate $w$ has degree $2$ for $Y=Y_2$. Thus the exponents of $y,z,w$ are the same in every monomial of $Q_Y$ with nonzero coefficient. The equality $\deg Q_Y=2$ then determines the exponent of $x$. Hence $Q_Y$ is proportional to a single monomial.
Permutations of the coordinates of weight $1$ normalize $\overline H$. Since $B_Y$ is reduced, we may therefore take $B_Y=Y\cap\{xy=0\}$ or, for $Y=Y_2$, $B_Y=Y_2\cap\{w=0\}$.

If $B_Y=Y_2\cap\{w=0\}$, put $u=y/x$, $v=z/x$, and $r=w/x^2$. We have $\operatorname{div}_{Y_2}(r)=B_Y-2(Y_2\cap\{x=0\})$, so $\CC(Y_2)(\sqrt r)/\CC(Y_2)$ ramifies exactly along the components of $B_Y$. By uniqueness of the double cover with branch divisor $B_Y$, we may identify $\CC(T)$ with $\CC(Y_2)(\sqrt r)=\CC(u,v,\sqrt r)$. The matrices $D_{2,y},D_{2,z},D_{2,w}$ multiply $u,v,r$ by $\xi_4,\xi_4,-1$, respectively, with each matrix fixing the other two functions. Hence every $g\in H$ multiplies $u,v,\sqrt r$ by powers of $\xi_4$. Since $H$ acts faithfully on $\CC(T)=\CC(u,v,\sqrt r)$, we have $\exp H\mid4$, contrary to $\exp H_4=8$.

Assume that $B_Y=Y\cap\{xy=0\}$. Write $a,b$ for the first two coordinates of $S_i$ in \eqref{eq:ranktwoK3models}, and let $m$ be the weight of $w$ on $Y$. The morphism
\[
 \overline\pi\colon S_i\longrightarrow Y,\qquad
 [a\colon b\colon z\colon w]\longmapsto[a^2\colon b^2\colon z\colon w]
\]
is well defined because $a^2,b^2,z,w$ have weights $2,2,2,2m$ on $S_i$ and no common zero. The morphism $\overline\pi$ is finite because $a,b$ are integral over the image of the homogeneous coordinate ring of $Y$.
The functions $z/a^2$ and $w/a^{2m}$ are the pullbacks of $z/x$ and $w/x^m$, so $\CC(S_i)=\CC(Y)(b/a)$ with $(b/a)^2=y/x$. Since $B_Y$ is reduced, the divisor of $y/x$ has coefficients $\pm1$ along the components of $B_Y$ and zero elsewhere. Thus the branch divisor of $\overline\pi$ is $B_Y$. By uniqueness of the normal double cover of $Y$ with branch divisor $B_Y$, we identify $S_i$ with the normalization of $Y$ in $\CC(T)$. Thus there is a birational map $\psi\colon T\dashrightarrow S_i$ fitting into the commutative diagram
\[
\begin{tikzcd}[column sep=large,row sep=large]
 T \arrow[r,dashed,"\psi"] \arrow[d,"\pi"']
 & S_i \arrow[d,"\overline\pi"] \\
 Q \arrow[r,"\beta"']
 & Y
\end{tikzcd}.
\]

The generators of $H_i$ multiplying $b,z,w$ separately in \eqref{eq:ranktwoK3models} induce the automorphisms defined by $D_{k,y},D_{k,z},D_{k,w}$ on $Y=Y_k$. Thus the actions of $H$ on $T$ and $H_i$ on $S_i$ induce $\overline H$ on $Y$ along $\beta\circ\pi$ and $\overline\pi$, respectively. Via $\psi^*\colon\CC(S_i)\to\CC(T)$, the groups $H$ and $H_i$ both act on $\CC(T)$ by lifts of $\overline H$ and have order $2|\overline H|$. By the bound in the proof of Lemma~\ref{lem:sexticK3}, both groups give all lifts of $\overline H$ to $\CC(T)$. The induced identification $H\simeq H_i$ makes the birational map $T\dashrightarrow T_i$ induced by $\psi$ $H_i$-equivariant. Since $T$ and $T_i$ are smooth projective K3 surfaces, the birational map $T\dashrightarrow T_i$ is an isomorphism.
\end{proof}

\begin{lemma}\label{lem:K3contractions}
Let $i\in\{3,4\}$. The target of every $H_i$-equivariant birational morphism from $T_i$ to a Du Val K3 surface is $H_i$-equivariantly isomorphic to $T_i$ or $S_i$. Moreover,
\begin{equation}\label{eq:invariantWeilclasses}
 \Cl(S_i)^{H_i}=\ZZ A_0,\qquad A_0=[\OO_{S_i}(1)].
\end{equation}
\end{lemma}
\begin{proof}
Put $(T,H)=(T_i,H_i)$ and let $e=3$ for $i=3$ and $e=2$ for $i=4$. The singular locus of $S_i$ consists of $e$ points of type $A_1$. They are cut out by $z^3+w^3=0$ on $\PP(2,2)$ for $i=3$ and by $z^4+w^2=0$ on $\PP(2,4)$ for $i=4$. Let $E_1,\ldots,E_e$ be the exceptional curves of $\mu_i$ and put $R=\sum_{j=1}^e E_j$. Since $H$ acts transitively on $\Sing S_i$, it also permutes $E_1,\ldots,E_e$ transitively.

The rational map $S_i\dashrightarrow\PP^1$ given by $[x\colon y]$ is defined on $S_i\setminus\Sing S_i$. All points of $\Sing S_i$ lie in $z\ne0$. Let $\widetilde U_i\subset\AA^3$ be the surface obtained by setting $z=1$ in the equation of $S_i$ in \eqref{eq:ranktwoK3models}. The morphism $\pi_i\colon\widetilde U_i\to S_i\cap\{z\ne0\}$ given by $(x,y,w)\mapsto[x\colon y\colon 1\colon w]$ is the quotient by the involution $\delta_i(x,y,w)=(-x,-y,w)$. For $P\in\Sing S_i$, the unique point $\widetilde P$ with $\pi_i(\widetilde P)=P$ has the form $(0,0,w_P)$ with $w_P\ne0$. The partial derivative with respect to $w$ of the defining polynomial of $\widetilde U_i$ at $\widetilde P$ is $3w_P^2$ for $i=3$ and $2w_P$ for $i=4$. Hence $x,y$ are analytic coordinates at $\widetilde P$.

Let $W_P$ be the blowup of $(\widetilde U_i,\widetilde P)$ at $\widetilde P$, with exceptional curve $\widetilde E_P=\PP(T_{\widetilde P}\widetilde U_i)$, and let $\delta_P$ be the lift of $\delta_i$ to $W_P$. In the coordinates $x,y$, we have $d_{\widetilde P}\delta_i=-\id$, so $\delta_P|_{\widetilde E_P}=\PP(d_{\widetilde P}\delta_i)=\id_{\widetilde E_P}$. The map $[x\colon y]$ extends to $W_P$ and restricts to the isomorphism $[dx_{\widetilde P}\colon dy_{\widetilde P}]$ from $\widetilde E_P$ to $\PP^1$. Since $\delta_i$ preserves $[x\colon y]$, the morphism $W_P\to\PP^1$ descends to $f_P\colon W_P/\langle\delta_P\rangle\to\PP^1$. The quotient $W_P/\langle\delta_P\rangle$ is the minimal resolution of $(S_i,P)$, and $f_P$ maps the image of $\widetilde E_P$ isomorphically onto $\PP^1$. Thus $[x\colon y]$ defines a morphism $f\colon T\to\PP^1$. Each $E_j$ is a section of $f$. Let $F$ be a general fiber of $f$. We have
\begin{equation}\label{eq:FRintersections}
 F^2=0,\qquad FR=e,\qquad R^2=-2e.
\end{equation}
The classes of $F$ and $R$ are $H$-invariant. The classes $F,F+R$ generate $U(e)$. By \eqref{eq:K3lattices}, the sublattice $\ZZ F\oplus\ZZ(F+R)$ has the same rank and determinant as $\NS(T)^H$. Consequently
\begin{equation}\label{eq:FRbasis}
 \NS(T)^H=\ZZ F\oplus\ZZ(F+R).
\end{equation}

Let $\mu\colon T\to S$ be a nontrivial $H$-equivariant contraction to a Du Val K3 surface. Let $D$ be the reduced divisor on $T$ whose components are the curves contracted by $\mu$ other than $E_1,\ldots,E_e$. The divisor $D$ is $H$-invariant because $\mu$ is $H$-equivariant and $H$ permutes $E_1,\ldots,E_e$. If $D\ne0$, then $FD\ge0$ by nefness of $F$ and $RD\ge0$ because $R$ and $D$ have no common component. Write $D\sim aF+b(F+R)$ using \eqref{eq:FRbasis} and $\Pic(T)=\NS(T)$. By \eqref{eq:FRintersections}, we have $b\ge0$ and $a-b\ge0$, so $D^2=2eab\ge0$, contrary to negative definiteness of the intersection matrix of the curves contracted by $\mu$. Hence $D=0$. Since $\mu$ contracts at least one curve and $H$ permutes $E_1,\ldots,E_e$ transitively, the morphisms $\mu$ and $\mu_i$ contract exactly the curves $E_1,\ldots,E_e$. Therefore $S$ and $S_i$ are $H$-equivariantly isomorphic.

The homomorphism $(\mu_i)_*\colon\Pic(T)\to\Cl(S_i)$ is surjective and its kernel is generated by the classes of $E_1,\ldots,E_e$. Since the $E_j$ are disjoint and $E_j^2=-2$, every class $A\in\Cl(S_i)$ has a unique lift $M\in\Pic(T)$ with $(\mu_i)_*M=A$ and $ME_j\in\{0,1\}$ for all $j$. If $A$ is $H$-invariant, then $h^*M$ satisfies the same conditions for every $h\in H$ because $H$ permutes the $E_j$. Hence $M$ is $H$-invariant. Since $\Pic(T)=\NS(T)$ and $(\mu_i)_*R=0$, equation~\eqref{eq:FRbasis} gives $\Cl(S_i)^H=\ZZ[(\mu_i)_*F]$. For $F=f^{-1}([a\colon b])$, the divisor $(\mu_i)_*F$ is cut out by $bx-ay$ on $S_i$. Thus $[(\mu_i)_*F]=[\OO_{S_i}(1)]=A_0$.
\end{proof}

\section{Reconstruction}\label{sec:reconstruction}
In this section, we reconstruct the four threefolds $X_i$ from their anticanonical K3 surfaces and complete the proofs of Theorems~\ref{thm:complexity}, \ref{thm:cremona}, and~\ref{thm:rc-models}.

\subsection{Anticanonical section rings}\label{sec:ringextension}

Fix a $G$-invariant element $S\in|-K_X|$ and a nonzero $G$-eigensection $s\in H^0(X,-K_X)$ defining $S$. Define $A\in\Cl(S)$ by
$
 \OO_S(A)=\bigl(\OO_X(-K_X)\otimes_{\OO_X}\OO_S\bigr)^{**}.
$
For a positive integer $r$ such that $-rK_X$ is Cartier, the reflexive sheaf $\OO_S(rA)$ agrees with $\OO_X(-rK_X)|_S$ on $S\cap X_{\reg}$ and hence on $S$. This restriction is an ample line bundle, so $A$ is an ample integral $\QQ$-Cartier Weil divisor class on $S$.

Put
\[
 R(X,-K_X)=\bigoplus_{m\ge0}H^0(X,-mK_X),
 \qquad R(S,A)=\bigoplus_{m\ge0}H^0(S,mA).
\]
For each $m\ge0$, we have $\OO_S(mA)\simeq(\OO_X(-mK_X)\otimes_{\OO_X}\OO_S)^{**}$, since the two reflexive sheaves agree on $S\cap X_{\reg}$ and $S\setminus X_{\reg}$ is finite. For $u\in H^0(X,-mK_X)$, reflexivity of $\OO_S(mA)$ gives a unique extension $\rho(u)\in H^0(S,mA)$ of $u|_{S\cap X_{\reg}}$. The extensions respect multiplication on $S\cap X_{\reg}$, and hence define a graded ring homomorphism $\rho\colon R(X,-K_X)\to R(S,A)$.

\begin{lemma}\label{lem:restrictionring}
The homomorphism $\rho$ induces an isomorphism of graded rings
\begin{equation}\label{eq:restrictionring}
 R(X,-K_X)/(s)\simeq R(S,A).
\end{equation}
In particular, we have $\plur_{-1}(X)=1+h^0(S,A)$.
\end{lemma}
\begin{proof}
Let $m\ge1$. 
Applying \cite[Proposition~5.26]{KM}, we obtain the exact sequence
\[
0\longrightarrow\OO_X(-(m-1)K_X)
 \xrightarrow{\,s\,}\OO_X(-mK_X)
 \xrightarrow{\,\rho_m\,}\OO_S(mA)\longrightarrow0.
\]
For every $m\ge1$, the degree-$m$ component $H^0(\rho_m)$ of $\rho$ is surjective by Kawamata--Viehweg vanishing, and exactness gives $(\ker\rho)_m=(s)_m$. In degree $m=0$, the map $\rho$ is the identity on $\CC$. Thus $\rho$ induces the isomorphism~\eqref{eq:restrictionring}.
\end{proof}

Recall that the inertia group $I_S=\ker(G\to\Aut(S))$ is cyclic by Lemma~\ref{lem:inertiatypes}, so $I_S\simeq\ZZ/d\ZZ$ for some $d\ge1$.

\begin{lemma}\label{lem:cyclicring}
Suppose that there is a graded isomorphism $\psi\colon\CC[Y_1,\ldots,Y_4]/(F)\xrightarrow{\sim}R(S,A)$ sending $Y_j$ to $y_j$, where $\deg Y_j=a_j>0$ and $F$ is homogeneous of degree $d$. Then for some $c\in\CC$ there is a commutative diagram of graded rings
\begin{equation}\label{eq:cyclicextension}
\begin{tikzcd}[column sep=large,row sep=large]
 \CC[Y_1,\ldots,Y_4,z]/(F+c z^d) \arrow[r,"\varphi","\sim"'] \arrow[d,two heads,"z\mapsto0"']
 & R(X,-K_X) \arrow[d,two heads,"\rho"] \\
 \CC[Y_1,\ldots,Y_4]/(F) \arrow[r,"\psi" above,"\sim"']
 & R(S,A).
\end{tikzcd}
\end{equation}
Here $\deg z=1$ and $\varphi(z)=s$. If the $y_j$ are $G$-eigenvectors, then the diagram can be chosen $G$-equivariant for diagonal actions on the polynomial rings.
\end{lemma}
\begin{proof}
Choose a smooth point $Q\in S\cap X_{\reg}$. Let $\nu\colon I_S\to\CC^\times$ be the character of the action on $N^*_{S/X,Q}$. Since $I_S$ fixes $S$ pointwise, $\nu$ is faithful by Lemma~\ref{lem:tangent}. As in the proof of Lemma~\ref{lem:distinctinertia}, each $g\in I_S$ fixes $s$ and acts on $\OO_X(-K_X)|_Q$ by $\nu(g)^{-1}$. For $m\ge0$, the rank-one reflexive sheaf $\OO_S(mA)$ has only scalar automorphisms because $S$ is normal and projective. Thus $g$ acts on $H^0(S,mA)$ by $\nu(g)^{-m}$, the scalar on $\OO_S(mA)|_Q$.

Put $R=R(X,-K_X)$ and $R_m=H^0(X,-mK_X)$. If $0<m<d$, each section in $R_m^{I_S}$ has zero restriction to $S$ because $\nu^{-m}$ is nontrivial. By \eqref{eq:restrictionring}, such a section has the form $sv$ with $v\in R_{m-1}$. Since $s$ is $I_S$-invariant and is a nonzero divisor in $R$, we have $v\in R_{m-1}^{I_S}$. Induction from $R_0=\CC$ gives $R_m^{I_S}=\CC s^m$ for $0\le m<d$.

By semisimplicity, for each $j$ choose a $G$-equivariant linear splitting of the degree-$a_j$ component of $\rho$ and let $\widetilde y_j\in R_{a_j}$ be the image of $y_j$ under this splitting. The group $I_S$ acts on $\widetilde y_j$ by $\nu^{-a_j}$. Since $\deg F=d$, the element $F(\widetilde y_1,\ldots,\widetilde y_4)$ is $I_S$-invariant. Its image under $\rho$ is $F(y_1,\ldots,y_4)=0$, so
\[
 F(\widetilde y_1,\ldots,\widetilde y_4)\in sR_{d-1}^{I_S}=\CC s^d.
\]
Thus $F(\widetilde y_1,\ldots,\widetilde y_4)+cs^d=0$ for some $c\in\CC$. The assignments $z\mapsto s$ and $Y_j\mapsto\widetilde y_j$ define the homomorphism $\varphi$ in \eqref{eq:cyclicextension}.

Via \eqref{eq:restrictionring}, reduction of $\varphi$ modulo $(z)$ in its source and $(s)$ in its target gives $\psi$. Hence $R=\operatorname{im}\varphi+sR$, and induction on degree gives surjectivity. If $u\in\ker\varphi$, then $u=zv$ because $\psi$ is injective. Since $s\varphi(v)=0$, we have $v\in\ker\varphi$. Thus $\ker\varphi=z\ker\varphi$, and the nonnegative grading gives $\ker\varphi=0$. If the $y_j$ are $G$-eigenvectors, then $y_j$ and $\widetilde y_j$ have the same character by the choice of $\widetilde y_j$. Letting $G$ act on $Y_j$ and $z$ by the characters of $\widetilde y_j$ and $s$, respectively, makes \eqref{eq:cyclicextension} $G$-equivariant.
\end{proof}

\subsection{The four Fano threefolds}\label{sec:Fanoreconstruction}

\begin{theorem}\label{thm:Fanomodels}
Let $G\simeq\Gamma_i$ for some $i\in\{1,2,3,4\}$ and suppose that $G$ acts faithfully on a $G\QQ$-Fano threefold $X$ with $\plur_{-1}(X)>0$. Then there exist an isomorphism $\alpha\colon G\to\Gamma_i$ and a $G$-equivariant isomorphism $X\to X_i$, where $G$ acts on $X_i$ via $\alpha$ and the standard diagonal action in Example~\ref{ex:exceptional-actions}.
\end{theorem}
\begin{proof}
Suppose first that $i\in\{1,2\}$. If $i=1$, then $S$ is smooth and $X$ is Gorenstein by \cite[Proposition~2.11]{Lo26}. If $i=2$, then $S$ is smooth by Lemma~\ref{lem:sexticK3}. Since $H_S\simeq H_i$ by Lemma~\ref{lem:inertiatypes}, \eqref{eq:K3lattices} gives $\Pic(S)^{H_S}=\ZZ h$ with $h$ ample. We have $h^2=4$ for $i=1$ and $h^2=2$ for $i=2$. Thus $A\sim kh$ for an integer $k\ge1$. Riemann--Roch and Kodaira vanishing on $S$, together with Lemma~\ref{lem:restrictionring}, give
\[
 \plur_{-1}(X)=1+h^0(S,kh)=3+\frac{k^2h^2}{2}.
\]
The bounds in \eqref{eq:nonemptysectionbounds} give $k=1$. For $i=2$, the ring $R(S,A)\simeq R(S,h)$ is given by \eqref{eq:sexticK3ring}.

For $i=1$, we have $\plur_{-1}(X)=5$. Extend $s$ to a $G$-eigenbasis $s,s_1,\ldots,s_4$ of $H^0(X,-K_X)$. By Lemma~\ref{lem:restrictionring}, the sections $y_j=\rho(s_j)$ form a basis of $H^0(S,A)$ and give homogeneous coordinates on $\PP^3=\PP(H^0(S,A)^\vee)$. Since $|-K_X|$ is base point free \cite[Proposition~3.4]{Lo26} and $\rho$ is surjective in degree $1$, the linear system $|A|$ defines a morphism $f_A\colon S\to\PP^3$. Since $A$ is ample, $f_A$ is finite onto $Z=f_A(S)$. The surface $Z$ is nondegenerate and $(\deg f_A)\deg Z=A^2=4$, so $\deg f_A\in\{1,2\}$.

Let $F_1\in\CC[y_1,\ldots,y_4]$ be an irreducible homogeneous polynomial defining $Z$. For $j=2,3,4$, let $D_j$ be the diagonal matrix that multiplies $y_j$ by $\xi_4$ and fixes the other coordinates. By Lemma~\ref{lem:distinctinertia}, the class of $D_j$ belongs to the image of $G$ in $\Aut(\PP^3)$. Since $Z$ is $G$-invariant, we have $D_j^*F_1\in\CC^\times F_1$. Thus the exponents of $y_j$ in the monomials of $F_1$ with nonzero coefficient are congruent modulo $4$. Since $F_1$ is irreducible of degree at least $2$, we have $y_j\nmid F_1$, so the common residue is $0$.

If $\deg f_A=2$, then $\deg F_1=2$ and $F_1\in\CC^\times y_1^2$, contrary to the irreducibility of $F_1$. Hence $\deg f_A=1$ and $F_1=\sum_{j=1}^4 c_j y_j^4$. The quartic $Z$ is normal by \cite[\S6.5.2]{SaintDonat}, so the finite birational morphism $f_A\colon S\to Z$ is an isomorphism. If $c_j=0$, then the $j$th coordinate point lies on $Z$ and all partial derivatives of $F_1$ vanish there, contradicting smoothness of $Z$. Choose $r_j\in\CC^\times$ with $r_j^4=c_j$. In the coordinates $y'_j=r_j y_j$, the surface $Z$ is the Fermat quartic $S_1=\{\sum_{j=1}^4(y'_j)^4=0\}\subset\PP^3$. The homogeneous coordinate ring of $S_1$ is $R(S,A)$ by \cite[Theorem~6.1(ii)]{SaintDonat}.

Now suppose that $i\in\{3,4\}$. By Lemma~\ref{lem:inertiatypes}, we have $H_S=G/I_S\simeq H_i$. Lemmas~\ref{lem:restrictionring} and \ref{lem:distinctinertia} give
\begin{equation}\label{eq:smallrestriction}
 h^0(S,A)=\plur_{-1}(X)-1\le2.
\end{equation}
By Proposition~\ref{prop:markedK3}, we can identify $H_S$ with $H_i$ so that the minimal resolution of $S$ is $H_i$-equivariantly isomorphic to $T_i$. If $S$ were smooth, then Riemann--Roch and Kodaira vanishing for the ample line bundle $\OO_S(A)$ would give $h^0(S,A)=2+A^2/2\ge3$, contrary to \eqref{eq:smallrestriction}. Hence Lemma~\ref{lem:K3contractions} gives an $H_i$-equivariant isomorphism $S\simeq S_i$. Under this isomorphism, $A$ is an ample $H_i$-invariant class. By \eqref{eq:invariantWeilclasses}, we have $A=kA_0$ for some $k\in\ZZ_{>0}$, where $A_0=[\OO_{S_i}(1)]$. If $k\ge2$, then $x^k,x^{k-1}y,y^k$ are linearly independent sections of $\OO_{S_i}(k)$ because $x/y$ is nonconstant on $S_i$. This contradicts \eqref{eq:smallrestriction}, so $k=1$.

Let $S_2\subset\PP(1,1,1,3)$ be the smooth sextic defined by \eqref{eq:sexticK3ring}. Table~\ref{tab:reconstruction} gives $(S,A)$, the order $d$ of $I_S$, and the degrees $a_j$ of the homogeneous coordinates $y_j$ of $S_i$. Using the identifications in the table, these coordinates generate $R(S,A)$ with the single relation $\sum_{j=1}^4 y_j^{d/a_j}=0$.

\begin{table}[H]
\centering
\renewcommand{\arraystretch}{1.15}
\begin{tabular}{@{}clcc@{}}
\toprule
$G$ & $(S,A)$ & $d=|I_S|$ & $(a_1,a_2,a_3,a_4)$\\
\midrule
$\Gamma_1$ & $(S_1,[\OO_{S_1}(1)])$ & $4$ & $(1,1,1,1)$\\
$\Gamma_2$ & $(S_2,[\OO_{S_2}(1)])$ & $6$ & $(1,1,1,3)$\\
$\Gamma_3$ & $(S_3,[\OO_{S_3}(1)])$ & $6$ & $(1,1,2,2)$\\
$\Gamma_4$ & $(S_4,[\OO_{S_4}(1)])$ & $8$ & $(1,1,2,4)$\\
\bottomrule
\end{tabular}
\par\vspace{2pt}
\caption{Reconstruction of the Fano threefolds $X_1,\ldots,X_4$ from anticanonical K3 surfaces}\label{tab:reconstruction}
\end{table}

The action of $G$ on $S$ is diagonal in the coordinates $y_1,\ldots,y_4$ and therefore preserves the Weil divisor of zeros of each section $y_j\in H^0(S,a_jA)$. For $g\in G$, the quotient $g^*y_j/y_j$ is a rational function on $S$ with zero divisor. Since $S$ is normal and projective, $g^*y_j/y_j\in\CC^\times$. Thus each $y_j$ is a $G$-eigenvector.

Apply Lemma~\ref{lem:cyclicring} with $F=\sum_{j=1}^4Y_j^{d/a_j}$. Taking $\Proj$ of the isomorphism $\varphi$ in \eqref{eq:cyclicextension} identifies $X$ with the hypersurface $F+cz^d=0$ in $\PP(1,a_1,\ldots,a_4)$, with a diagonal $G$-action. If $c=0$, the point $[1\colon0\colon\cdots\colon0]$ on this hypersurface is fixed by $G$, contrary to Proposition~\ref{prop:invariants}.

On the affine chart $z\ne0$ of $\PP(1,a_1,\ldots,a_4)$, put $u_j=Y_j/z^{a_j}$ and $F_0=c+\sum_{j=1}^4u_j^{d/a_j}$. Let $\chi_s$ and $\chi_j$ be the characters of $s$ and $\widetilde y_j=\varphi(Y_j)$, respectively. For $g\in G$, we have $g^*u_j=\beta_j(g)u_j$, where $\beta_j=\chi_j\chi_s^{-a_j}$. Since $g$ preserves $X\cap\{z\ne0\}$, the pullback $g^*$ preserves its ideal $(F_0)$ in $\CC[u_1,\ldots,u_4]$. Thus $g^*F_0$ is a scalar multiple of $F_0$. Comparing the nonzero constant terms gives $g^*F_0=F_0$, and the coefficients of $u_j^{d/a_j}$ give $\beta_j(g)^{d/a_j}=1$ for all $j$. Write $\beta_j(g)=\xi_{d/a_j}^{b_j(g)}$ with $b_j(g)\in\ZZ/(d/a_j)\ZZ$. The homomorphism $\alpha\colon G\to\Gamma_i$ given by $g\mapsto(b_1(g),\ldots,b_4(g))$ is injective because its kernel fixes the dense open subset $X\cap\{z\ne0\}$ pointwise and $G$ acts faithfully on $X$. Since $|G|=|\Gamma_i|$, the homomorphism $\alpha$ is an isomorphism. For $\zeta\in\CC^\times$ with $\zeta^d=c$, the change of coordinate $z'=\zeta z$ identifies $X$ with $X_i$. On $z'\ne0$, the coordinates $Y_j/(z')^{a_j}=\zeta^{-a_j}u_j$ have characters $\beta_j$, so the isomorphism $X\to X_i$ is $G$-equivariant for the action via $\alpha$.
\end{proof}

\subsection{Proofs of the main results}\label{sec:mainproofs}

\begin{proof}[Proof of Theorem~\ref{thm:rc-models}]
By Lemma~\ref{lem:complexitytypes}, none of $\Gamma_1,\ldots,\Gamma_4$ is of product type, so the two alternatives are mutually exclusive. 
Suppose that $G$ is not of product type. By Theorem~\ref{thm:nonemptymodel}, there is a $G\QQ$-Fano threefold $Y$ with nonempty anticanonical linear system and a $G$-equivariant birational map $X\dashrightarrow Y$. Theorem~\ref{thm:nonempty} gives $G\simeq\Gamma_i$ for some $i\in\{1,2,3,4\}$. The groups $\Gamma_1,\ldots,\Gamma_4$ are pairwise nonisomorphic, so $i$ is unique. By Theorem~\ref{thm:Fanomodels}, there is a $G$-equivariant isomorphism $Y\to X_i$, where $G$ acts on $X_i$ via an isomorphism $G\to\Gamma_i$. Composing with $X\dashrightarrow Y$ gives the required birational map. By Example~\ref{ex:exceptional-actions}, the hypersurface $X_i$ is nonrational, so $X$ is nonrational.
\end{proof}

\begin{proof}[Proof of Theorem~\ref{thm:cremona}]
Let $G\subset\Cr_3(\CC)$ be finite abelian. The given birational action of $G$ on $\PP^3$ regularizes to a faithful action on a smooth projective threefold $X$ birational to $\PP^3$. Theorem~\ref{thm:rc-models} implies that $G$ is of product type. Conversely, suppose that $G\simeq G_1\times G_2$ is of product type. Choose a smooth projective rational surface $S$ with a faithful $G_2$-action. The actions of $G_1$ on $\PP^1$ and of $G_2$ on $S$ give a faithful action of $G$ on the rational threefold $\PP^1\times S$. Thus $G$ embeds into $\Cr_3(\CC)$.
The list \eqref{eq:productlist} is established in the proof of Lemma~\ref{lem:complexitytypes}.
\end{proof}

\begin{proof}[Proof of Theorem~\ref{thm:complexity}]
For $n=1,2$, rationally connected varieties are rational, so the assertion follows from Lemma~\ref{lem:complexitytypes}\ref{item:complexity-low-dimension}. For $n=3$, it follows from Theorem~\ref{thm:rc-models} and Lemma~\ref{lem:complexitytypes}\ref{item:complexity-three}. This completes the proof.
\end{proof}

\end{document}